%% file: 2005-46.tex
\documentclass{gtart_h}  

\input gtspec

\lognumber{502}
\received{3 October 2004}
\volumenumber{9}\papernumber{46}\volumeyear{2005}
\pagenumbers{2013}{2078}   
\revised{24 July 2005}
\published{26 October 2005}
\accepted{17 September 2005}
\proposed{Yasha Eliashberg}
\seconded{Peter Ozsv\'ath, Tomasz Mrowka}

\def\psfraga <#1,#2> #3#4{%
\psfrag {#3}{\smash{\rlap{\kern #1 \raise #2\hbox{#4}}}}}

\def\fref#1{\hyperlink{#1anchor}{\ref*{#1}}}
\def\figref#1{\hyperlink{#1anchor}{Figure~\ref*{#1}}}
\def\anchor#1{\noindent\hypertarget{#1anchor}{\smash{$\phantom{99}$}}\newline}

\usepackage{amssymb}
\usepackage{amsthm}
\usepackage{amsmath}
\usepackage[mathscr]{eucal}
\usepackage{graphicx}
\usepackage{psfrag}
\usepackage{amscd}
\usepackage{enumerate}

\theoremstyle{plain}
\newtheorem{tetel}{Theorem}[section]
\newtheorem{all}[tetel]{Proposition}
\newtheorem{lemma}[tetel]{Lemma}
\newtheorem{kov}[tetel]{Corollary}
\theoremstyle{definition}\newtheorem{Def}[tetel]{Definition}
\theoremstyle{remark}\newtheorem{megj}[tetel]{Remark}
\newtheorem{pelda}[tetel]{Example}

\newcommand*{\R}{\ensuremath{\mathbf R}}

\newcommand*{\Z}{\ensuremath{\mathbf Z}}

\newcommand*{\di}{\mathrm d}

\DeclareMathOperator{\Emb}{Emb}
\DeclareMathOperator{\id}{id}
\DeclareMathOperator{\im}{im}

\DeclareMathOperator{\Aut}{Aut}
\DeclareMathOperator{\Leg}{Leg}

\usepackage{array}
\usepackage{tabularx}

\newcommand*{\leg}{\Leg(S^1,\R^3)}
\newcommand*{\Zt}{\Z[t,t^{-1}]}

\DeclareMathOperator{\Crit}{Crit}

\begin{document}

\title{Contact homology and one parameter families\\of Legendrian
knots}
\authors{Tam\'{a}s K\'{a}lm\'{a}n}
\coverauthors{Tam\noexpand\'{a}s K\noexpand\'{a}lm\noexpand\'{a}n}
\asciiauthors{Tamas Kalman}

\address{Department of Mathematics,
University of Southern California\\Los Angeles, CA 90089, USA}
\email{tkalman@usc.edu}

\begin{abstract}
We consider $S^1$--families of Legendrian knots in the standard
contact $\R^3$. We define the monodromy of such a loop, which is
an automorphism of the Che\-ka\-nov--Eli\-ash\-berg contact homology of
the starting (and ending) point. We prove this monodromy is a
homotopy invariant of the loop (Theorem \ref{thm:hurkok}). We also
establish techniques to address the issue of Reidemeister moves of
Lagrangian projections of Legendrian links. As an application, we
exhibit a loop of right-handed Legendrian torus knots which is
non-contractible in the space $\Leg(S^1,\R^3)$ of Le\-gend\-rian
knots, although it is contractible in the space $\Emb(S^1,\R^3)$
of smooth knots. For this result, we also compute the contact
homology of what we call the Legendrian closure of a positive
braid (Definition \ref{def:lezaras}) and construct an augmentation
for each such link diagram.
\end{abstract}

\asciiabstract{%
We consider S^1-families of Legendrian knots in the standard contact
R^3.  We define the monodromy of such a loop, which is an automorphism
of the Chekanov-Eliashberg contact homology of the starting (and
ending) point.  We prove this monodromy is a homotopy invariant of
the loop.  We also establish techniques to address the issue of
Reidemeister moves of Lagrangian projections of Legendrian links.  As
an application, we exhibit a loop of right-handed Legendrian torus
knots which is non-contractible in the space Leg(S^1,R^3) of
Legendrian knots, although it is contractible in the space
Emb(S^1,R^3) of smooth knots.  For this result, we also compute the
contact homology of what we call the Legendrian closure of a positive
braid and construct an augmentation for
each such link diagram.}



\maketitle

\section{Introduction}

A Legendrian knot is an embedding $S^1\hookrightarrow\R^3_{xyz}$
which is everywhere tangent to the kernel of the standard contact
$1$--form $\alpha=\di z - y\di x$. Two such knots are equivalent
if they are homotopic through Legendrian knots. This is a
refinement of the obvious subdivision according to classical knot
type; just how much finer is a topic of current research. (We may
refer to this as the $\pi_0$--problem of the space
$\Leg(S^1,\R^3)$ of Legendrian knots). The integer-valued
invariants $tb$ (Thurston--Bennequin number) and $r$ (rotation
number) guarantee that every classical knot type contains
infinitely many Legendrian equivalence classes. It was a major
break-through when Chekanov \cite{chek} introduced a much more
complicated, Floer-theoretical invariant which was able to
distinguish between Legendrian knots of the same $tb$, $r$, and
classical type (since then, Etnyre and Honda \cite{EH1,EH2,EH3}
have been able to prove similar results with different methods).
The form in which we'll consider Chekanov's invariant is that of a
system of homology groups, which we will call contact homology.

In this paper we use contact homology to tackle another, so called
$\pi_1$--problem; namely, we use monodromy induced on this
homology to provide the first examples that Legendrian and
classical knots are different at the level of $1$--parameter
families.

\subsection{Statements of results}

In its most general formulation, Chekanov's invariant of
Legendrian knots \cite{chek} takes the form of a stable tame
isomorphism class of differential graded algebras, or DGA's for
short, $(\mathscr A_*,\partial)$. What's more important for us is
that in particular, the homology $\ker\partial/\im\partial$ is an
invariant of the Legendrian knot type. This homology, which is
actually also a graded algebra over $\Z_2$, will be referred to as
contact homology and will be denoted by $H(L)$, where $L$ is a
(representative of a) Legendrian knot type. In \cite{ENS}, Etnyre,
Ng, and Sabloff defined a natural extension of Chekanov's work,
namely a DGA and corresponding contact homology with
$\Zt$--coefficients. We will carry out our proofs in this more
general setting, even though in all of the applications known to
the author the original $\Z_2$--version (which will be used in
sections \ref{sec:BC}--\ref{sec:moreaug}) would be sufficient.

The main result of the paper is the following:

\begin{tetel}\label{thm:hurkok}
For every Legendrian knot type $\mathscr L\subset\Leg(S^1,\R^3)$ and
generic $L\in\mathscr L$, there exists a multiplicative homomorphism
\[\mu\colon\pi_1(\mathscr L,L)\to\Aut(H(L)),\]
defined by continuation on the Chekanov--Eliashberg contact homology
$H$. This is true with either $\Z_2$ or $\Z[t,t^{-1}]$--coefficients.
\end{tetel}

\begin{figure}[ht!]\small\anchor{fig:hurok}
\cl{\includegraphics[width=.97\hsize]{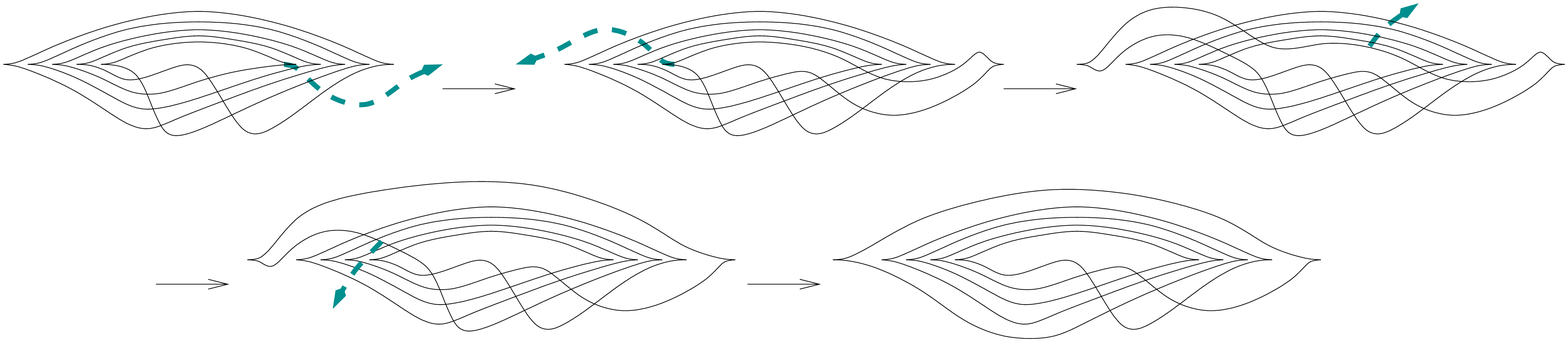}}
\caption{The loop $\Omega_{p,q}$ when $p=2$ and $q=5$, shown in the
front projection}\label{fig:hurok}
\end{figure}

In particular, monodromy calculations are carried out in the
homology and neither in its linearized \cite{chek} or abelianized
\cite{ENS} versions, nor in the DGA itself. Augmentations (which
are the objects also needed to define linearized homology) do
however play an important role in proving that $\mu$ is
non-trivial. This will be demonstrated by the following example.
Let $\mathscr L$ be the space of positive Legendrian $(p,q)$ torus
knots with maximal Thurston--Bennequin number (where $p,q\ge2$ are
relatively prime). $\mathscr L$ is connected by a result of Etnyre
and Honda \cite{EH1}. Consider the loop
$\Omega_{p,q}\subset\mathscr L$ of \figref{fig:hurok} (a more
elaborate definition is in section \ref{sec:loop}). We claim the
following:

\begin{tetel}\label{thm:trefli}
$[\Omega_{p,q}]\in\pi_1(\mathscr L)$ is either an element of infinite
order or its order is divisible by $p+q$.
\end{tetel}

By Theorem \ref{thm:hurkok}, this follows immediately from the 
following:

\begin{tetel}\label{thm:monorendje}
The restricted monodromy
$\mu_0(\Omega_{p,q})=\mu(\Omega_{p,q})\big|_{H_0(L)}$ of the loop
$\Omega_{p,q}$, where $H_0(L)$ is the index $0$ part of the
$\Z_2$--coefficient contact homology $H(L)$, has order $p+q$.
\end{tetel}

In particular, $\mu(\Omega_{p,q})$ is different from the identity
automorphism, and because $p+q$ doesn't divide $2p$, so is
$(\mu(\Omega_{p,q}))^{2p}$. Thus by Theorem \ref{thm:hurkok},
$[\Omega_{p,q}]^{2p}\in\pi_1(\mathscr L)$ is non-trivial. On the other
hand, if $\mathscr K$ is the space of smooth positive $(p,q)$ torus
knots, then $[\Omega_{p,q}]^{2p}\in\pi_1(\mathscr K)$ \emph{is}
trivial\footnote{In fact, an elementary application of Hatcher's
work \cite{hatch} shows that if $p+q$ is odd, then $\Omega_{p,q}$ is
itself contractible in $\mathscr K$, and if $p+q$ is even, then 
$[\Omega_{p,q}]\in\pi_1(\mathscr K)$ has order $2$.}, because 
$\Omega_{p,q}$ is homotopic to a rotation by 
$2\pi/p$ radians and $\pi_1(SO(3))=\Z_2$. So Theorem \ref{thm:trefli}
had to be demonstrated via a specifically ``Legendrian'' technique,
which is the monodromy invariant $\mu$. We also have

\begin{kov}\label{kov}
There exist Legendrian knot types $\mathscr L$ so that for the
corresponding smooth knot type $\mathscr K\supset\mathscr L$, the
homomorphism $\pi_1(\mathscr L)\to\pi_1(\mathscr K)$ induced by the
inclusion is not injective.
\end{kov}

It is an open problem whether this homomorphism, for fixed $\mathscr
K$ \emph{and} $\mathscr L$, is always surjective. It is also
worthwhile to compare our results to those of Benham, Lin, and
Miller in \cite{lin}. They prove that if $\mathscr W$ is the set of
knots in $\mathscr K$ with writhe $tb(\mathscr L)$ (so that
$\mathscr L\subset\mathscr W\subset\mathscr K$) then the inclusion
$\mathscr W\hookrightarrow\mathscr K$ is a weak homotopy
equivalence. For results about the homotopy type of $\mathscr K$,
see Hatcher \cite{hatch}.

There are two widely used ways to represent Legendrian knots in
the standard contact $\R^3$. In \figref{fig:hurok}, we used
the front projection. For most of the paper, however, we will use
Lagrangian projections, for the simple reason that $H$ is
naturally defined in terms of those. (Ng \cite{computable} gives a
definition in terms of fronts.) Lagrangian projections are also
more symmetric, but on the other hand less robust than fronts. By
this I mean that certain linear constraints need to be satisfied
both for the existence of an actual Legendrian lift of the diagram
and for doing Reidemeister moves as well. While \cite[section
11]{chek} contains a thorough description of the former
phenomenon, our Theorem \ref{thm:linprog} appears to be the first
account in the literature on the latter problem. We formulate a
solution in terms of techniques and language of linear
programming.

The paper is organized as follows. In section \ref{sec:floer} we
give some motivation and describe related work. We review some basic
material in section \ref{sec:prep}. Section \ref{sec:ind} contains
the definition of monodromy and the proof of Theorem
\ref{thm:hurkok}. Section \ref{sec:simplex} concludes the general
discussion of one parameter families of Legendrian links by addressing
Reidemeister moves.

Then we turn to examples; section \ref{sec:trefli} contains the
simplest non-trivial one known to the author, a loop of Legendrian
trefoils. In section \ref{sec:loop} we include it in a
much more general class of examples. Those loops consist of so
called Legendrian closures of positive braids. Section \ref{sec:BC}
examines their base points; in particular, we compute their DGA. In
section \ref{sec:aug} we construct an augmentation for that DGA.
Finally, we specialize our discussion again to the case when the 
Legendrian closure
is a positive torus knot and use the last two sections to prove
Theorem \ref{thm:monorendje}.

\medskip

{\bf Acknowledgments}\qua This paper is almost identical to the
author's PhD thesis that was submitted in 2004 at the
University of California at Berkeley. It is my pleasure to thank
the institution, and in particular my advisors, Rob Kirby and
Michael Hutchings, for all the generous support that they
provided. It was in fact Mike Hutchings who gave me the idea of
this project. I would also like to thank Lenny Ng and Josh Sabloff
for checking earlier versions of the manuscript, as well as them,
Tobias Ekholm, Michael Sullivan and Andr\'as Sz\H ucs for
stimulating conversations.  I was supported by NSF grant
DMS-0244558 during part of the research.

\subsection{Floer theory}\label{sec:floer}

Loosely speaking, the idea of Floer theory is to formulate a
Morse-type homology theory of certain real-valued smooth functionals
on certain in\-fi\-nite-di\-men\-sion\-al manifolds. The manifold
and the functional are typically associated to an object such as a
pair of Lagrangians, a symplectomorphism, a contact structure, a
Legendrian submanifold etc. Some special cases of this project have
had great success (alas, analytical details are sometimes missing)
and we refer to these as ``Floer theories.'' We should mention at
least one case when the analogy with finite-dimensional Morse theory
starts breaking down: in symplectic field theory \cite{field}, of
which contact homology is a special case \cite{ENS}, generators of
Floer homology are indeed critical points of an action functional,
but flow lines (which are pseudoholomorphic curves in an appropriate
sense) starting from a critical point may split and end at a finite
collection of critical points (as opposed to just one endpoint). The
natural way to include this phenomenon in our algebraic formulation
of the theory is to write the multiple endpoints as a product. This
is how the chain complex in this case takes the form of a
differential graded \emph{algebra} (DGA).

Besides finding new Floer theories, it is also an important and
promising direction of research to generalize more and more features
of finite-dimensional Morse theory to these new settings. One such
feature can be the so called \emph{continuation map}. Suppose
$(f_t,g_t)$, $t\in[0,1]$ is a $1$--parameter family of smooth
functions and Riemannian metrics on the finite-dimensional smooth
manifold $F$. Suppose that the family is generic. This in particular
means that the pairs $(f_t,g_t)$ are Morse--Smale (ie, they have
associated Morse homologies $H^{(t)}$) with finitely many exceptions
$0<t_1<\ldots<t_k<1$. Note that we can \emph{not} speak of homologies
$H^{(t_i)}$. Yet we may define a vector field on $[0,1]\times F$ by
$V=t(1-t)\partial_t + V_t$, where $V_t$ is the negative gradient of
$f_t$ with respect to $g_t$, and after observing that
$\Crit_i(V)=\left(\{\,0\,\}\times\Crit_{i-1}(V_0)\right)\cup
\left(\{\,1\,\}\times\Crit_i(V_1)\right)$, we may count flow lines of
$V$ to define a chain map from the Morse complex of $(f_0,g_0)$ to
that of $(f_1,g_1)$. The fact that it induces an isomorphism (which we
call the continuation map) $H^{(0)}\to H^{(1)}$ can be established by
a certain construction of chain homotopies \cite{floer}, similar in
flavor to the above\footnote{What has just been outlined is what I
call the `analytic approach.' Alternatively, one may study the
\emph{bifurcations} that take place at $t=t_1,\ldots,t_k$ and define
`combinatorial' continuation maps $H^{(t_i-\varepsilon)}\to
H^{(t_i+\varepsilon)}$, as well as chain homotopies to guarantee that
these are isomorphisms. The approach taken in this paper will in fact
be a generalization of the latter.}. This provides an alternative
proof that the isomorphism class of Morse homology depends only on the
manifold, thus making Morse theory a self-contained homology theory.


But the usefulness of the continuation map doesn't stop here.
Consider a fibration $\pi\colon Z\to S^1$ with smooth manifold
fiber $F$, along with a smooth function $f\colon Z\to\R$ and a
Riemannian metric $g$ on $Z$. Suppose that for $1\in S^1$, the
pair $(f\big|_{\pi^{-1}(1)},g\big|_{\pi^{-1}(1)})$ is Morse--Smale
and the family $(f\big|_{\pi^{-1}(t)},g\big|_{\pi^{-1}(t)})$ is
generic. Then by continuation, we obtain an automorphism (called
\emph{monodromy}) $H^{(1)}\to H^{(1)}$. This is naturally
equivalent to the automorphism of the singular homology of
$\pi^{-1}(1)$ induced by the gluing map of the fibration (which is
well defined up to homotopy, therefore its action is well
defined). In particular, a non-trivial monodromy implies the
non-triviality of the fibration.

As already demonstrated by Seidel \cite{Seidel} and Bourgeois
\cite{burzsuj}, the continuation (monodromy) map can prove
non-triviality theorems, in a way analogous to the above, about
the fundamental group of the space of objects (``object'' means a
symplectomorphism in \cite{Seidel} and a contact structure in
\cite{burzsuj}), to which a version of Floer theory is
associated\footnote{\label{lab}Strictly speaking, Floer homology
is associated only to generic objects. Once we have proven that
its isomorphism class is invariant within a path-component of the
objects, this distinction becomes immaterial, provided that we
only wish to use Floer theory to separate such path-components
(ie, we do $\pi_0$--theory). But the presence of non-generic
objects is potentially an issue again when we investigate the
topology of the path-components (like their $\pi_1$, as in this
paper).}. To get such conclusions, it is of course a key point to
prove that monodromy is \emph{invariant} under homotopies of the
loop of objects. In this paper we obtain such a result about the
space of Legendrian knots $\leg$.

There is also a generalized invariant of families of objects
parametrized by an arbitrary finite-dimensional manifold, recently
discovered by Hutchings \cite{Mike}. The value of this is a spectral
sequence. \cite{Mike} contains its description in the case of
finite-dimensional Morse theory, where the spectral sequence
associated to a generic $4$--tuple $(Z,B,f,g)$ (as above, except that
$S^1$ is replaced by the arbitrary finite-dimensional manifold $B$)
turns out to be isomorphic to the Leray--Serre spectral sequence of
the fibration. A paper \cite{Mike2} detailing the case of the Floer
theory of symplectomorphisms, thus generalizing Seidel's results, is
in preparation. In other Floer theories, apart from Bourgeois's
results (where $B=S^n$) and the present paper (where $B=S^1$), the
spectral sequence invariant is yet to be either established (see the
next paragraph for more on this when $B=S^1$) or exploited.

The motivation behind the group of results concerning one-parameter
families is that Floer homology should not just be viewed as an
invariant of certain objects, but should also be treated as a twisted
coefficient system over the space of those objects; its non-trivial
twisting is then capable of detecting the fundamental group of the
moduli space. The existence of this twisted coefficient system, modulo
an issue of fixing the signs in the continuation maps (see \cite{Mike}
for more), is an immediate consequence of what's already in the
literature for many versions of Floer theory. These are the versions
where they prove topological invariance by constructing continuation
maps and chain homotopies using the analytic methods introduced in
Floer's seminal paper \cite{floer} (as opposed to a study of
bifurcations, which also appears in earlier work of Floer). In these
cases the existence of a monodromy invariant, which is equivalent to
the spectral sequence when $B=S^1$, is automatic too.

However, the invariance of contact homology was established in a
different way, namely through a certain combinatorial ``bifurcation
analysis'' \cite{chek}. Therefore we can use the analytic
considerations mentioned so far only as a source of motivation that
the monodromy invariant (and the spectral sequence invariant) exists
in this context, even though one suspects that once monodromy is
defined analytically, it will agree with our combinatorially defined
version\footnote{It \emph{is} proven in \cite{ENS} that Chekanov's
combinatorially defined homology is isomorphic to the relative contact
homology whose analytical definition was outlined in \cite{field} and
completed in \cite{EES}. For proof of invariance, however, the authors
use bifurcation analysis in \cite{EES}.}. Following this ``hint''
then, we streamline Chekanov's original proof in section \ref{sec:cm}
(without changing its combinatorial nature) to emphasize the role of
continuation more, and then in section \ref{sec:ch}, we continue the
combinatorial study of bifurcations to provide the chain homotopies
needed to prove that monodromy is a well-defined invariant.

\section{Preliminaries}\label{sec:prep}

For a more thorough treatment and proofs of what follows see 
\cite{chek,ENS,knotsurv}.

Throughout this paper, we work in the standard contact $\R^3$; we 
write the $1$--form, whose kernel field is the contact distribution, 
as $\alpha=\di z - y\di x$. A curve $\lambda\colon S^1\to\R^3$ is 
called a \emph{Legendrian knot} if it is embedded (in particular, it 
has an everywhere non-zero tangent) and $\lambda^*(\alpha)=0$. The 
space of such knots, modulo reparametrization, is denoted by $\leg$. 
We will refer to connected components of $\leg$ (with respect to the 
quotient of the $C^\infty$ topology) as \emph{Legendrian knot 
types}\footnote{We will use $\lambda$ for actual maps and $L$ for 
their equivalence classes modulo reparametrization, while $\mathscr L$ 
will be used to denote Legendrian knot types.}. When talking about 
\emph{classical} or \emph{smooth knots}, we simply mean that the 
Legendrian assumption is not made. Modulo reparametrization, these 
form the moduli space $\Emb(S^1,\R^3)$, whose connected components are 
called \emph{classical} or \emph{smooth knot types}.

A (Legendrian or classical) \emph{link} is a finite disjoint union of 
knots.

\subsection{Link projections}

The \emph{Lagrangian projection} is $\pi\colon(x,y,z)\mapsto(x,y)$. If
$L$ is a Legendrian link, it admits an $xy$--projection that is an
immersion with the additional properties described below in
Proposition \ref{teruletes}. We call $L$ and $\gamma=\pi(L)$
\emph{generic} if $\gamma$ satisfies the usual assumptions for
classical link diagrams, ie,\ it has no singularities other than
finitely many transverse self-crossings.

\begin{figure}[ht!]\small\anchor{fig:frontrei}
\psfrag{11}{I$^{-1}$}
\psfrag{1}{I}
\psfrag{2}{II}
\psfrag{22}{II$^{-1}$}
\psfrag{3}{III}
\cl{\includegraphics[width=.7\hsize]{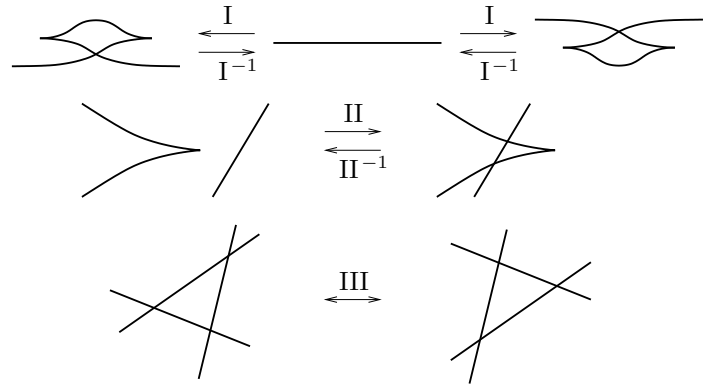}}
\caption{Reidemeister moves in the front projection; the list is
completed by three more moves, which are just the reflections of the
middle move in the vertical, horizontal, and in both,
axes.}\label{fig:frontrei}
\end{figure}

The \emph{front projection} is $(x,y,z)\mapsto(x,z)$. The image of
a generic Legendrian link under this is a smooth curve with
finitely many semicubical cusp points, no self-tangencies and no
vertical tangents. As $y=\di z/\di x$, at every crossing the
strand with smaller slope is the overcrossing one. Generic
homotopies are composed of isotopies and a list of Reidemeister
moves (see \figref{fig:frontrei}). Unlike for Lagrangian
projections, these moves can always be carried out without further
restrictions. For example, the definition of $\Omega_{p,q}$ in
\figref{fig:hurok} (with $2(q-1)$ Reidemeister II moves and
$2(q-1)$ Reidemeister II$^{-1}$ moves) is sound. For the rest of
the paper, we won't need to use fronts, except occasionally for
motivation.

Fix an orientation of the Legendrian knot $L$. Let $a$ be a simple
transverse self-intersection (crossing) in its Lagrangian projection
$\gamma=\pi(L)$ and denote the difference in $z$--coordinate between
the two preimages of $a$ by $h(a)$. We will refer to the positive
number $h(a)$ as the \emph{height} of $a$. Starting at the
undercrossing, follow $\gamma$ until we reach $a$ again, this time on
the upper strand. This path is called the \emph{capping path} of $a$
and we denote it by $\gamma_a$.

\begin{figure}[ht!]\small\anchor{fig:Reeb}
\cl{\includegraphics[width=.5\hsize]{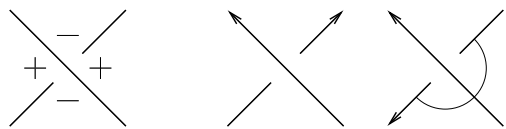}}
\caption{Reeb signs and orientation signs}\label{fig:Reeb}
\end{figure}

We use capping paths to define the \emph{grading}, or \emph{index} of
a crossing, as follows. Assume that all crossings happen at a right
angle. Then the rotation $r(\gamma_a)$ of the tangent to the capping
path $\gamma_a$ (with respect to the orientation $\di x\wedge\di y$ of
the plane) is an odd multiple of $1/4$. The grading of the crossing
$a$ is the integer
\begin{equation}\label{eq:grad}|a|=-2r(\gamma_a)-\frac12.\end{equation}
The \emph{rotation number} $r$ of a Legendrian knot is the Whitney index 
of its
(oriented) Lagrangian projection. Its sign does depend on the choice of 
orientation. (However the modulo $2r$ residue of the grading of a
crossing does not depend on the orientation. This residue
is the value of Chekanov's original grading.) 
The \emph{Thurston--Bennequin number} $tb$ of a Legendrian knot is
the writhe of its Lagrangian projection, ie,\ the number of
crossings with even grading minus the number of those with odd
grading\footnote{\label{lab'}It is easy to verify that $(-1)^{|a|}$
coincides with the sign of the crossing $a$ in classical knot
theory. This fact is assumed in \figref{fig:Reeb}, too.}. It doesn't 
depend on the orientation of the knot. 
The parities of $r$
and $tb$ are opposite. The list of the three \emph{classical
invariants} of a Legendrian knot is completed by its \emph{classical
knot type}.

The \emph{Reeb signs} of the four quadrants surrounding a crossing are
defined according to the scheme on the left of \figref{fig:Reeb}.
For the definition of the DGA with $\Zt$ coefficients in the next
subsection, we will also need the concept of the so called
\emph{orientation sign}. This is $+1$ for all quadrants surrounding an
oddly graded crossing and at evenly graded crossings, the two marked
quadrants of \figref{fig:Reeb} have orientation sign $-1$, the
other two have $+1$. The only property of the orientation
signs that matters for the purposes of defining a homology theory is
the following:

\begin{lemma}{\rm\cite{ENS}}\label{elojellemma}\qua
Around the crossing $a$, the product of the orientation signs of any
two opposite quadrants is $-(-1)^{|a|}$.
\end{lemma}

Let $a_1,a_2,\ldots,a_n$ be the crossings of a \emph{classical}
link diagram $\gamma$. Let us denote the bounded components of the
complement by $U_1,U_2,\ldots,U_{n+1}$. Let
$S_\gamma=\bigoplus_{j=1}^n\R a_j$ and
$V_\gamma=\bigoplus_{i=1}^{n+1}\R U_i$ be vector spaces with
labeled bases and let $S_\gamma^+$ and $V_\gamma^+$ be their
respective (open) positive cones. Chekanov defined the linear map
$\Psi_\gamma\colon S_\gamma\to V_\gamma$ by the matrix $[E_{ij}]$,
where $E_{ij}\in\{-2,-1,0,1,2\}$ is the number of quadrants at
$a_j$ that face $U_i$, counted with Reeb signs.

\begin{all}\label{teruletes}
Suppose $\gamma$ is the Lagrangian projection of a generic Legendrian
link $L$. For the areas $f_i$ of $U_i$, heights $h_j$ of $a_j$, and
incidence coefficients $E_{ij}$ as described above, we have
\begin{equation}\label{eq:teri}f_i=\sum_{j=1}^n E_{ij}h_j.\end{equation}
\end{all}

\begin{proof}
First, it is an elementary fact that any arc
$g=(x,y)\colon[0,1]\to\R^2_{xy}$ has Legendrian lifts
$\lambda(t)=(g(t),z(t))$ in $\R^3$: let $z(0)=0$ and define
$z(t)=\int_0^t y(t)\di x(t)$. Further, any two such lifts differ by a
translation in the $z$ direction.

For each $i$, consider a piecewise continuously differentiable closed
path $\kappa=(x,y)\colon[0,1]\to\R^2_{xy}$ that traverses the boundary
of $U_i$ once counterclockwise. Suppose it meets the crossings
$b_1,\ldots,b_s$ in this order (repetition in the list is allowed).
Construct a Legendrian lift $\lambda(t)=(\kappa(t),z(t))$ of
$\kappa(t)$. This is not a closed curve; in fact by Stokes' theorem,
$z(1)-z(0)=\int_0^1 y(t)\di x(t)=\int_{U_i}\di y\wedge\di x=-f_i$. On
the other hand, the same $\lambda$ can also be constructed by gluing
together $z$--translates of pieces of $L$. Inspecting 
\figref{fig:Reeb}, it is easy to see that the total change of $z$ along
the latter is exactly $\pm h(b_1)\pm\cdots\pm h(b_s)$, where the
counting happens with the opposite of the Reeb sign.
\end{proof}

We will also need the following observation from \cite{chek} later:


\begin{lemma}\label{lem:fokok}
In the situation of Proposition \ref{teruletes},
\[\sum_{j=1}^n E_{ij}|a_j| \equiv 2 -
\#(\text{positive corners of } U_i)\pmod{2r}.\]
\end{lemma}

\begin{megj}
The equations \eqref{eq:teri} imply the well known fact that a
Lagrangian projection bounds zero area. Indeed, at each crossing, the
sum of the indices with respect to $\gamma$ of either pair of opposite
regions (sharing the same Reeb sign) is the same.
\end{megj}

\begin{megj}\label{rem:kup}
By equation \eqref{eq:teri}, the \emph{height vector}
$h=(h_1,h_2,\ldots,h_{n})\in S_\gamma^+$ of a generic Lagrangian
projection determines the \emph{area vector}
$f=(f_1,f_2,\ldots,f_{n+1})$. Moreover, $h$ is such that
$f=\Psi_\gamma h$ lies in $V_\gamma^+$. In addition to $h_i>0$ for all
$i$, this forces $n+1$ homogeneous linear inequalities for $h$. By
\emph{satisfying the linear constraints} we either mean that $h$ is
such that $h\in S_\gamma^+$ and $\Psi_\gamma h\in V_\gamma^+$, or that
the pair $(h,f)\in S_\gamma^+\times V_\gamma^+$ is such that
\eqref{eq:teri} holds.
\end{megj}

\begin{Def}
Let $\gamma$ be a classical link diagram and consider the open convex
cone
\[\mathscr C_\gamma=\{\,h\in S^+_\gamma\mid\Psi_\gamma h\in
V^+_\gamma\,\}.\]
We say that $\gamma$ is a \emph{Lagrangian diagram} if $\mathscr
C_\gamma$ is non-empty.
\end{Def}

By the previous remark, Lagrangian projections of generic Legendrian
links are Lagrangian diagrams. Lagrangian diagrams are the objects
represented in our figures; in particular, we do not draw areas of
regions to scale. All Lagrangian diagrams can, though, be
\emph{isotoped} (in the smooth sense) to actual Lagrangian projections
(or, equivalently, we can use orientation-preserving diffeomorphisms
of $\R^2$ to take one diagram to another). Moreover, \emph{any}
element of $\mathscr C_\gamma$ can be realized this way and because
$\mathscr C_\gamma$ is convex, it is easy to see that all Lagrangian
projections with the same underlying Lagrangian diagram can be
isotoped into each other (this time in the Legendrian sense, ie,\
keeping the linear constraints intact throughout). In particular, a
Lagrangian diagram represents a well-defined Legendrian link.

\begin{pelda}\label{ex:trefli1}
Consider the knot diagram $\gamma$ on the left side of
\figref{fig:trefli}. We claim that $\mathscr C_\gamma$ is non-empty,
thus $\gamma$ is the Lagrangian diagram of a Legendrian knot $L$.
The equations \eqref{eq:teri} take the following form (note the
Reeb signs in \figref{fig:trefli}):
\[\begin{array}{rcrclclclcl}
f_1&=&&&h(a_2)&&&&&&\\
f_2&=&-2h(a_1)&+&h(a_2)&+&h(b_1)&&&+&h(b_3)\\
f_3&=&h(a_1)&&&&&&&&\\
f_4&=&h(a_1)&&&-&h(b_1)&-&h(b_2)&-&h(b_3)\\
f_5&=&&&&&h(b_1)&+&h(b_2)&&\\
f_6&=&&&&&&&h(b_2)&+&h(b_3).
\end{array}\]
It is easy to check that the choice of $h(a_1)=4$, $h(a_2)=7$,
$h(b_1)=h(b_2)=h(b_3)=1$ yields the values $f_1=7$, $f_2=1$, $f_3=4$,
$f_4=1$, and $f_5=f_6=2$. As these are all positive, our linear
constraints can indeed be satisfied.

\begin{figure}[ht!]\small\anchor{fig:trefli}
\psfrag{a1}{$a_1$}
\psfrag{a2}{$a_2$}
\psfrag{b1}{$b_1$}
\psfrag{b2}{$b_2$}
\psfrag{b3}{$b_3$}
\psfrag{U1}{$U_1$}
\psfrag{U2}{$U_2$}
\psfrag{U3}{$U_3$}
\psfrag{U4}{$U_4$}
\psfrag{U5}{$U_5$}
\psfrag{U6}{$U_6$}
\cl{\includegraphics[width=.9\hsize]{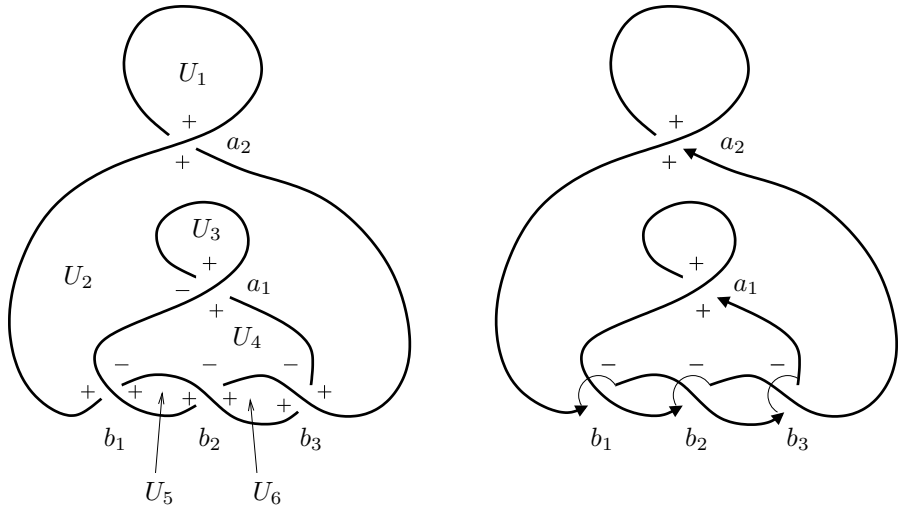}}
\caption{Right-handed Legendrian trefoil knot}\label{fig:trefli}
\end{figure}

The claim can be proven in a less direct way, too: the diagram
$\gamma$ is the result of the resolution \cite{computable} of a
certain front diagram, analogous to the base point of the loop of
\figref{fig:hurok}. See section \ref{sec:BC} for more.

The rotation number of this knot (with either orientation, of course)
is $r(L)=0$. The gradings of the crossings (with either orientation,
because $r=0$) are $|a_1|=|a_2|=1$ and $|b_1|=|b_2|=|b_3|=0$. Hence
the Thurston--Bennequin number is $tb(L)=1$. Finally, the classical
knot type of $L$ is the right-handed trefoil knot.
\end{pelda}

\begin{megj}\label{lokalis}
Let $p$ and $q$ be arbitrary points on the same component $C$ of
the Legendrian link $L$. Then $\pi(p)$ and $\pi(q)$ divide
$\gamma=\pi(C)$ into two arcs, $\gamma_1$ and $\gamma_2$. If
$\gamma_1'$ is another immersed arc with the same endpoints and
starting and ending tangent vectors as $\gamma_1$, and so that
$\gamma_1'\cup(-\gamma_1)$ bounds zero area, then
$\gamma_1'\cup\gamma_2$ possesses a Legendrian lift (in general,
an immersion) which coincides with $C$ over $\gamma_2$. This
diagrammatical method can be used to construct modifications of
Legendrian links that only affect a small segment of the link,
although possibly quite dramatically. We will use this idea in
sections \ref{sec:simplex}, \ref{sec:BC}, and \ref{sec:loop}.
\end{megj}

\subsection{Knot DGA and homology}

\begin{Def}
Let $R$ be a commutative ring. A \emph{differential graded algebra},
or \emph{DGA}, is a graded associative algebra
$\mathscr A=\bigoplus_{i=-\infty}^\infty\mathscr A_i$ over $R$ with
identity $1\in\mathscr A_0$, together with a differential
$\partial\colon\mathscr A\to\mathscr A$. We will also refer to the
grading as the \emph{index}, which is multiplicative and takes values
in a cyclic group $\Gamma$. We require that $\partial$ be an index
$-1$, $R$--linear map so that
$\partial(ab)=\partial(a)b+(-1)^{|a|}a\partial(b)$ for all
$b\in\mathscr A$ and $a\in\mathscr A$ with pure index $|a|$
(consequently, $\partial(1)=0$), and finally that $\partial^2=0$.

A DGA $\mathscr A$ is called \emph{semi-free} (with generators
$a_1,\ldots,a_n$) if its underlying algebra is the ring (tensor
algebra) $T(a_1,\ldots,a_n)$ of non-commutative polynomials in a
finite collection of pure-index elements $a_1,\ldots,a_n\in\mathscr
A$. We will use the term \emph{degree} to refer to the degree of
polynomials in this situation.
\end{Def}

Throughout the paper, we will use the term polynomial to mean
non-com\-mu\-ta\-tive polynomial.

Let $L$ be an oriented, generic Legendrian knot in $\R^3$ with
Lagrangian projection $\gamma\subset\R^2_{xy}$. Let the crossings of
$\gamma$ be $a_1,a_2,\ldots,a_n$. We associate a semi-free DGA to $L$
as follows.

\textbf{The graded algebra}\qua Let $\mathscr A$ be the non-commutative,
associative, unital polynomial algebra over $R=\Zt$ generated by the
symbols $a_1,a_2,\ldots,a_n$. Grade each generator using the value
$|a_i|$ defined by \eqref{eq:grad} and let $|t|=2r$ (this value, ie\
twice the rotation number, is often called the \emph{Maslov number} of
the knot). Extend the grading to monomials multiplicatively, obtaining
a $\Z$--grading so that we may write $\mathscr A_*$. Note that this
grading, also called index, is different from the grading given by the
degree of the polynomials.

For each positive integer $k$, fix an oriented disc $\Pi_k$ (thus the
boundary $\partial\Pi_k$ is also oriented) with $k$ points
$x_0^k,x_1^k,\ldots,x_{k-1}^k$ marked on its boundary in the given
cyclic order. A continuous map $f\colon\Pi_k\to\R^2_{xy}$ is
\emph{admissible with respect to $\gamma$} (sometimes we'll say $f$ is
an \emph{admissible disc in $\gamma$}) if it is an immersion away from 
the marked points, it preserves the orientation
of the disc and moreover it takes $\partial\Pi_k$ to $\gamma$ so that each 
marked
point is mapped to a crossing in such a way that locally, the image of
the disc forms an angle less than $180^\circ$. We also require that
the quadrant occupied near the crossing $f(x_0^k)$ (the so called
\emph{positive corner}) has Reeb sign $+1$ and all other quadrants, at
the so called \emph{negative corners} $f(x_1^k),\ldots,f(x_{k-1}^k)$,
have Reeb sign $-1$. We'll also say that the disc \emph{turns} at its
corners $f(x_0^k),f(x_1^k),\ldots,f(x_{k-1}^k)$, or that it turns at
certain quadrants located at those points.

Note that we didn't require that the orientations of the
curves $\gamma$ and $f(\partial\Pi_k)$ match. Let us call an
admissible immersion $f$ \emph{compatible with the orientation of
$\gamma$} (or a \emph{compatible disc}) if those two orientations
either agree at all points of $\partial\Pi_k$ or disagree at all
such points. We will need the following lemma in section
\ref{sec:BC}.

\begin{lemma}\label{lem:kovet}
Let $\gamma$ be an oriented Lagrangian diagram. An admissible
immersion $f$ is compatible with the orientation of $\gamma$ if and
only if its positive corner has odd index and all of its negative
corners have even indices.
\end{lemma}

\begin{proof}
An immersion is compatible if and only if the orientation requirement is
satisfied in arbitrarily small neighborhoods of the punctures. Then
the statement follows very easily from the observation (also mentioned
in footnote \ref{lab'}) that $(-1)^{|a|}$ coincides with the sign of
the crossing $a$ in classical knot theory. Indeed, the middle diagram
of \figref{fig:Reeb} shows an oddly graded crossing and we see
(comparing with the diagram on the left) that a compatible disc may
only turn at its positive quadrants. The diagram on the right shows an
evenly graded crossing and a similar examination yields that
compatible discs may only turn at its negative quadrants.
\end{proof}

The following lemma is a consequence of Proposition \ref{teruletes}.

\begin{lemma}\label{lem:ugras}
Let $f\colon\Pi_k\to\R^2_{xy}$ be an admissible immersion. Then
\[h(f(x_0^k))=\int_{\Pi_k}f^*(\di x\wedge\di y)
+\sum_{i=1}^{k-1}h(f(x_i^k)).\]
In particular, $h(f(x_0^k))>h(f(x_i^k))$ for all $i=1,\ldots,k-1$.
\end{lemma}

\textbf{The differential}\qua To compute $\partial(a)$ for a given
crossing $a$, one considers all admissible immersions, modulo
reparametrization, with positive corner $a$. A typical one described
above contributes the monomial
\[t^{-l}\cdot\varepsilon_0\varepsilon_1\ldots\varepsilon_{k-1}
f(x_1^k)\ldots f(x_{k-1}^k),\]
where $\varepsilon_0,\varepsilon_1,\ldots,\varepsilon_{k-1}$ are the
orientation signs (see \figref{fig:Reeb}) of the quadrants
occupied by $f(\Pi_k)$ at its corners, and $l$ is defined as follows.
Take the union of the oriented segments of $f(\partial\Pi_k)$ (the
orientation is induced from that of $\partial\Pi_k$) with the capping
path $\gamma_{f(x_0^k)}=\gamma_a$ and the reversed capping paths
$-\gamma_{f(x_1^k)},\ldots,-\gamma_{f(x_{k-1}^k)}$. This is a cycle
whose homology class in $H(L;\Z)\cong\Z$ is identified with the
integer $l$ via the chosen orientation of $L$.

$\partial(a)$ is the sum of these contributions over all admissible
immersions. We extend $\partial$ to $\mathscr A$ using
$\Zt$--linearity and the signed Leibniz rule
$\partial(ab)=\partial(a)b+(-1)^{|a|}a\partial(b)$.

\begin{tetel}[Chekanov--Etnyre--Ng--Sabloff]
The sum taken when computing $\partial(a)$ for any crossing $a$ is
finite, hence $\partial\colon\mathscr A\to\mathscr A$ is well defined.
It has index $-1$ and satisfies $\partial^2=0$.
\end{tetel}

\begin{Def}
The semi-free DGA $(\mathscr A_*,\partial)$ is called the
\emph{Che\-ka\-nov--Eli\-ash\-berg differential graded algebra}
associated to the generic Legendrian knot $L$. Its homology
$H(L)=\ker(\partial)/\im(\partial)$ is the
\emph{Chekanov--Eliashberg contact homology}, or simply
\emph{contact homology} of $L$.
\end{Def}

Chekanov's main theorem is that $H(L)$ (and further, a certain
equivalence class of $(\mathscr A_*,\partial)$, called its stable tame
isomorphism class) is unchanged when one passes to another generic
Legendrian knot $L'$ that is homotopic to $L$ through Legendrian
knots. (In this case we also say that $L'$ is \emph{Legendrian
isotopic} to $L$.) We will state this result in Theorem
\ref{thm:moves}.

\begin{pelda}\label{ex:trefli2}
The trefoil DGA: please refer to \figref{fig:trefli} and Example
\ref{ex:trefli1}. We choose the orientation as on the right side of
the figure. Quadrants with negative orientation sign are marked. Note
that $|t|=0$ and thus for index reasons,
$\partial(b_1)=\partial(b_2)=\partial(b_3)=0$. On the other hand,
\begin{equation}\label{eq:del1}
\partial(a_1)=1-b_1-b_3-tb_1b_2b_3.
\end{equation}
The admissible discs contributing these terms are $U_3$, $U_4\cup
U_5$, $U_4\cup U_6$, and $U_4$, respectively. Similarly,
\begin{equation}\label{eq:del2}
\partial(a_2)=1-b_2+t^{-1}+b_2b_3+b_1b_2+tb_2b_3b_1b_2,
\end{equation}
with terms contributed by\footnote{This way of describing the
actual immersions works fine in this case. In general however it
is possible that two different admissible discs cover the same
regions with the same multiplicities. As to whether we indeed
listed all admissible discs, the reader can check that as an
exercise or can refer to section \ref{sec:BC}, where we shall
determine all admissible discs in a certain class of Lagrangian
diagrams generalizing $\gamma$.} $U_1$, $U_2\cup U_3\cup U_4$,
$U_2\cup2U_4\cup U_5\cup U_6$, $U_2\cup2U_4\cup U_5$,
$U_2\cup2U_4\cup U_6$, and $U_2\cup2U_4$, respectively. Note that
all admissible immersions in $\gamma$ are compatible with the
orientation of the diagram. The following observation will be
useful later in the proof of Theorem \ref{thm:trefli} for $p=3$,
$q=2$:
\begin{equation}\label{trick1}\begin{split}
\mspace{-4mu}\partial(a_2b_3+(t^{-1}+b_2b_3)a_1)
&=b_3-b_2b_3+t^{-1}b_3+b_2b_3^2+b_1b_2b_3+tb_2b_3b_1b_2b_3\\
&+t^{-1}-t^{-1}b_1-t^{-1}b_3-b_1b_2b_3\\
&+b_2b_3-b_2b_3b_1-b_2b_3^2-tb_2b_3b_1b_2b_3\\
&=b_3+t^{-1}-t^{-1}b_1-b_2b_3b_1.
\end{split}\end{equation}
In other words, $b_3+t^{-1}-t^{-1}b_1-b_2b_3b_1=0$ in the contact
homology.
\end{pelda}

We may recover Chekanov's original DGA by keeping the same
generators and admissible discs, but putting $t=1$, reducing
coefficients modulo $2$ and the grading modulo $2r$. For example,
in this theory, \eqref{eq:del1} and \eqref{eq:del2} reduce to
$\partial(a_1)=1+b_1+b_3+b_1b_2b_3$ and
$\partial(a_2)=b_2+b_2b_3+b_1b_2+b_2b_3b_1b_2$, respectively.

\subsection{DGA maps and the product structure}

Chain complexes of homology theories are primarily Abelian groups with
a differential. The presence of the (non-commutative) product
structure in Che\-ka\-nov--Eli\-ash\-berg theory is an `extra' feature. In
this section we point out how this does not actually complicate
matters; namely, how we can concentrate on generators only and then
let the algebra take care of itself. For more on general DGA theory,
see \cite{handbook}.

\begin{Def}
Let $(\mathscr A,\partial)$ and $(\mathscr B,\partial')$ be
differential graded algebras (over the same commutative ring $R$ and
graded by the same cyclic group $\Gamma$). A linear map
$\varphi\colon\mathscr A\to\mathscr B$ is a \emph{chain map} if it is
index-preserving and intertwines the differentials:
$\varphi(\partial(x))=\partial'(\varphi(x))$ for all $x\in\mathscr A$.
If such a $\varphi$ is also an algebra homomorphism, then we call it a
\emph{DGA morphism}.
\end{Def}

\begin{Def}
Let $(\mathscr A,\partial)$ and $(\mathscr B,\partial')$ be
differential graded algebras (over the same $R$, graded by the same
$\Gamma$), and $\varphi,\psi\colon\mathscr A\to\mathscr B$ chain maps
between them. An index $r$ map $S\colon\mathscr A\to\mathscr B$ is
called a \emph{$(\varphi,\psi)$--derivation of index $r$} if for all
$a\in\mathscr A$ of pure index $|a|$ and for all $b\in\mathscr A$, we
have
\[S(ab)=S(a)\psi(b)+(-1)^{r|a|}\varphi(a)S(b).\]
If $(\mathscr A,\partial)=(\mathscr B,\partial')$ and
$\varphi=\psi=\id_{\mathscr A}$, then we simply speak of a derivation
of index $r$.
\end{Def}

In this paper we consider three types of DGA mappings. They are all
$R$--linear, but each has a different relation to the product
structure of the DGA. Namely, by definition, the \emph{differential}
itself is a derivation of index $-1$. All \emph{chain maps} that we
consider are going to be DGA morphisms (in particular, they fix
constants), and all \emph{chain homotopies} connecting chain maps
$\varphi$, $\psi$ will be $(\varphi,\psi)$--derivations of index $1$
(from which it follows that they take constants to $0$). The following
two lemmas show that these self-imposed restrictions are very natural
and useful.

\begin{lemma}\label{lem:eta}
Let $(\mathscr A,\partial)$ and $(\mathscr B,\partial')$ be semi-free
differential graded algebras over a commutative ring $R$, graded by
$\Gamma$. If an index $0$ algebra homomorphism
$\eta\colon\mathscr A\to\mathscr B$ is such that
$\eta(\partial(a))=\partial'(\eta(a))$ for all generators of
$\mathscr A$, then the same holds for all elements $a\in\mathscr A$
(ie, $\eta$ is a chain map).

Moreover, any assignment of values $\eta(a)\in\mathscr B$ to the
generators $a$ of $\mathscr A$ which satisfies $|\eta(a)|=|a|$ for all
$a$ can be uniquely extended to an index $0$ algebra homomorphism.
\end{lemma}

\begin{proof}
To complete a trivial induction proof, we need to show that if $a$ and
$b$ are homogeneous elements of $\mathscr A$ that satisfy the desired
property, then such is $ab$, too:
\begin{multline*}
\eta(\partial(ab))
=\eta\left(\partial(a)b+(-1)^{|a|}a\partial(b)\right)
=\eta(\partial(a))\eta(b)+(-1)^{|a|}\eta(a)\eta(\partial(b))\\
=\partial'(\eta(a))\eta(b)+(-1)^{|\eta(a)|}\eta(a)\partial'(\eta(b))
=\partial'(\eta(a)\eta(b))
=\partial'(\eta(ab)).
\end{multline*}
In the last assertion, the existence of an extension is obvious, since
$\mathscr A$ is semi-free with the given generators and $\mathscr B$
is associative. We just need to check that the unique extension is
also of index $0$, which is established with the same induction
argument based on noting that if $|\eta(a)|=|a|$ and $|\eta(b)|=|b|$,
then $|\eta(ab)|=|\eta(a)\eta(b)|=|\eta(a)|+|\eta(b)|=|a|+|b|=|ab|$.
\end{proof}

\begin{lemma}\label{lem:K}
Let $(\mathscr A,\partial)$ and $(\mathscr B,\partial')$ be semi-free
differential graded algebras over the commutative ring $R$, graded by
$\Gamma$, and let $\varphi,\psi\colon\mathscr A\to\mathscr B$ be DGA
morphisms. If $K\colon\mathscr A\to\mathscr B$ is a
$(\varphi,\psi)$--derivation of index $1$ so that
$K(\partial(a))+\partial'(K(a)) = \varphi(a) - \psi(a)$ for all
generators of $\mathscr A$, then the same holds for all elements
$a\in\mathscr A$ (ie, $K$ is a chain homotopy between $\varphi$ and
$\psi$).

Moreover, any assignment of values $K(a)\in\mathscr B$ to the
generators $a$ of $\mathscr A$ which satisfies $|K(a)|=|a|+1$ for all
$a$ can be uniquely extended to a $(\varphi,\psi)$--derivation of
index $1$.
\end{lemma}

\begin{proof}
The first part of the proof is a similar, but longer computation as in
the previous argument:
\begin{equation*}\begin{split}
(K\partial&+\partial'K)(ab)
=K(\partial(ab))+\partial'(K(ab))\\
&=K\left(\partial(a)b+(-1)^{|a|}a\partial(b)\right)
+\partial'\left(K(a)\psi(b)+(-1)^{|a|}\varphi(a)K(b)\right)\\
&=K(\partial(a))\psi(b)+(-1)^{|a|-1}\varphi(\partial(a))K(b)\\
&+(-1)^{|a|}K(a)\psi(\partial(b))+(-1)^{2|a|}\varphi(a)K(\partial(b))\\
&+\partial'(K(a))\psi(b)+(-1)^{|a|+1}K(a)\partial'(\psi(b))\\
&+(-1)^{|a|}\partial'(\varphi(a))K(b)
+(-1)^{|a|+|\varphi(a)|}\varphi(a)\partial'(K(b))\\
&=\left(K(\partial(a))+\partial'(K(a))\right)\psi(b)
+(-1)^{|a|}\left(-\varphi(\partial(a))+\partial'(\varphi(a))\right)K(b)\\
&+(-1)^{|a|}K(a)\left(\psi(\partial(b))-\partial'(\psi(b))\right)
+\varphi(a)\left(K(\partial(b))+\partial'(K(b))\right)\\
&=(\varphi(a)-\psi(a))\psi(b)+0+0+\varphi(a)(\varphi(b)-\psi(b))\\
&=\varphi(a)\varphi(b)-\psi(a)\psi(b)
=\varphi(ab)-\psi(ab).
\end{split}\end{equation*}
In the second part, to show that the extension is unique, we need
to establish that $K((ab)c)=K(a(bc))$, which is due to the
following direct computation:
\begin{gather*}
\begin{gathered}[t]
\shoveleft{K((ab)c)=K(ab)\psi(c)+(-1)^{|a|+|b|}\varphi(ab)K(c)}\\
\shoveright{=K(a)\psi(b)\psi(c)+(-1)^{|a|}\varphi(a)K(b)\psi(c)+(-1)^{|a|+|b|}
\varphi(a)\varphi(b)K(c).}
\end{gathered}\\
\begin{gathered}[t]
\shoveleft{K(a(bc))=K(a)\psi(bc)+(-1)^{|a|}\varphi(a)K(bc)}\\
\shoveright{=K(a)\psi(b)\psi(c)+(-1)^{|a|}\varphi(a)K(b)\psi(c)+(-1)^{|a|}(-1)^{|b|}
\varphi(a)\varphi(b)K(c).}
\end{gathered}
\end{gather*}
The extension is also index $1$ by the already usual induction
argument: If $a$ and $b$ are homogeneous with $|K(a)|=|a|+1$ and
$|K(b)|=|b|+1$, then
$|K(ab)|=|K(a)\psi(b)+(-1)^{|a|}\varphi(a)K(b)|$, where
$|\psi(b)|=|b|$ and $|\varphi(a)|=|a|$, so $K(ab)$ is homogeneous
of index $|a|+|b|+1=|ab|+1$.
\end{proof}

Table \fref{algebra} summarizes the relevant algebraic machinery. For
all three rows, it is true that if the `map' is defined on the
generators so that the `index' requirement is satisfied and then if it
is extended to the whole DGA using linearity and the `product
identity,' then
\begin{enumerate}
\item such an extension is unique, it satisfies the `index'
requirement, and
\item if the `property' requirement holds for the generators, then it
holds for all elements of the DGA.
\end{enumerate}

\begin{table}\anchor{algebra}
\small\def\strutt{\vrule width 0pt height 12pt}
\def\botstrutt{\vrule width 0pt depth 7pt}
\cl{\begin{tabular}{|c|c|c|c|}
\hline\strutt
map & property & \hfill index & product identity\botstrutt \\
\hline\strutt
differential $\partial$ & $\partial^2=0$ & $-1$ &
$\partial(ab)=\partial(a)b+(-1)^{|a|}a\partial(b)$\botstrutt \\
\hline\strutt
chain map $\eta$ & $\eta\partial=\partial'\eta$ & $\enskip0\enskip$ &
$\eta(ab)=\eta(a)\eta(b)$\botstrutt \\
\hline\strutt
chain homotopy $K$&$K\partial +\partial'K$ &&\\
between chain&$= \varphi - \psi$ &1&
$K(ab)=K(a)\psi(b)+(-1)^{|a|}\varphi(a)K(b)$\\
maps $\varphi$, $\psi$ &&&\botstrutt\\
\hline
\end{tabular}}
\caption{Product identities}\label{algebra}
\end{table}

\section{The monodromy invariant}\label{sec:ind}

Let us first outline the contents of the section. Let
$\mathscr L\subset\leg$ be a Legendrian knot type. We call an element 
of $\mathscr L$ \emph{generic}, as in section \ref{sec:prep}, if its 
Lagrangian projection is generic (this is equivalent to saying that 
the Legendrian knot doesn't have degenerate Reeb chords). These are 
the ones to which a graded chain complex $(\mathscr A_*,\partial)$, 
and a corresponding Floer homology $H$ are actually associated. 
Non-generic objects form a codimension $1$ discriminant $\mathscr 
D\subset\mathscr L$.

Generic homotopies, ie paths in $\mathscr L$, meet only the
codimension $1$ (codimension is always meant with respect to
$\mathscr L$) stratum $\mathscr D_1\subset\mathscr D$ and 
diagrammatically, this corresponds to Reidemeister II, II$^{-1}$, or 
III moves (see \figref{fig:moves}). If $L$ and $L'$ are in 
adjacent chambers of $\mathscr L\setminus\mathscr D$, then following 
\cite{chek} (and \cite{ENS} for the case of 
$\Z[t,t^{-1}]$--coefficients), one is able to write down a chain map 
$\varphi\colon\mathscr A(L)\to\mathscr A(L')$, corresponding to the 
Reidemeister move, that induces an isomorphism 
$\varphi_*\colon H(L)\to H(L')$. (In \cite{chek,ENS} they use these 
chain maps to prove that the homology is independent of the actual 
projection.) Composing these chain maps/isomorphisms, every generic 
path $\Phi(t)$ that connects generic objects $L_0$ and $L_1$ has an 
induced isomorphism $\Phi_*\colon H(L_0)\to H(L_1)$, called the 
holonomy of the path. Details are found in section \ref{sec:cm}.

Now to prove that the holonomy only depends on the homotopy class 
(with fixed endpoints) of $\Phi(t)$, one has to consider two kinds of 
events: paths that are tangent to $\mathscr D_1$ (this is easy) and 
paths that meet $\mathscr D_2$. To deal with the latter, we take a 
small generic loop $\Psi(t)$, starting and ending at the generic 
element $L$, and linking the codimension $2$ discriminant
$\mathscr D_2\subset \mathscr D$ once. While it is not necessarily the 
case that the induced \emph{chain map} 
$\psi\colon\mathscr A_*(L)\to\mathscr A_*(L)$ is the identity, we will 
establish that there is always a \emph{chain homotopy} 
$K\colon\mathscr A_*(L)\to\mathscr A_{*+1}(L)$ connecting $\psi$ and 
$\id_{\mathscr A_*(L)}$, thereby proving that \emph{on the homology 
level} $\psi_*=\Psi_*=\id_{H(L)}$, as desired.

\begin{figure}[ht!]\small\anchor{fig:hasab}
\psfraga <-2pt,0pt> {f'}{$\varphi'$}
\psfraga <-2pt,0pt> {f''}{$\varphi''$}
\psfraga <-2pt,0pt> {f'*}{$\varphi'_*$}
\psfraga <-2pt,0pt> {f''*}{$\varphi''_*$}
\psfraga <-2pt,0pt> {p'}{$\psi'$}
\psfraga <-2pt,0pt> {p''}{$\psi''$}
\psfraga <-2pt,0pt> {p'*}{$\psi'_*$}
\psfraga <-2pt,0pt> {p''*}{$\psi''_*$}
\psfraga <-2pt,0pt> {F'}{$\Phi'$}
\psfraga <-2pt,0pt> {F''}{$\Phi''$}
\psfraga <-2pt,0pt> {P'}{$\Psi'$}
\psfraga <-2pt,0pt> {P''}{$\Psi''$}
\psfraga <-2pt,0pt> {D1}{$\mathscr D_1$}
\psfraga <-2pt,0pt> {D2}{$\mathscr D_2$}
\psfraga <-2pt,0pt> {K}{$K$}
\psfraga <-2pt,0pt> {f}{$\varphi$}
\psfraga <-2pt,0pt> {p}{$\psi$}
\psfraga <-2pt,0pt> {L1}{$L_{11}$}
\psfraga <-2pt,0pt> {L2}{$L_{21}$}
\psfraga <-2pt,0pt> {L3}{$L_{22}$}
\psfraga <-2pt,0pt> {L4}{$L_{12}$}
\psfraga <-2pt,0pt> {A1}{$\mathscr A(L_{11})$}
\psfraga <-2pt,0pt> {A2}{$\mathscr A(L_{21})$}
\psfraga <-2pt,0pt> {A3}{$\mathscr A(L_{12})$}
\psfraga <-2pt,0pt> {A4}{$\mathscr A(L_{22})$}
\psfraga <-2pt,0pt> {H1}{$H(L_{11})$}
\psfraga <-2pt,0pt> {H2}{$H(L_{21})$}
\psfraga <-2pt,0pt> {H3}{$H(L_{12})$}
\psfraga <-2pt,0pt> {H4}{$H(L_{12})$}
\cl{\includegraphics[width=.7\hsize]{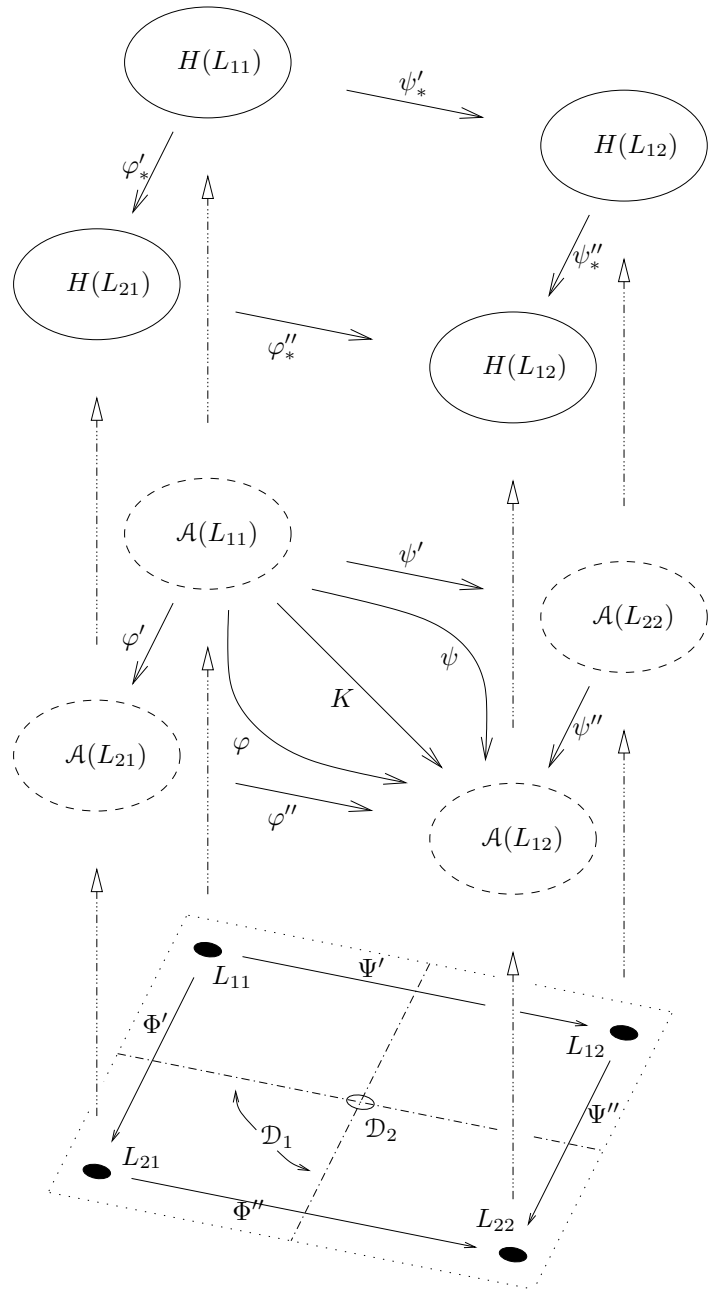}}
\caption{Chain complexes and homologies of different
objects}\label{fig:hasab}
\end{figure}

In \figref{fig:hasab}, we represented a hypothetical case when two 
strands of the discriminant meet transversally. (As it is shown in 
\cite[Figure 6]{arnold} and as we'll see in section \ref{sec:ch}, the 
local description at other points of $\mathscr D_2$ may be slightly 
less or more complicated, but the approach is the same.) As we know 
that all arrows on the `top level' represent isomorphisms, we may 
choose any two objects, which may even coincide, and direct all arrows 
from the first to the second. If this choice of directions can be 
accompanied with a chain homotopy on the `middle level,' and in 
section \ref{sec:ch} we prove that this is always the case, then we 
indeed get that the diagram of holonomies on the `top level' commutes.

In section \ref{sec:z2} we review our formulas in their simplified, 
$\Z_2$--coefficient versions.

\subsection{Chain maps}\label{sec:cm}

\begin{figure}[ht!]\small\anchor{fig:moves}
\psfraga <-2pt,0pt> {a}{$a$}
\psfraga <-2pt,0pt> {b}{$b$}
\psfraga <-2pt,0pt> {c}{$c$}
\psfraga <-2pt,0pt> {3a}{{III$_{\text{a}}$}}
\psfraga <-2pt,0pt> {3b}{{III$_{\text{b}}$}}
\psfraga <-2pt,0pt> {-1}{{II}}
\psfraga <-2pt,0pt> {2}{{II$^{-1}$}}
\psfraga <-2pt,0pt> {a'}{$a'$}
\psfraga <-2pt,0pt> {b'}{$b'$}
\psfraga <-2pt,0pt> {c'}{$c'$}
\psfraga <-2pt,0pt> {e1}{$\varepsilon'_a$}
\psfraga <-2pt,0pt> {e2}{$\varepsilon_a$}
\psfraga <-2pt,0pt> {e3}{$\varepsilon_b$}
\psfraga <-2pt,0pt> {e4}{$\varepsilon'_b$}
\cl{\includegraphics[width=.9\hsize]{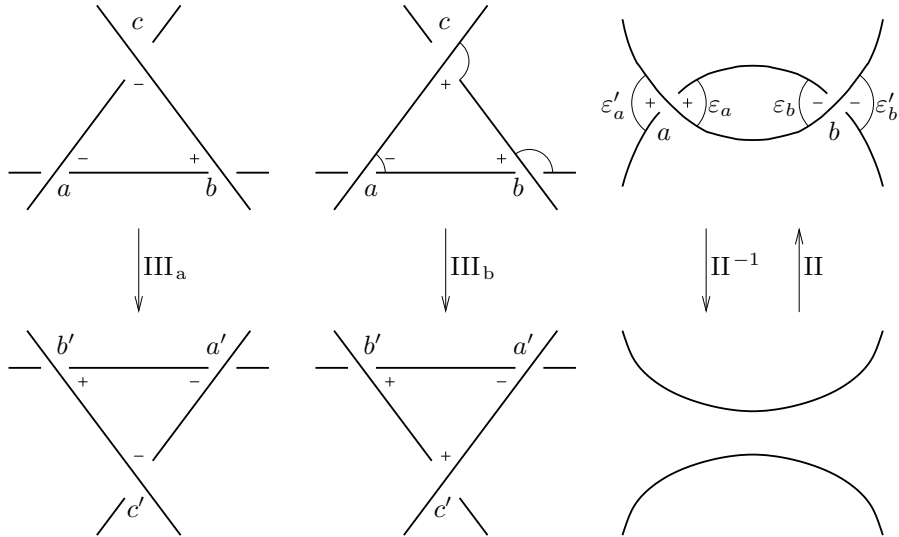}}
\caption{Legendrian Reidemeister moves in the Lagrangian projection. The 
$+$ and $-$ signs are Reeb signs and $\varepsilon_a$, $\varepsilon'_a$, 
$\varepsilon_b$,
$\varepsilon'_b$ refer to orientation signs.}\label{fig:moves}
\end{figure}

Generic paths in $\leg$ are isotopies with a finite number of
Reidemeister moves. If we work with Lagrangian projections, then there 
are four different kinds of the latter,
depicted in \figref{fig:moves}. III$_{\text{a}}$ and
III$_{\text{b}}$ moves are also called \emph{triangle moves}. Move II
is sometimes referred to as a \emph{birth} and move II$^{-1}$ as a
\emph{death} (of a pair of crossings). 

Let us denote the DGA's of the
diagrams in the upper row of \figref{fig:moves} by $(\mathscr 
A,\partial)$ and in the lower
row, by $(\mathscr A',\partial')$. The crossings not affected by the
moves are in an obvious one-to-one correspondence. If $x\in\mathscr A$
is one of these, then the corresponding crossing in $\mathscr A'$ will
be denoted by $x'$.

Before stating Theorem \ref{thm:moves}, let us review Chekanov's
analysis, although in slightly different terms, of the bifurcation
that takes place on the right side of \figref{fig:moves}. The
re-phrasing is possible because we shall only be concerned with the
homology, and not the stable tame isomorphism class of the DGA, as our
invariant. Let the crossings that generate $\mathscr A$ be arranged by
height as follows:
\[h(a_l)\ge\ldots\ge h(a_1)\ge h(a)>h(b)\ge h(b_1)\ge\ldots\ge
h(b_m).\]
Let $\partial(a)=\varepsilon_a\varepsilon_bb+v$, where $\varepsilon_a$
and $\varepsilon_b$ are the orientation signs from
\figref{fig:moves} (it is easy to prove that the exponent of $t$ is $0$
in the term coming from the vanishing bigon) and by Lemma
\ref{lem:ugras}, $v\in T(b_1,\ldots,b_m)$. Define
\[\tau\colon\mathscr A\to\mathscr A'\]
by putting $\tau(a_i)=a_i'$ ($1\le i\le l$), $\tau(b_j)=b_j'$
($1\le j\le m$), $\tau(a)=0$, and
$\tau(b)=\varepsilon_a'\varepsilon_b'v'$, where $\varepsilon_a'$ and
$\varepsilon_b'$ are the other two orientation signs in 
\figref{fig:moves} and $v'$ is obtained from $v$ in the obvious way of
putting primes on all the $b_j$'s in the expression. Note that these
assignments preserve the grading, so by the second statement of Lemma
\ref{lem:eta}, we can extend $\tau$ to a well-defined index $0$
algebra homomorphism. The geometric content of the bifurcation is
captured by the following result of Chekanov:

\begin{lemma}\label{lem:tau}
$\tau\colon\mathscr A\to\mathscr A'$ is a chain map.
\end{lemma}

One of the statements of Theorem \ref{thm:moves} is that $\tau$
induces an isomorphism of homologies. Our goal now is to define what
will be its homotopy inverse, that is another chain map\footnote{All
of the considerations leading up from this point to Theorem
\ref{thm:moves} are purely algebraic. The chain map $\varphi$ is less
explicit in \cite{chek} and \cite{ENS} than in what follows.}
$\varphi\colon\mathscr A'\to\mathscr A$. While it will be the case
that $\tau\circ\varphi=\id_{\mathscr A'}$, $\varphi\circ\tau$ and
$\id_{\mathscr A}$  will only be chain homotopic via a certain
$K\colon\mathscr A\to\mathscr A$. In fact, we will build $\varphi$ and
$K$ with a somewhat subtle simultaneous construction\footnote{In
Remark \ref{kiszoroz}, we express $\varphi$ in even more explicit
terms and without referring to $K$. That description is needed for
applications, but would be impractical in the theoretical
discussion.}, using the filtration that we will introduce next.

Define
\[\mathscr A_i=T(b_1,\ldots,b_m,b,a,a_1,\ldots,a_i)\]
$$\mathscr A'_i=T(b_1',\ldots,b_m',a_1',\ldots,a_i')\leqno{\rm and}$$
for all $i=0,1,\ldots,l$. Note that by Lemma \ref{lem:ugras},
$\partial(a_i)\in\mathscr A_{i-1}$ and
$\partial(a_i')\in\mathscr A_{i-1}'$. It follows that
$\big(\mathscr A_i,\partial\big|_{\mathscr A_i}\big)$ and
$\big(\mathscr A_i',\partial'\big|_{\mathscr A_i'}\big)$ are DGA's
themselves\footnote{It is not true though that, for example,
$\partial$ would map $\mathscr A_i$ to $\mathscr A_{i-1}$;
$\partial(a_i^2)$ is usually not in $\mathscr A_{i-1}$.} for all
$i=0,1,\ldots,l$. Furthermore, for all such $i$, the restriction
\[\tau\colon\mathscr A_i\to\mathscr A_i'\]
is a DGA morphism. Define the index $0$ algebra homomorphism
\[\varphi\colon\mathscr A_0'\to\mathscr A_0\]
by $\varphi(b'_j)=b_j$ ($1\le j\le m$). By Lemma \ref{lem:eta}, this
is also a chain map, for
$\varphi(\partial'(b_j'))=\varphi(\partial'(\tau(b_j)))
=\varphi(\tau(\partial(b_j)))=\partial(b_j)=\partial(\varphi(b_j'))$.
Here we used Lemma \ref{lem:tau} and that by Lemma \ref{lem:ugras},
$\partial(b_j)$ can't contain $a$ or $b$, and finally, the obvious
fact that
$\varphi\circ\tau\big|_{T(b_1,\ldots,b_m)}=\id_{T(b_1,\ldots,b_m)}$.

Now, $\varphi\circ\tau$ and $\id_{\mathscr A_0}$ are both DGA
endomorphisms of $\mathscr A_0$. By the second statement of Lemma
\ref{lem:K}, the last three rows of the assignment
\begin{equation}\label{eq:K}
\begin{array}{rcll}
a_i&\mapsto&0&(i=1,\ldots,l)\\
a&\mapsto&0&\\
b&\mapsto&\varepsilon'_a\varepsilon'_b a&\\
b_j&\mapsto&0&(j=1,\ldots,m)
\end{array}\end{equation}
can be uniquely extended to a
$(\varphi\circ\tau,\id_{\mathscr A_0})$--derivation
\[K\colon\mathscr A_0\to\mathscr A_0\]
of index $1$; note that indeed, because $\partial$ is of index $-1$,
we have $|a|=|b|+1$. Then we use the first part of Lemma \ref{lem:K}
to show that $K$ is a chain homotopy between $\varphi\circ\tau$ and
$\id_{\mathscr A_0}$: For the generators $b_1,\ldots,b_m$, the
condition holds by the observation at the end of the previous
paragraph, together with the fact, following from the definition, that
the restriction of $K$ to the DGA $T(b_1,\ldots,b_m)$ is $0$. For $a$
and $b$, it is established below:
\begin{multline*}
K(\partial(a))+\partial(K(a))=K(\varepsilon_a\varepsilon_b b+v)+0\\
=\varepsilon_a\varepsilon_b\varepsilon_a'\varepsilon_b'a+K(v)
=(-1)^{|a|}(-1)^{|b|}a+0=-a
\end{multline*}
by Lemma \ref{elojellemma}, and
\[K(\partial(b))+\partial(K(b))=0+\partial(\varepsilon'_a\varepsilon'_b a)
=\varepsilon'_a\varepsilon'_b(\varepsilon_a\varepsilon_b b+v)
=-b+\varepsilon'_a\varepsilon'_b v.\]
On the other hand,
\[\varphi(\tau(a))-a=\varphi(0)-a=-a\text{ and }
\varphi(\tau(b))-b=\varphi(\varepsilon'_a\varepsilon'_b v')-b
=\varepsilon'_a\varepsilon'_b v-b.\]
Next, we use a recursive process to define extensions of $\varphi$ and
$K$ to the whole of $\mathscr A'$ and $\mathscr A$, respectively.
Suppose that for some value $0\le i<l$, the DGA morphism
$\varphi\colon\mathscr A_i'\to\mathscr A_i$ and the
$(\varphi\circ\tau,\id_{\mathscr A_i})$--derivation
$K\colon\mathscr A_i\to\mathscr A_i$ of index $1$ are already defined
so that $K$ is a chain homotopy between
$\varphi\circ\tau\big|_{\mathscr A_i}$ and $\id_{\mathscr A_i}$.
Extend $\varphi$ to $\mathscr A_{i+1}$ by letting
\begin{equation}\label{aelm}
\varphi(a'_{i+1})=a_{i+1}+K(\partial(a_{i+1})).
\end{equation}
As this preserves the grading of $a'_{i+1}$, by Lemma \ref{lem:eta}
(second part) there does exist an extension, which is an index $0$
algebra homomorphism $\varphi\colon\mathscr A_{i+1}'\to\mathscr
A_{i+1}$. To see that it's also a DGA morphism, we use Lemma
\ref{lem:eta} (first part), the inductive hypothesis, and the
following computation:
\begin{equation*}\begin{split}
\varphi(\partial'(a'_{i+1}))&=\varphi(\partial'(\tau(a_{i+1})))=
\varphi(\tau(\partial(a_{i+1})))\\
&=(K\circ\partial+\partial\circ K)(\partial(a_{i+1})
+\partial(a_{i+1})\\
&=\partial\left(K(\partial(a_{i+1}))+a_{i+1}\right)
=\partial(\varphi(a_{i+1}')),
\end{split}\end{equation*}
which works because $\partial(a_{i+1})\in\mathscr A_i$ and
$\partial^2=0$.

Now we are in a position to use the DGA morphisms
$\varphi\circ\tau\colon\mathscr A_{i+1}\to\mathscr A_{i+1}$ and
$\id_{\mathscr A_{i+1}}$, along with the relevant values from
\eqref{eq:K} to define the
$(\varphi\circ\tau,\id_{\mathscr A_{i+1}})$--derivation
$K\colon\mathscr A_{i+1}\to\mathscr A_{i+1}$ of index $1$, which is
obviously an extension of the previous $K$. To make sure that $K$ is
also a chain homotopy between $\varphi\circ\tau$ and
$\id_{\mathscr A_{i+1}}$, we still have to check that for $a_{i+1}$,
but it's a tautology:
$(K\circ\partial+\partial\circ K)(a_{i+1})=K(\partial(a_{i+1}))
=\varphi(a'_{i+1})-a_{i+1}=\varphi(\tau(a_{i+1}))-a_{i+1}$.

\begin{tetel}[Chekanov]\label{thm:moves}
The following maps, when extended to $\mathscr A$ as algebra
homomorphisms, are chain maps that induce isomorphisms (as graded 
algebras) of the
corresponding Chekanov--Eliashberg contact homologies. (Recall that
$(\mathscr A,\partial)$ is the DGA of the upper row, and
$(\mathscr A',\partial')$ is the DGA of the lower row in 
\figref{fig:moves}.)
\begin{enumerate}
\item Move III$_{\text{a}}$:
Let $a\mapsto a'$, $b\mapsto b'$, $c\mapsto c'$, and $x\mapsto x'$,
where $x$ is any generator of $\mathscr A$ not shown in 
\figref{fig:moves}.
\item Move III$_{\text{b}}$:
Let
\[a\mapsto a'-\varepsilon_a\varepsilon_b\varepsilon_c t^k c'b',\]
where $\varepsilon_a$, $\varepsilon_b$, and $\varepsilon_c$ are the
orientation signs of the quadrants indicated in \figref{fig:moves}
and $k=k(a;b,c)\in\Z$ is the winding, with respect to the chosen
orientation of the knot, of the union of capping paths
$-\gamma_a+\gamma_b+\gamma_c$ (thought of as starting and ending at
the single point that the triangle shrinks to). Other generators are
mapped trivially: $b\mapsto b'$, $c\mapsto c'$, and $x\mapsto x'$.
\item Move II$^{-1}$:
Let $\partial(a)=\pm b+v$ (recall that $v$ does not contain neither
$a$ nor $b$). Define $x\mapsto x'$, which gives rise to the obvious
re-labeling $v\mapsto v'$. Then let
\begin{eqnarray*}
a&\mapsto&0\\
b&\mapsto&\mp v'
\end{eqnarray*}
(this is the map $\tau$ of the preceding discussion).
\item Move II:
The map $\varphi\colon\mathscr A'\to\mathscr A$, constructed in the
preceding paragraphs (for a more explicit description, see Remark
\ref{kiszoroz}).
\end{enumerate}
\end{tetel}

\begin{proof}
We cite \cite{chek} and \cite{ENS} that the maps associated to
III$_{\text{a}}$, III$_{\text{b}}$, and II$^{-1}$ moves are indeed
chain maps (cf Lemma \ref{lem:tau}), and we have already seen that
$\varphi$ is one. In each of the four cases, the reverse of the move
is again a Reidemeister move, with a corresponding chain map
associated to it. It is enough to prove that the relevant pairwise
compositions of these chain maps are chain homotopic to the identity.
This is obvious for the pair of III$_{\text{a}}$ moves on the left
side of \figref{fig:moves}. For the pair of III$_{\text{b}}$ moves
in the middle, note that the reverse move sends
\[a'\mapsto a - \varepsilon'_a\varepsilon'_b\varepsilon'_c t^k cb,\]
where $\varepsilon'_a$, $\varepsilon'_b$, and $\varepsilon'_c$ are the
orientation signs of the quadrants opposite to the previously
considered ones, while the value of $k$ is the same. So, since
$\varepsilon_a\varepsilon_b\varepsilon_c\varepsilon'_a\varepsilon'_b
\varepsilon'_c=[-(-1)^{|a|}][-(-1)^{|b|}][-(-1)^{|c|}]=-1$ by Lemmas
\ref{elojellemma} and \ref{lem:fokok}, the chain maps themselves are
inverses of each other.

Next, we claim that $\tau\circ\varphi=\id_{\mathscr A'}$. As both
sides are DGA morphisms, it's enough to check for generators. The
claim is obvious for the $b_j'$ and follows from the definition for
the $a'_i$:
$\tau(\varphi(a_i'))=\tau(a_i+K(\partial(a_i)))=\tau(a_i)=a_i'$,
because $\tau\circ K=0$. This last assertion is true because it's
satisfied by all generators and $\tau\circ K$ and $0$ are both
$(\tau\circ\varphi\circ\tau,\tau)$--derivations of index $1$.

Finally, it is only the composition $\varphi\circ\tau$ that requires a
non-zero chain homotopy, but we have already constructed such a $K$
before stating the theorem.
\end{proof}

\begin{Def}
We will refer to the graded algebra isomorphisms induced by the maps 
described in
Theorem \ref{thm:moves} as the \emph{holonomies} of the
corresponding Reidemeister moves. Sometimes the chain maps
themselves will be called holonomies, too. Along generic paths,
finitely many Reidemeister moves happen and the composition of the
corresponding holonomies is called the holonomy of the path. When
such a sequence of Reidemeister moves returns to the original
starting knot, the holonomy is referred to as the \emph{monodromy}
of the loop.
\end{Def}

\begin{megj}\label{kiszoroz}
For the purposes of computing holonomies (and monodromies), one should
understand the map $\varphi\colon\mathscr A'\to\mathscr A$ in more
concrete terms. Of course, we have
\begin{equation}\label{bgyak}
\varphi(b_j')=b_j\text{ for all }j=1,\ldots,m.
\end{equation}
Next, recall that $\varphi(a_1')=a_1+K(\partial(a_1))$. Let us write
\[\partial(a_1)=\sum\kappa B_1 b B_2 b B_3 b B_4b\ldots B_k b A,\]
where $k\ge0$, $\kappa\in\Zt$,
$B_1,B_2,\ldots,B_k\in T(b_1,\ldots,b_m)$ are monomials, and in the
monomial $A\in\mathscr A_0$, every $b$ factor is preceded by an $a$
factor (the point here is that the $b$ factor right before $A$ is the
\emph{last $b$ before the first $a$}; if there is no $a$ at all, then
it's just the last $b$ in the word). The sum is taken over all
admissible discs with positive corner at $a_1$. Then it is not hard to
compute, applying the Leibniz-type rule of Table \ref{algebra} that we
used to define $K$ (recall also that
$\partial(a)=\varepsilon_a\varepsilon_b
b+v=-\varepsilon_a'\varepsilon_b' b+v$), that
\begin{equation}\label{eq:undor}\begin{split}
\varphi(a_1')&=a_1+\sum\kappa\Big(
(-1)^{|B_1|}\varepsilon_a'\varepsilon_b'B_1aB_2bB_3bB_4b\ldots
B_kbA\\
&+(-1)^{|B_1|+|B_2|}(-1)^{|b|}B_1vB_2aB_3bB_4b\ldots B_kbA\\
&+(-1)^{|B_1|+|B_2|+|B_3|}\varepsilon_a'\varepsilon_b'
B_1vB_2vB_3aB_4b\ldots B_kbA\\
&+(-1)^{|B_1|+|B_2|+|B_3|+|B_4|}(-1)^{|b|}B_1vB_2vB_3vB_4a\ldots
B_kbA\\
&+\ldots\\
&+(-1)^{|B_1|+\ldots+|B_k|}(-1)^{|b|}
\big(\varepsilon_a'\varepsilon_b'(-1)^{|b|}\big)^k
B_1vB_2vB_3vB_4v\ldots B_kaA\Big).
\end{split}\end{equation}
With this expression in hand, we may compute
$\varphi(a_2')=a_2+K(\partial(a_2))$. In fact, we find the following
recursion: If $\partial(a_i)=\sum\kappa B_1 b B_2 b B_3 b B_4b\ldots
B_k b A$, where
$\kappa\in\Zt$, $B_1,\ldots,B_k\in
T(b_1,\ldots,b_m,a_1,\ldots,a_{i-1})$, and $A$ doesn't contain any
copy of $b$ which is not preceded by a copy of $a$, then
\begin{equation}\label{megundor}\begin{split}
\varphi(a_i')&=a_i+\sum\kappa\Big(
(-1)^{|B_1|}\varepsilon_a'\varepsilon_b'\bar{B_1}aB_2bB_3bB_4b\ldots
B_kbA\\
&+(-1)^{|B_1|+|B_2|}(-1)^{|b|}\bar{B_1}v\bar{B_2}aB_3bB_4b\ldots B_kbA\\
&+(-1)^{|B_1|+|B_2|+|B_3|}\varepsilon_a'\varepsilon_b'
\bar{B_1}v\bar{B_2}v\bar{B_3}aB_4b\ldots B_kbA\\
&+(-1)^{|B_1|+|B_2|+|B_3|+|B_4|}(-1)^{|b|}\bar{B_1}v\bar{B_2}v\bar{B_3}v
\bar{B_4}a\ldots B_kbA\\
&+\ldots\\
&+(-1)^{|B_1|+\ldots+|B_k|}(-1)^{|b|}
\big(\varepsilon_a'\varepsilon_b'(-1)^{|b|}\big)^k\bar{B_1}v\bar{B_2}v
\bar{B_3}v\bar{B_4}v\ldots
\bar{B_k}aA\Big),
\end{split}\end{equation}
where $\bar{B_1},\ldots,\bar{B_k}$ are obtained from $B_1,\ldots,B_k$
by replacing each of the symbols $a_1,\ldots,a_{i-1}$ by the corresponding
polynomial from the already constructed list
$\varphi(a_1'),\ldots,\varphi(a'_{i-1})$. Note that if $\partial(a_i)$
doesn't contain any of $a_1,\ldots,a_{i-1}$, then \eqref{megundor} is
analogous to \eqref{eq:undor}.
\end{megj}

The following is an easy consequence of \eqref{bgyak} and either
\eqref{aelm} or \eqref{megundor}:

\begin{all}\label{pro:easy}
If $\varphi\colon\mathscr A'\to\mathscr A$ is the holonomy of a
Reidemeister II move and $x'\in\mathscr A'$ is a generator so that for
the corresponding $x\in\mathscr A$, we have $\partial(x)=0$, then
$\varphi(x')=x$.
\end{all}

\begin{figure}[ht!]\small\anchor{fig:loop}
\psfraga <-2pt,0pt> {a1}{$a_1$} \psfraga <-2pt,0pt> {a2}{$a_2$}
\psfraga <-2pt,0pt> {b1}{$b_1$} \psfraga <-2pt,0pt> {b2}{$b_2$}
\psfraga <-2pt,0pt> {b3}{$b_3$} \psfraga <-2pt,0pt> {c1}{$c_1$}
\psfraga <-2pt,0pt> {c2}{$d$} \psfraga <-2pt,0pt> {3}{{III$_{\text{b}}$}}
\psfraga <-2pt,0pt> {2}{{II}} \psfraga <-2pt,0pt> {2'}{{II$^{-1}$}}
\cl{\includegraphics[width=.9\hsize]{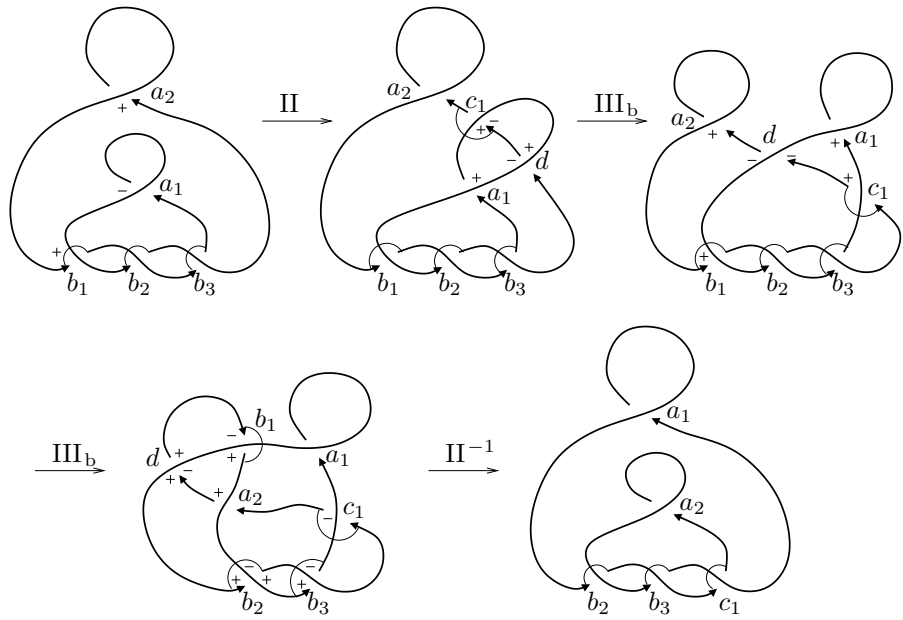}}
\caption{The loop $\Omega_{3,2}$ of trefoils, shown in the
Lagrangian projection}\label{fig:loop}
\end{figure}

\begin{pelda}\label{ex:trefli3}
\figref{fig:loop} shows the succession of four Reidemeister
moves. Using Theorem \ref{thm:linprog} it will be very easy to check
that all four of them are possible; see Example \ref{ex:trefli4}.
Note that we dropped primed labels from our notation which we hope
will not lead to confusion. Let us compute the images of the cycles
$b_1$, $b_2$, and $b_3$ under the composition of the four
holonomies. (We could apply the chain maps to $a_1$ and $a_2$ as
well but because those are not cycles in the chain complex (see
Example \ref{ex:trefli2}), the results in themselves are not very
informative. Let us just mention as an illustration that under the
first of the four moves, the image of $a_2$ is $a_2+d+tb_2b_3d$.)

The first move is of type II. But because after the move, the
diagram still only contains crossings with grading $0$ or $1$
($|c_1|=0$ and $|d|=1$) and $|t|=0$, the boundary of either one of
the index $0$ crossings $b_1$, $b_2$, and $b_3$ is still $0$ after
the move. Thus by Proposition \ref{pro:easy}, all of them are mapped
trivially. The two type III$_{\text{b}}$ moves that follow only
affect $d$ ($d\mapsto d+a_1c_1\mapsto d-b_1a_2+a_1c_1$) but this is
not relevant for our purpose. Finally, the holonomy of the II$^{-1}$
move takes $b_1$ to $-t^{-1}-b_3c_1$ (and $d$ to $0$).
\end{pelda}

Next, we prove that a non-trivial monodromy implies that the loop
itself is non-trivial.

\subsection{Chain homotopies}\label{sec:ch}

It is a folk theorem of singularity theory (see for example 
\cite[p.\ 18]{arnold}, where the more general situation of plane 
curves with semicubical cusp singularities is considered) that there 
are six different codimension $2$ strata (bifurcations) of degenerate 
immersions of $1$--manifolds into $2$--manifolds. These are the three 
different possibilities of two distant Reidemeister moves happening at 
the same time, and three more types: quadruple points with pairwise 
transverse branches; degenerate triple points when exactly two 
branches are tangent; and cubic tangency of two branches. The 
neighborhoods (versal deformations) of points on each of these strata 
are also understood (see \cite[Figure 6]{arnold}). 

For our purposes, each case has to be further separated according to 
the actual position of the branches in $3$--space. Thus, we have six 
types of quadruple points corresponding to the six cyclic orderings of 
the branches. Degenerate triple points have three types distinguished 
by whether the third branch crosses over, under, or between the two 
whose projections are tangent. Cubic tangencies are separated into two 
obvious types. As to the cases of simultaneous Reidemeister moves, let 
us note that if two degenerate Reeb chords appear at the same time, 
their lengths, which we'll call the \emph{height of the move}, can 
generically be assumed different.

\begin{tetel}\label{thm:inv}
Let $\mathscr L\subset\leg$ be a Legendrian knot type. Let $L_0$ and 
$L_1$ be generic elements of $\mathscr L$ and $\Phi(t)$, $\Psi(t)$ two 
generic paths joining them, inducing the holonomies 
$\Phi_*,\Psi_*\colon H(L_0)\to H(L_1)$. If these paths are homotopic 
with their endpoints fixed, then $\Phi_*=\Psi_*$.
\end{tetel}

\begin{proof} 
As $\Phi(t)$ and $\Psi(t)$ are homotopic, there is also a generic 
homotopy $\Phi_s(t)$ connecting them. This map meets the codimension 
$2$ discriminant $\mathscr D_2$ transversely in finitely many points. 
If the number of such points is zero, then the only time the sequence 
of Reidemeister moves changes is those finitely many values of $s$ 
when $\Phi_s(t)$ is tangent to $\mathscr D_1$ (for some value of $t$). 
At these instances, a pair of Reidemeister moves appears or 
disappears. Because these two moves are inverses of each other, their 
holonomies cancel each other out by what we have shown in the proof of 
Theorem \ref{thm:moves}. The rest of the proof consists of an analysis 
of the codimension $2$ strata of the discriminant, as described in the 
introduction of the section. 

{\bf Case 1}\qua  {\sl Simultaneous Reidemeister moves where the two or three 
crossings affected by the moves form disjoint sets (``far away'' 
moves)}\qua Each such case fits the scheme (differentials noted for 
reference later)
\begin{equation}\label{eq:dia}\begin{CD}
\text{before, before }(\partial)@>{\varphi}>>
\text{after, before }(\partial_2)\\
@V{\psi}VV @VV{\hat\psi}V\\
\text{before, after }(\partial_1)@>>{\hat\varphi}>
\text{after, after }(\partial')
\end{CD}\end{equation}
and we can arrange that if one or two pairs of Reidemeister II, 
II$^{-1}$ moves are involved, then the arrows always point in the 
direction of the II$^{-1}$ move. Then it suffices to check that in 
each case, the two compositions of chain maps from upper left to lower 
right agree. Therefore it will not be necessary to introduce a chain 
homotopy, or in other words, it can be chosen to be $0$. This 
phenomenon will re-occur in each case except the very last one. 

By our last assumption and the formulas of Theorem \ref{thm:moves}, 
the claim is clear for generators not directly affected by the moves. 
If any of the moves is a triangle move, it is straightforward to check 
for the vertices of the triangle, too. Suppose the horizontal arrows 
represent II$^{-1}$ moves, with disappearing crossings $a$ and $b$, so 
that $\partial(a)=\pm b+v$. The crossing $a$ gets mapped to $0$ under 
either composition. Let $\partial_1(a)=\pm b+w$. Then 
$\hat\varphi(\psi(b))=\mp w$ and 
$\hat\psi(\varphi(b))=\mp\hat\psi(v)$. But because $\psi(a)=a$ and 
hence $\psi(\partial(a))=\partial_1(a)$, and $\psi(b)=b$, we have 
$\psi(v)=w$. So to show that 
$\hat\varphi(\psi(b))=\hat\psi(\varphi(b))$, we need to prove 
that\footnote{Note that, with a slight abuse of notation, we sometimes 
suppress the chain map relating DGA's of different diagrams when it 
simply replaces all generators in one expression by corresponding 
generators from the other diagram. Our ``excuse'' is that these new 
generators are denoted by the same symbols as the old ones. A more 
exact version of the claim would be that $\psi(v)$ and $\hat\psi(v)$ 
are expressed by identical polynomials in the two DGA's.} 
$\psi(v)=\hat\psi(v)$. This follows by a case-by-case analysis. If 
$\psi$ and $\hat\psi$ are induced by triangle moves then they clearly 
act the same way on $v$. If it is a II$^{-1}$ move with a bigger 
height than that of the one inducing $\varphi$, then 
$\psi(v)=\hat\psi(v)=v$. If the height is lower than that of the first 
move, then we have to make sure that the same expression computes the 
boundary of the higher-index vanishing crossing of the second move in 
each diagram of the upper row of \eqref{eq:dia}, because this is used 
to define $\psi$ and $\hat\psi$. But this is clear since $\varphi$ is 
a chain map and the expressions in question can't contain neither $a$ 
nor $b$.

{\bf Case 2}\qua {\sl Quadruple points}\qua We are going to draw circles linking the 
singularity once and compute monodromies around these loops, starting 
and ending at the top left in each of the cases depicted in Figures 
\fref{fig:quad3} and \fref{fig:quad1}. Such a loop meets the triple 
point discriminant eight times and in each case we will find that the 
composition of the corresponding eight chain maps is the identity. 

\begin{figure}[ht!]\small\anchor{fig:threes}
\psfraga <-2pt,0pt> {a}{$a$}
\psfraga <-2pt,0pt> {b}{$b$}
\psfraga <-2pt,0pt> {c}{$c$}
\psfraga <-2pt,0pt> {d}{$d$}
\psfraga <-2pt,0pt> {a1}{$a_1$}
\psfraga <-2pt,0pt> {a2}{$a_2$}
\psfraga <-2pt,0pt> {b1}{$b_1$}
\psfraga <-2pt,0pt> {b2}{$b_2$}
\psfraga <-2pt,0pt> {x}{$x$}
\psfraga <-2pt,0pt> {K}{$K$}
\psfraga <-2pt,0pt> {p'}{$\psi'$}
\psfraga <-2pt,0pt> {p''}{$\psi''$}
\psfraga <-2pt,0pt> {f'}{$\varphi'$}
\psfraga <-2pt,0pt> {f''}{$\varphi''$}
\psfraga <-2pt,0pt> {ea}{$\varepsilon_a$}
\psfraga <-2pt,0pt> {ea'}{$\varepsilon'_a$}
\psfraga <-2pt,0pt> {ea''}{$\varepsilon''_a$}
\psfraga <-2pt,0pt> {eb}{$\varepsilon_b$}
\psfraga <-2pt,0pt> {eb'}{$\varepsilon'_b$}
\psfraga <-2pt,0pt> {ec}{$\varepsilon_c$}
\psfraga <-2pt,0pt> {ec'}{$\varepsilon'_c$}
\psfraga <-2pt,0pt> {ed}{$\varepsilon_d$}
\psfraga <-2pt,0pt> {ed'}{$\varepsilon'_d$}
\cl{\includegraphics[width=.9\hsize]{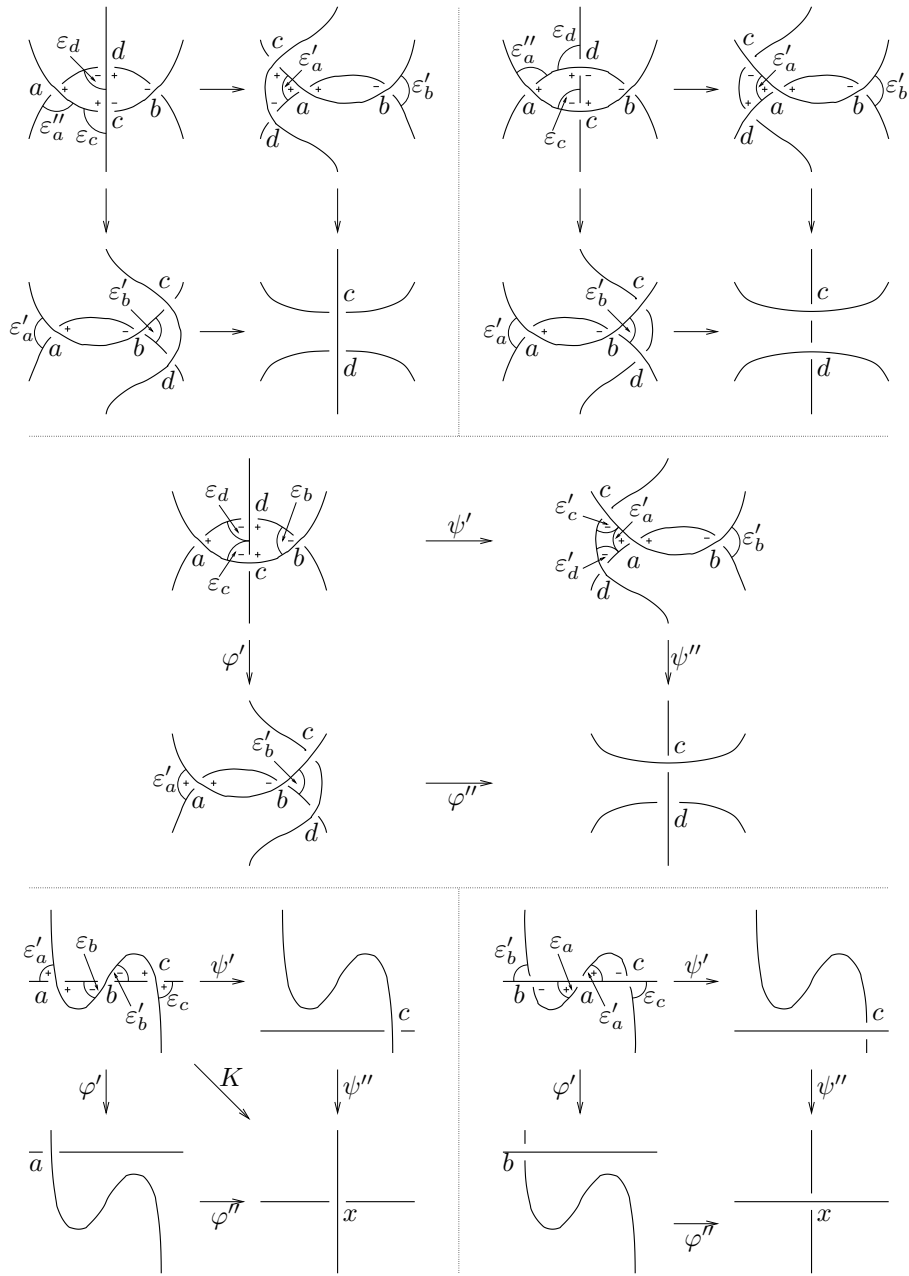}}
\caption{Degenerate triple points and cubic self-tangencies}
\label{fig:threes}
\end{figure}

\begin{figure}[ht!]\small\anchor{fig:quad3}
\psfraga <-2pt,0pt> {a}{$a$}
\psfraga <-2pt,0pt> {b}{$b$}
\psfraga <-2pt,0pt> {c}{$c$}
\psfraga <-2pt,0pt> {d}{$d$}
\psfraga <-2pt,0pt> {e}{$e$}
\psfraga <-2pt,-1.5pt> {f}{$f$}
\psfraga <-2pt,0pt> {ea}{$\varepsilon_a$}
\psfraga <-2pt,0pt> {ea'}{$\varepsilon'_a$}
\psfraga <-2pt,0pt> {eb}{$\varepsilon_b$}
\psfraga <-2pt,0pt> {eb'}{$\varepsilon'_b$}
\psfraga <-2pt,0pt> {ec}{$\varepsilon_c$}
\psfraga <-2pt,0pt> {ec'}{$\varepsilon'_c$}
\psfraga <-2pt,0pt> {ed}{$\varepsilon_d$}
\psfraga <-2pt,0pt> {ed'}{$\varepsilon'_d$}
\psfraga <-2pt,0pt> {ee}{$\varepsilon_e$}
\psfraga <-2pt,0pt> {ee'}{$\varepsilon'_e$}
\psfraga <-2pt,0pt> {ef}{$\varepsilon_f$}
\psfraga <-2pt,0pt> {ef'}{$\varepsilon'_f$}
\cl{\includegraphics[width=.97\hsize]{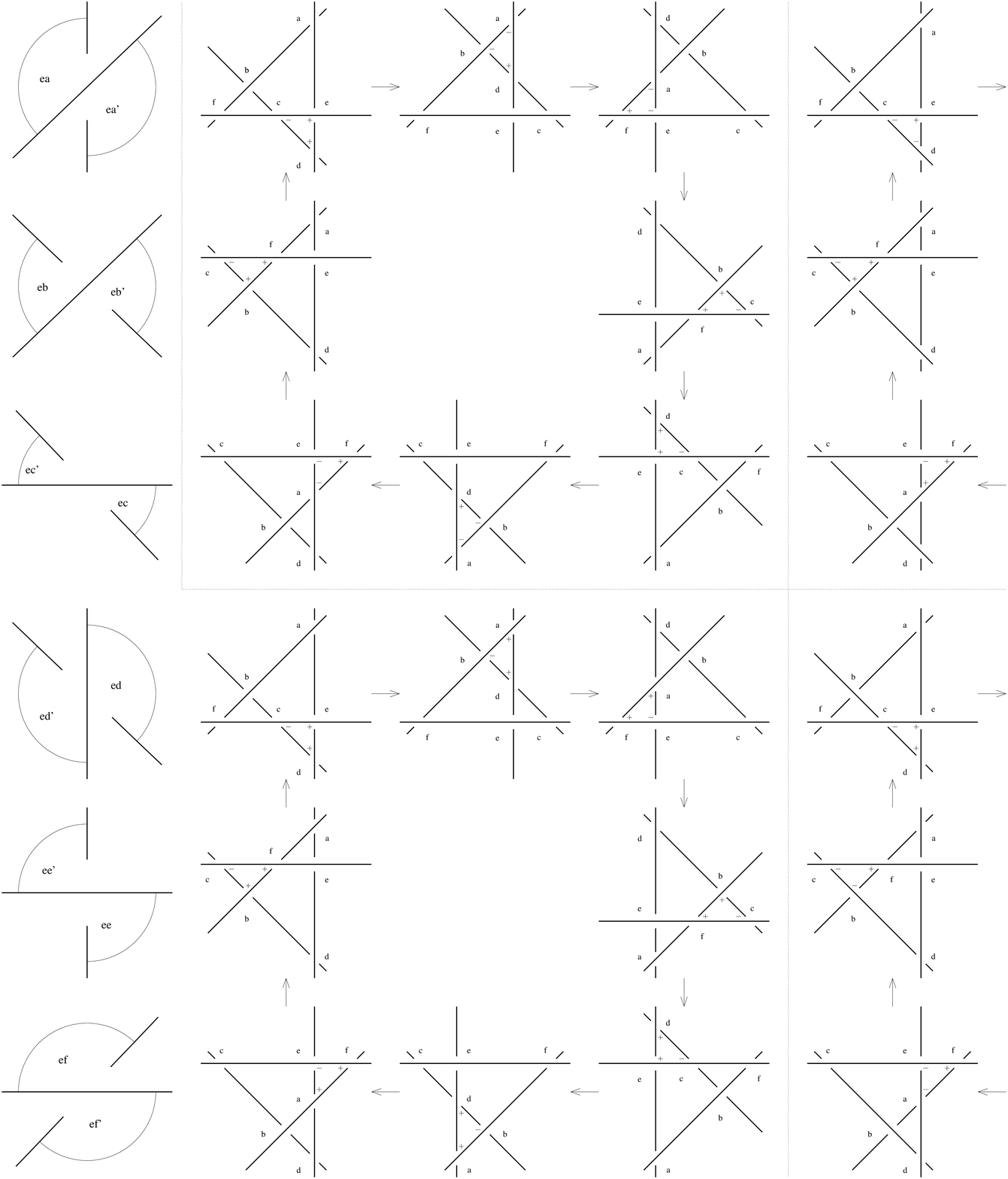}}
\caption{Quadruple points I}\label{fig:quad3}
\end{figure}

\begin{figure}[ht!]\small\anchor{fig:quad1}
\psfraga <-2pt,0pt> {a}{$a$}
\psfraga <-2pt,0pt> {b}{$b$}
\psfraga <-2pt,0pt> {c}{$c$}
\psfraga <-2pt,0pt> {d}{$d$}
\psfraga <-2pt,0pt> {e}{$e$}
\psfraga <-2pt,-1.5pt> {f}{$f$}
\cl{\includegraphics[width=.97\hsize]{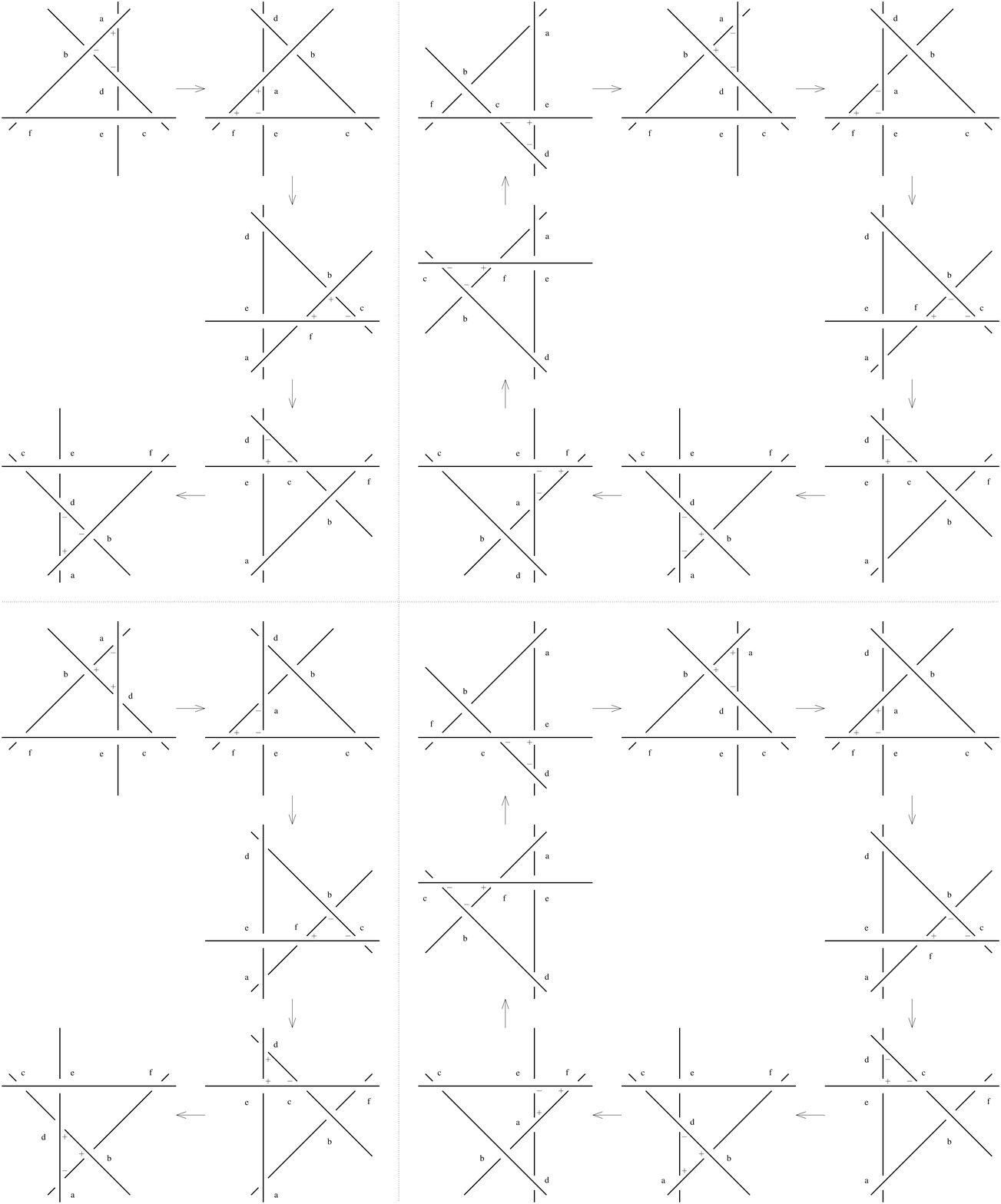}}
\caption{Quadruple points II}\label{fig:quad1}
\end{figure}

Out of the six cases, the one on the bottom of \figref{fig:quad3} 
is the most complicated, due to the fact that all eight moves are of 
type III$_\text{b}$. Even in this case, only the crossings $b$, $c$, 
and $e$ need checking, because others are always mapped trivially. On 
the left side of \figref{fig:quad3}, we indicated the orientation 
signs used in the following. In the second move, $b$'s image is 
$b-\varepsilon'_b\varepsilon_d\varepsilon'_at^{k(b:d,a)}ad$. The three 
generators in this expression are mapped trivially until the sixth 
move, when the image becomes 
$(b-\varepsilon_b\varepsilon'_d\varepsilon_at^{k(b:d,a)}ad)
-\varepsilon'_b\varepsilon_d\varepsilon'_at^{k(b:d,a)}ad=b$ (this is 
by the same argument applied in the proof of Theorem \ref{thm:moves}). 
The generator $e$ returns to itself through two completely analogous 
moves. Finally, we have the following for $c$:
\begin{equation*}\begin{split}
&c\mapsto c-\varepsilon_c\varepsilon_d\varepsilon_et^{k(c;d,e)}ed
\mapsto\text{same}\mapsto 
c-\varepsilon_c\varepsilon_d\varepsilon_et^{k(c;d,e)}
(e-\varepsilon'_e\varepsilon_a\varepsilon_ft^{k(e;a,f)}fa)d\\
&\mspace{-7mu}\mapsto c-\varepsilon'_c\varepsilon_b\varepsilon_ft^{k(c;b,f)}fb
-\varepsilon_c\varepsilon_d\varepsilon_et^{k(c;d,e)}
(e-\varepsilon'_e\varepsilon_a\varepsilon_ft^{k(e;a,f)}fa)d\\
&\mspace{-7mu}\mapsto
c-\varepsilon'_c\varepsilon'_d\varepsilon'_et^{k(c;d,e)}ed
-\varepsilon'_c\varepsilon_b\varepsilon_ft^{k(c;b,f)}fb
-\varepsilon_c\varepsilon_d\varepsilon_et^{k(c;d,e)}
(e-\varepsilon'_e\varepsilon_a\varepsilon_ft^{k(e;a,f)}fa)d\\
&\mspace{-7mu}=c-\varepsilon'_c\varepsilon_b\varepsilon_ft^{k(c;b,f)}fb
-(-1)^{|e|}\varepsilon_c\varepsilon_d\varepsilon_a\varepsilon_f
t^{k(c;d,e)+k(e;a,f)}fad\\
&\mspace{-7mu}\mapsto c-\varepsilon'_c\varepsilon_b\varepsilon_ft^{k(c;b,f)}f
(b-\varepsilon_b\varepsilon'_d\varepsilon_at^{k(b;d,a)}ad)
-(-1)^{|e|}\varepsilon_c\varepsilon_d\varepsilon_a\varepsilon_f
t^{k(c;d,e)+k(e;a,f)}fad\\
&\mspace{-7mu}=c-\varepsilon'_c\varepsilon_b\varepsilon_ft^{k(c;b,f)}fb\mapsto
\text{same}\mapsto c-\varepsilon_c\varepsilon'_b\varepsilon'_ft^{k(c;b,f)}fb
-\varepsilon'_c\varepsilon_b\varepsilon_ft^{k(c;b,f)}fb=c.
\end{split}\end{equation*}
The equality stated in the fourth row, as well as the very last one, 
holds by the argument in the proof of Theorem \ref{thm:moves}. The 
other one in the sixth row is true because 
\begin{itemize}
\item $k(c;b,f)+k(b;d,a)$ and $k(c;d,e)+k(e;a,f)$ both equal the total 
winding of $-\gamma_c+\gamma_d+\gamma_a+\gamma_f$, with each capping 
path thought of as starting and ending at the quadruple point.
\item $\varepsilon'_c\varepsilon_f\varepsilon'_d\varepsilon_a
(-1)^{|e|}\varepsilon_c\varepsilon_d\varepsilon_a\varepsilon_f
=(-1)^{|e|}(-1)^{|c|}(-1)^{|d|}=1$ by Lemma \ref{lem:fokok}.
\end{itemize}

The other five quadruple point bifurcations can be handled with very 
similar, and simpler, computations.

{\bf Case 3}\qua {\sl Degenerate triple points}\qua 
For the three bifurcations on the 
top of \figref{fig:threes}, we claim that the chain maps 
corresponding to the two ways of moving from upper left to lower right 
in the diagrams induce the same chain maps. Please refer to
\figref{fig:threes} for the positions of the orientation signs that 
appear in the calculations and note that $\varepsilon$ and 
$\varepsilon'$ always denote orientation signs of opposite quadrants. 
Consider the third case (middle of \figref{fig:threes}). Denote 
the differential of the upper left diagram by $\partial$, of the lower 
left diagram by $\partial_1$ and that of the upper right diagram by 
$\partial_2$. Applying the chain maps of Theorem \ref{thm:moves} to 
the generator $b$, we find the following:
\[b\overset{\varphi'}\longmapsto b - 
\varepsilon_b\varepsilon_d\varepsilon_c t^{k(b;c,d)} 
cd\overset{\varphi''}\longmapsto \varepsilon'_b\varepsilon'_a v - 
\varepsilon_b\varepsilon_d\varepsilon_c t^{k(b;c,d)} cd\]
$$b\overset{\psi'}\longmapsto b\overset{\psi''}\longmapsto 
\varepsilon'_b\varepsilon'_a
(\varepsilon'_a\varepsilon'_d\varepsilon'_c t^{k(a;d,c)} cd + w).\leqno{\rm and}$$
Here, $v$ and $w$ are defined by the equations 
$\partial_1(a)=\varepsilon_a\varepsilon_b b + v$ and 
$\partial_2(a)=\varepsilon_a\varepsilon_b b + 
\varepsilon'_a\varepsilon'_d\varepsilon'_c t^{k(a;d,c)} cd + w$. The 
coefficients of $cd$ agree in the two images of $b$ by Lemma 
\ref{lem:fokok} and the observation that the capping paths $\gamma_a$ 
and $\gamma_b$ used to compute the exponents of $t$ coincide at the 
moment of self-tangency. Thus, we need to prove that $v=w$. For this, 
note that $(\varphi')^{-1}(a)=a$ and $\psi'(a)=a$, therefore 
$\partial(a)=(\varphi')^{-1}(\partial_1(a))=\varepsilon_a\varepsilon_b
(b-\varepsilon'_b\varepsilon'_d\varepsilon'_c t^{k(b;d,c)} cd) + v$, 
which is also the expression of $\partial_2(a)=\psi'(\partial(a))$. 
Now, $v=w$ follows because 
$(-\varepsilon_a\varepsilon_b\varepsilon'_b\varepsilon'_d
\varepsilon'_c)(\varepsilon'_a\varepsilon'_d\varepsilon'_c)
=-(-1)^{|a|+|b|}=1$. The claim is clear for other generators of the 
upper left diagram. 

In the other two cases (top row of \figref{fig:threes}), the proof 
is very similar. It is important to note though that in the first 
diagram of the top row, $h(d)>h(a)$, and in the second, $h(c)>h(a)$ in 
a neighborhood of the codimension $2$ discriminant. Due to this fact, 
the boundary of $a$ is given by the same polynomial in both pairs of 
upper right and lower left diagrams.

{\bf Case 4}\qua {\sl Cubic self-tangencies}\qua
Note that on the two remaining diagrams 
(bottom row of \figref{fig:threes}), we only have two Reidemeister 
moves; the chain maps $\varphi''$ and $\psi''$ are simple 
re-labelings. In the case depicted on the second diagram, we still 
don't need to introduce a non-zero chain homotopy. This is because 
$\varphi''(\varphi'(a))= \psi''(\psi'(a))=0$, 
$\psi''(\psi'(b))=\psi''(\varepsilon'_a\varepsilon'_b
(\varepsilon'_a\varepsilon'_c c))=x=\varphi''(\varphi'(b))$ (it is an 
easy fact that $\varepsilon'_b=\varepsilon'_c$), and a similar 
computation for $c$. 

Finally, it is the situation on the bottom left of
\figref{fig:threes} which does require a non-trivial chain homotopy $K$, 
mapping between the two indicated DGA's. Let 
$K(b)=\varepsilon_a\varepsilon_b x$ and let $K$ map all other 
generators to $0$. We extend $K$ as in Lemma \ref{lem:K}. Let 
$\varphi=\varphi''\circ\varphi'$ and $\psi=\psi''\circ\psi'$. Then 
$\varphi(a)-\psi(a)=x$, $\varphi(c)-\psi(c)=-x$ and 
$\varphi(b)-\psi(b)=
\varepsilon_b\varepsilon_cv_c-\varepsilon'_a\varepsilon'_bv_a$, where 
the polynomials $v_a$ and $v_c$ are defined in the upper left diagram 
by the equations $\partial(a)=\varepsilon_a\varepsilon_bb+v_a$ and 
$\partial(c)=\varepsilon'_c\varepsilon'_bb+v_c$. Finally, it is easy 
to check the following, too:
\[(K\circ\partial+\partial'\circ K)(a)
=K(\varepsilon_a\varepsilon_bb+v_a)=x\]
(because $v_a$ doesn't contain $b$);
\[(K\circ\partial+\partial'\circ K)(c)
=K(\varepsilon'_c\varepsilon'_bb+v_c)=\varepsilon'_c\varepsilon'_b
\varepsilon_a\varepsilon_bx=(-1)^{|a|+|b|}x=-x\]
(because $\varepsilon'_c=\varepsilon'_a$), and 
\[(K\circ\partial+\partial'\circ K)(b)
=\partial'(\varepsilon_a\varepsilon_b 
x)=\varepsilon_a\varepsilon_b v_a + \varepsilon_a\varepsilon_b v_c,\]
where, indeed, $\varepsilon_a=\varepsilon_c$ and 
$\varepsilon_a\varepsilon_b=-\varepsilon'_a\varepsilon'_b$.
\end{proof}

\subsection[Z_2-coefficients]{\boldmath $\Z_2$--coefficients}\label{sec:z2}

In this subsection we briefly summarize the modifications needed to
reduce our discussion to the original $\Z_2$--coefficient theory of
Chekanov. (We will work in that context in sections \ref{sec:BC}
through \ref{sec:moreaug}.) Please note again that this means that
we are going to describe the same geometry (Legendrian knots, their
families, Reeb chords, and holomorphic discs) using slightly less
sophisticated algebra: The generators are still the same Reeb
chords, ie\ crossings, and the product is non-commutative, but we
substitute $t=1$, reduce integers modulo $2$, and reduce the grading
modulo $2r$, where $2r$ is the common Maslov number of each knot in
the family. The only simplification we get in Table \ref{algebra} is
that now, $(-1)^{|a|}=1$ (and $\varphi-\psi=\varphi+\psi$).
Otherwise, we keep all notation introduced in this section and still
refer to \figref{fig:moves} (orientation signs can now be
ignored).

A list of maps that, when extended from the generators to the DGA
as algebra morphisms, become the chain maps that are used to
define holonomies and monodromies of (sequences of) Reidemeister
moves, is as follows.

Move III$_{\text{a}}$: Let $a\mapsto a'$, $b\mapsto b'$, $c\mapsto
c'$, and $x\mapsto x'$, where $x$ is any other crossing of the
upper diagram.

Move III$_{\text{b}}$: Let
\[a\mapsto a'+c'b',\]
while other generators are mapped trivially: $b\mapsto b'$,
$c\mapsto c'$, and $x\mapsto x'$.

Move II$^{-1}$: Let $\partial(a)=b+v$. Define $x\mapsto x'$, which
gives rise to the obvious re-labeling $v\mapsto v'$. Then let
\begin{eqnarray*}
a&\mapsto&0\\
b&\mapsto&v'.
\end{eqnarray*}
Move II: The map $\varphi\colon\mathscr A'\to\mathscr A$ still
takes
\begin{equation*}
b'_j\mapsto\varphi(b_j')=b_j\end{equation*} for all
$j=1,\ldots,m$. Formula \eqref{eq:undor} of Remark \ref{kiszoroz}
becomes
\begin{equation*}
\begin{split}
a_1'\mapsto\varphi(a'_1)=a_1+\sum\Big(&
B_1aB_2bB_3bB_4b\ldots B_kbA\\
+&B_1vB_2aB_3bB_4b\ldots B_kbA\\
+&B_1vB_2vB_3aB_4b\ldots B_kbA\\
+&B_1vB_2vB_3vB_4a\ldots B_kbA\\
+&\ldots\\
+&B_1vB_2vB_3vB_4v\ldots B_kaA\Big),
\end{split}\end{equation*}
and formula \eqref{megundor} takes the form
\begin{equation*}
\begin{split}
a_i'\mapsto\varphi(a_i')=a_i+\sum\Big(&
\bar{B_1}aB_2bB_3bB_4b\ldots B_kbA\\
+&\bar{B_1}v\bar{B_2}aB_3bB_4b\ldots B_kbA\\
+&\bar{B_1}v\bar{B_2}v\bar{B_3}aB_4b\ldots B_kbA\\
+&\bar{B_1}v\bar{B_2}v\bar{B_3}v\bar{B_4}a\ldots B_kbA\\
+&\ldots\\
+&\bar{B_1}v\bar{B_2}v\bar{B_3}v\bar{B_4}v\ldots\bar{B_k}aA\Big).
\end{split}\end{equation*}
In the rest of the paper, each time we compute the holonomy of a
Reidemeister II move, Proposition \ref{pro:easy} applies. In other
words, the detailed description of move II was included only for
completeness.

\section{Legendrian Reidemeister moves}\label{sec:simplex}

In this section, we investigate a very basic question, which arises
naturally in the theory of Legendrian knots and links. The question is
this: if a Lagrangian diagram of a Legendrian link contains a part
that is \emph{topologically} fit for a Reidemeister II,
II$^{-1}$, or III move, can we carry the move out in the Legendrian
category? (If so, we will call the move \emph{consistent}.) As it
turns out, the answer is not always affirmative and a necessary and
sufficient condition can be formulated in terms of a linear program.

We call two regions (of the complement of a knot diagram)
\emph{adjacent at a vertex} if they share the vertex (crossing) but
don't share any edges.

\begin{tetel}\label{thm:linprog}
Each Reidemeister move of \figref{fig:moves} is possible if and
only if the Lagrangian diagram can be isotoped (in the Legendrian
sense) so that a certain inequality is satisfied. The list is as
follows:
\begin{itemize}
\item Move III$_{\text{a}}$:
the area of the triangle is smaller than the area of any of the three
regions adjacent to it at either of its vertices.
\item Move III$_{\text{b}}$:
the area of the triangle is smaller than the area of the region
adjacent to it at the vertex $a$.
\item Move II$^{-1}$:
the area of the $2$--gon is smaller than the sum of the areas of the
two regions adjacent to it at its vertices.
\item Move II:
the sum of the heights of the crossings along the right side of the
pinching region, counted with Reeb signs, is positive. (If the region
to the right of the pinching, as it is depicted in 
\figref{fig:moves}, is the 
unbounded one, then there is no
obstruction to the move.)
\end{itemize}
\end{tetel}

Recall from section \ref{sec:prep} that every Lagrangian diagram
$\gamma$ has the non-empty cone $\mathscr C_\gamma$ associated to
it. $\mathscr C_\gamma$ is defined by the positivity constraints
on the heights associated with $\gamma$ and finitely many
homogeneous linear inequalities (coming from bounded regions of
the complement) in terms of these heights (see Remark
\ref{rem:kup}). Each inequality of Theorem \ref{thm:linprog}
requires that a certain linear objective function\footnote{The
objective function is linear exactly because the areas of the
regions are linear functions of the heights. These areas can also
be thought of as slack variables.} defined on $\mathscr C_\gamma$
take on a positive value. That is, in each case the necessary and
sufficient condition can be re-phrased to require that a certain
linear programming problem be unbounded. This can be checked by
the simplex method or any other linear programming algorithm (for
small diagrams even by hand, as the author can attest).

\begin{proof}
The ``only if'' part of each statement is easy to justify. To see this
for move II, use Proposition \ref{teruletes}: at the moment of
self-tangency of the projection, the sum in question is the area of
the region pinching off on the right side plus the difference in
$z$--coordinate between the two preimages of the self-tangency. As
this is positive, our sum must have been positive for some period of
time before the move happened, too. The other three moves involve the
vanishing of the area of a region; clearly, for some time before that
happens, the particular area is smallest among all others and the said
inequalities hold.

Now to prove the ``if'' statement for move II, suppose that the sum
$S$ described in the statement is positive and extend a very narrow
`finger' that follows closely the right side of the region from the
lower edge to the upper one. The area of the finger can be made
arbitrarily small and therefore with a slight downward bump on the
lower edge (see Remark \ref{lokalis}) we can make sure that other
parts of the knot (in particular the value of $S$) are not affected.
As we are also able to arrange that the area that is pinching off be
less than $S$, an application of Proposition \ref{teruletes} shows
that at the moment the finger reaches the upper edge, it must do so so
that it crosses underneath it. If the `right side' that we used
doesn't exist due to the unboundedness of the pinching region, then
simply let the finger follow a path that encircles a large enough area
(on the finger's left) so that it dips down sufficiently.

\begin{figure}[ht!]\small\anchor{fig:btriangle}
\psfraga <-4pt,-4pt> {3b}{III$_{\text{b}}$}
\psfrag{a}{$a$}
\psfraga <-2pt,-2pt> {b}{$b$}
\psfraga <-2pt,-2pt> {c}{$c$}
\psfraga <-4pt,9pt> {a'}{\tiny$a'$}
\psfraga <-9pt,6pt> {b'}{\tiny$b'$}
\psfrag{c'}{\tiny$c'$}
\psfrag{V}{$V$}
\cl{\includegraphics[width=.7\hsize]{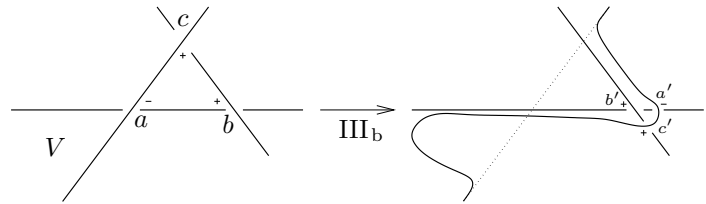}}
\caption{Reidemeister III$_{\text{b}}$ move}\label{fig:btriangle}
\end{figure}

In the other three cases, we again utilize the idea described in
Remark \ref{lokalis}. To carry out move III$_{\text{b}}$, replace the
uppermost strand as indicated in \figref{fig:btriangle}. It is
clear from the proof of Proposition \ref{teruletes} that if the new
arc, starting from the old crossing $c$, follows the side $cb$ close
enough, then the new crossing $c'$ has height arbitrarily close to
$h(c)$, and $h(a')$ is arbitrarily close to $h(c)+h(b)$. In
particular, the types of the crossings are as indicated. Then, if we
followed the sides $cb$ and $ba$ close enough, we can use the excess
area of the region $V$ to compensate for the modification.

Three very similar arguments establish the claim for move
III$_{\text{a}}$, and \figref{fig:trukk} is self-explanatory for
the case of move II$^{-1}$.\end{proof}

\begin{figure}[ht!]\small\anchor{fig:trukk}
\psfraga <-4pt,-4pt> {1}{II$^{-1}$}
\cl{\includegraphics[width=.7\hsize]{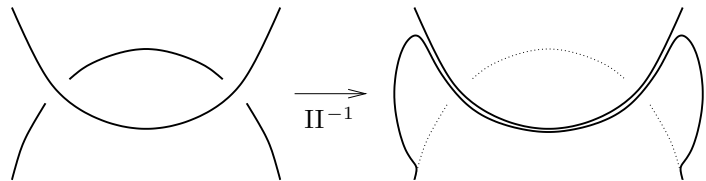}}
\caption{Reidemeister II$^{-1}$ move}\label{fig:trukk}
\end{figure}

\begin{pelda}\label{ex:trefli4}
We'll prove that the loop of trefoils in Example \ref{ex:trefli3} is
sound, ie,\ all four of the indicated Reidemeister moves can indeed
be carried out. Let us consider the first diagram in
\figref{fig:loop}. By adding ``bulges'' of equal area to the regions that
were denoted by $U_1$ and $U_2$ in \figref{fig:trefli}, we may
arrange that the areas of those, along with the height $h(a_2)$, be
arbitrarily larger than any other area or height occurring in the
projection. In particular, $h(a_2)+h(b_1)-h(a_1)$ can be made
positive, which by Theorem \ref{thm:linprog} means that the first move
is possible. The triangle move that follows is consistent because the
region adjacent to the triangle at $a'$ is the unbounded one. At this
point, we may finish the argument by pointing out that the remaining
two moves are mirror images of the first two, or we may continue as
follows: For the other triangle move, enlarge the areas of the
triangle $c_1a_1a'$ and of the region capping off $a_1$ until the area
of the triangle gets bigger than that of the triangle $a_2b_1a'$ (this
is meant in the third diagram of \figref{fig:loop}). Finally, the
same trick applied to the bounded regions adjacent to $a_1$ in the
fourth diagram makes sure that the appropriate condition from Theorem
\ref{thm:linprog} is satisfied and the last move of the loop is
consistent, too.
\end{pelda}

\section{First evidence of non-triviality}\label{sec:trefli}

At this point we are ready to prove that the monodromy invariant
of Theorem \ref{thm:hurkok} is non-trivial. We will use
essentially the same example to establish Corollary \ref{kov},
too.

Note that the first and last diagrams of \figref{fig:loop}
coincide, so in Examples \ref{ex:trefli3} and \ref{ex:trefli4} we
investigated a loop of Legendrian trefoil knots. It is not hard
to identify it with $\Omega_{3,2}$ (defined in the introduction). If
we re-label the fifth diagram like the first one, we can write that
the monodromy of the loop (which is the composition of the four
holonomies of Example \ref{ex:trefli3} and the re-labeling) acts on
the cycles $b_1$, $b_2$, and $b_3$ as follows:
\begin{eqnarray*}
b_1&\mapsto&-t^{-1}-b_2b_3\\
b_2&\mapsto&b_1\\
b_3&\mapsto&b_2.
\end{eqnarray*}

\begin{all}\label{pro:5}
The map $\mu_0$ defined by the formulas above on the index $0$
part of the contact homology (of the knot introduced in Example
\ref{ex:trefli1}) has order $5$. Thus Theorem \ref{thm:monorendje}
holds in the case $p=3$, $q=2$.
\end{all}

\begin{proof}
As $\mu_0(b_3)=b_2$ and $\mu_0^2(b_3)=b_1$, all three generators
are on the same orbit. Then we compute that
$\mu_0^4(b_3)=\mu_0(-t^{-1}-b_2b_3)=-t^{-1}-b_1b_2$ and that
\[\mu_0^5(b_3)=\mu_0(-t^{-1}-b_1b_2)=-t^{-1}-(-t^{-1}-b_2b_3)b_1
=-t^{-1}+t^{-1}b_1+b_2b_3b_1=b_3\] by equation \eqref{trick1} of
Example \ref{ex:trefli2}.

The argument for the order will be complete if we show that
$b_3\ne b_2$ in the homology, ie,\ that there isn't an element
$a\in\mathscr A$ so that $\partial(a)=b_3-b_2$. For this, notice
that $\im\partial$ is contained in the two-sided ideal generated
by $\partial(a_1)$ and $\partial(a_2)$. An examination of
equations \eqref{eq:del1} and \eqref{eq:del2} shows that after
substituting $t=-1$, they both contain an even number of terms
that are powers of $b_3$ (namely, two and zero terms). Therefore
the same property holds for the entire ideal. Since $b_3-b_2$ is
not such a polynomial, it cannot be in the image of the
differential.
\end{proof}

It is easy to see that if we fix the basepoint, but allow knots in
our one parameter families that are not Legendrian, then the
three-fold concatenation of $\Omega_{3,2}$ is homotopic to a
continuous rotation by $2\pi$ radians. Therefore, as
$\pi_1(SO(3))=\Z_2$ is generated by exactly this loop, the six-fold
concatenation of $\Omega_{3,2}$ is contractible in the space
$\mathscr K$ of smooth right-handed trefoil knots. 
But because
$\mu((\Omega_{3,2})^6)\ne\id$ (as
$\mu_0((\Omega_{3,2})^6)=\mu_0(\Omega_{3,2})\ne\id$),
$(\Omega_{3,2})^6$ is non-contractible in the space $\mathscr L$ of
Legendrian trefoil knots. This last implication, which establishes
Corollary \ref{kov}, holds by Theorem \ref{thm:inv}. (We remark that 
with a bit of extra work, it can be shown that $\Omega_{3,2}$ is itself 
contractible in $\mathscr K$, so taking its sixth power isn't really 
necessary in our argument.)

In sections \ref{sec:moreBC} and \ref{sec:moreaug}, we will
generalize the proof of Proposition \ref{pro:5} for an arbitrary $(p,q)$ 
torus knot.
(Except that to simplify the discussion, in those sections we will
work with $\Z_2$--coefficients.) Before we do that, we use sections
\ref{sec:BC} and \ref{sec:loop} to set up a more general picture.
Section \ref{sec:aug} contains a vital ingredient of the proof of
Theorem \ref{thm:monorendje}.

\section{Positive braid closures}\label{sec:BC}

Let us consider a positive braid $\beta$ on $q$ strands, as in
\figref{fig:fonat}. Label the left and right endpoints of the
strands from top to bottom with the first $q$ whole numbers. The
pair of left and right labels on each strand takes the form
$(i,\sigma(i))$, where $\sigma$ is the \emph{underlying permutation}
of $\beta$. Further, label the crossings of the braid with a pair of
numbers, the first one the left label $i$ of the overcrossing strand
and the second one the right label $j$ of the undercrossing one.
(Note that not all pairs of numbers between $1$ and $q$ occur as
labels of crossings; for example, $(i,\sigma(i))$ never does for any
$i$.)

\begin{figure}[ht!]\small\anchor{fig:pic}
\psfrag{a1}{$a_1$}
\psfrag{a2}{$a_2$}
\psfrag{a3}{$a_3$}
\psfrag{a4}{$a_4$}
\psfrag{t1}{$T_1$}
\psfrag{t2}{$T_2$}
\psfrag{t3}{$T_3$}
\psfrag{t4}{$T_4$}
\psfrag{u1}{$U_1$}
\psfrag{u2}{$U_2$}
\psfrag{u3}{$U_3$}
\psfrag{u4}{$U_4$}
\psfrag{x1}{$\alpha_1$}
\psfrag{x2}{$\alpha_2$}
\psfrag{x3}{$\alpha_3$}
\psfrag{x4}{$\alpha_4$}
\psfrag{p1}{$p_1$}
\psfrag{p2}{$p_2$}
\psfrag{p3}{$p_3$}
\psfrag{p4}{$p_4$}
\cl{\includegraphics[width=.8\hsize]{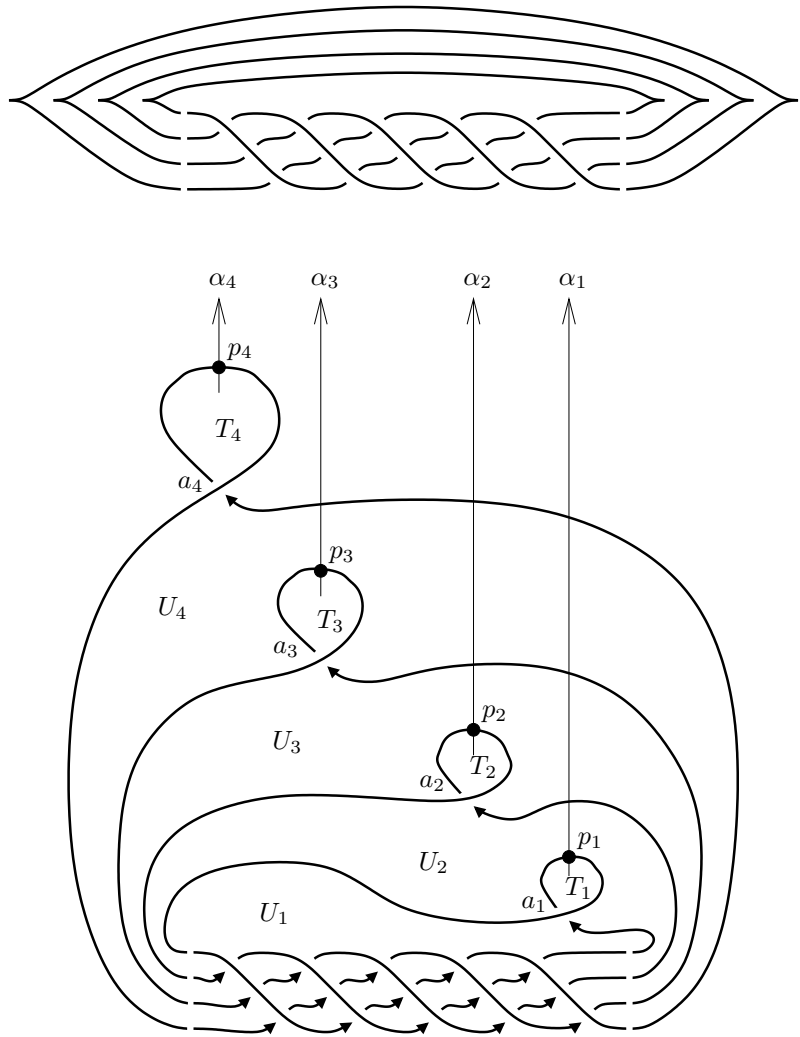}}
\caption{Front and Lagrangian diagrams of the closure of a positive
braid}\label{fig:pic}
\end{figure}

\begin{Def}\label{def:lezaras}
A positive braid $\beta$ defines the front diagram of an oriented
Legendrian link as on the upper half of \figref{fig:pic}
(orient each strand of $\beta$ from left to right). We call
this construction the \emph{Legendrian closure} of $\beta$ and denote it by
$L_\beta$.
\end{Def}

Note that the types of
crossings are correctly determined by the slopes of the branches that meet
there exactly because the braid is positive.
Applying Ng's \cite{computable} construction of
resolution to the front diagram, we obtain the Lagrangian diagram
$\gamma_\beta$ in the lower half of \figref{fig:pic}.

In section \ref{sec:prep}, we omitted the modifications needed to
define the invariants (such as rotation and Thurston--Bennequin
numbers and contact homology itself) for oriented multi-component
Legendrian links. We include an informal rundown here (see
\cite[section 2.5]{computable} for more). The rotation number is the
sum of the rotations of the components. The Thurston--Bennequin
number is still the writhe of the Lagrangian projection; note that
now, as opposed to the case of a single component, the orientation
matters in its definition.

From now on, we will always use $\Z_2$--coefficients in contact
homology. For a multi-component link, there is not a single
distinguished grading of its DGA, rather a family of so called
admissible gradings. We consider those introduced in
\cite{computable} and not the larger class of gradings described in
\cite[section 9.1]{chek}. In the $\Z_2$--coefficient theory, each of
these is defined modulo the greatest common divisor of the Maslov
numbers of the components (recall that the Maslov number is twice
the rotation number of a knot). Self-crossings of individual
components have the same well-defined index in any admissible
grading. We will refer to such generators as \emph{proper
crossings}. Other crossings' indices have the same parity in each
admissible grading. The parity of the index of any crossing
coincides with its sign in classical knot theory.

The differential $\partial$ is of index
$-1$ with respect to each admissible grading. (Furthermore, Lemma 
\ref{lem:fokok} holds for each admissible grading, too.) Contact homology 
is 
a well-defined invariant in the sense that if two diagrams are Legendrian
isotopic, then there is a one-to-one correspondence between their sets
of admissible gradings so that the corresponding contact homologies are
isomorphic as graded algebras.

It is easy to calculate that the Legendrian closure of any positive
braid has rotation number $r=0$ (in fact, every component has Maslov
number $0$). The Thurston--Bennequin number is
$tb(L_\beta)=(\text{word length of }\beta)-q$. The positive
Legendrian $(p,q)$ torus knot obtained as a special case is the one
with maximal Thurston--Bennequin number $p(q-1)-q$ (see \cite{EH1}
for a classification of Legendrian torus knots). 
The $(3,2)$ torus knot that we considered
earlier arises from this construction, too. The unique
Legendrian unknot with maximal Thurston--Bennequin number $tb=-1$
\cite{benn,unknot} is the Legendrian closure of the trivial braid on
a single strand.

\begin{figure}[ht!]\small\anchor{fig:indexek}
\psfraga <-1.5pt,-1.5pt> {bij1}{$b_{i,j,1}$}
\psfrag{bij2}{$b_{i,j,2}$}
\psfraga <-1.5pt,-1.5pt> {bij3}{$b_{i,j,3}$}
\psfraga <-2pt,-2pt> {bpq1}{$b_{p,r,1}$}
\psfraga <-1.5pt,-1.5pt> {bpq2}{$b_{p,r,2}$}
\psfrag{i}{$i$}
\psfrag{j}{$j$}
\psfrag{p}{$p$}
\psfrag{q}{$r$}
\cl{\includegraphics[width=.8\hsize]{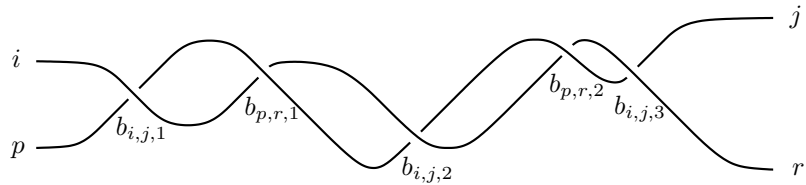}}
\caption{Labeling the crossings of a braid}\label{fig:indexek}
\end{figure}

In the Chekanov--Eliashberg DGA of $\gamma_\beta$, the grading
(any admissible grading) is integer-valued. The crossings $a_m$
($m=1,\ldots,q$) have index $1$. The rest of the generators are
the crossings of $\beta$ and they have index $0$. In the
multi-component case we should say instead that the grading that
assigns the index $0$ to each is admissible; from now on, we will
always work with this grading. These generators will be labeled by
$b_{i,j,t}$, where the integers $1\le i,j\le q$ are the ones
defined at the beginning of the section. The third label $t$ is
used to distinguish between multiple intersections of strands,
enumerating them from left to right, as shown in
\figref{fig:indexek}.

Since $\partial$ lowers the index by $1$,
\[\partial(b_{i,j,t})=0\text{ for all }i,j=1,\ldots,q\text{ and }t.\]
In fact, there are no admissible discs whose positive corner is an
index $0$ crossing. The boundaries $\partial(a_m)$,
$m=1,\ldots,q$, are polynomials in the non-commuting variables
$b_{i,j,t}$ only. Our goal for the remainder of this section is to
compute these polynomials.

The first observation is that for all $m=1,\ldots,q$, there is an
admissible disc covering the teardrop-shaped region $T_m$ right above $a_m$
once. We'll call it the \emph{$m^{\rm th}$ trivial disc}. It contributes $1$ to
$\partial(a_m)$. The $m^{\rm th}$ trivial disc is the only one that turns at the
upward-facing positive quadrant at $a_m$. The rest of the contributions to
$\partial(a_m)$ come from discs that turn at the positive quadrant facing
down.

Suppose $f\colon\Pi_k\to\R^2_{xy}$ is an admissible immersion with respect
to the projection $\gamma=\gamma_\beta$, with
positive corner $f(x_0^k)=a_m$ and
so that it is different from the $m^{\rm th}$ trivial disc. Fix
points $p_1,\ldots,p_q$ on $\gamma$ as shown in \figref{fig:pic}.

\begin{lemma}\label{lem:kicsi}
The curve $f(\partial\Pi_k)$ doesn't pass through the points
$p_m,\ldots,p_q$.
\end{lemma}

\begin{proof}
In fact, if we denote the region of the complement directly under
$a_i$ by $U_i$, then $f(\Pi_k)$ is disjoint from
$U_{m+1},\ldots,U_q$, as well as from $T_m,\ldots,T_q$. To see
this, first shrink $\gamma$ so that all areas and heights are
smaller than $\varepsilon^{q-m}$, where $\varepsilon>0$ is to be
chosen later. Then add ``bulges'' of equal area $1$ to $T_q$ and
$U_q$ on the outside of the diagram (see Remark \ref{lokalis}).
Next, add bulges of area $\varepsilon$ to $T_{q-1}$ and $U_{q-1}$
on their sides facing $U_q$. Continue all the way until adding
bulges of area $\varepsilon^{q-m-1}$ to $T_{m+1}$ and $U_{m+1}$ on
their sides facing $U_{m+2}$. The result is a Lagrangian
projection of the same link in which, after the choice of a small
enough $\varepsilon$, the claim is obvious by Lemma
\ref{lem:ugras}, except for the case of $T_m$. But that follows
from the easy observation that if any admissible disc passes
through $p_m$ so that locally, $f(\Pi_k)$ faces downward, then $f$
is the $m^{\rm th}$ trivial disc.
\end{proof}

Recall that $\gamma$ is oriented as shown in \figref{fig:pic}.
We have already noted that apart from its positive corner which is
of index $1$, $f$ turns only at certain index $0$ crossings
$b_{i,j,t}$. By Lemma \ref{lem:kovet}, this implies that $f$ is
compatible with the orientation of $\gamma$. Therefore, as the
orientation of $f(\partial\Pi_k)$ agrees with that of $\gamma$
near the positive corner $a_m$, these orientations agree at all
points of $\partial\Pi_k$. This last observation for example
implies that $f$ can't turn at the downward facing negative
quadrant at any of the $b_{i,j,t}$'s.

We may summarize our findings about $f$ so far as follows. The curve
$f(\partial\Pi_k)$
starts at $a_m$, follows $\gamma$ until it reaches the braid at the left
endpoint labeled $m$ (while $f$ extends this map toward $U_m$).
Then it travels through the braid,
possibly turning (left) at several crossings $b_{i,j,t}$ but always
heading
to the right, until it reaches a right endpoint labeled $i_1\le m$. Then
$f(\partial\Pi_k)$ climbs to $a_{i_1}$ and, unless $i_1=m$, beyond
$a_{i_1}$ (note that $f$ can't have a positive corner at $a_{i_1}$)
to $p_{i_1}$, at which point the extension $f$ is toward $U_{i_1+1}$. Then
$f(\partial\Pi_k)$ descends back to $a_{i_1}$ so that it doesn't turn
at the negative corner. So the process repeats with $(m,i_1)$
replaced with $(i_1,i_2)$ (where $i_2\le m$),
and so on until $i_{c+1}=m$ for some $c$.

\begin{Def}\label{def:sorozat}
A finite sequence of positive integers is called \emph{admissible} if for
all $s\ge0$, between any two appearances of $s$ in the sequence there is a
number greater than $s$ which appears between them. For $n\ge1$, let us
denote by $D_n$
the set of all admissible sequences that are composed of the numbers
$1,2,\ldots,n-1$.
\end{Def}

For example, $D_1=\{\:\varnothing\:\}$,
$D_2=\{\:\varnothing,\{\,1\,\}\:\}$, and
\[D_3=\{\:\varnothing,
\{\,1\,\},\{\,2\,\},\{\,1,2\,\},\{\,2,1\,\},\{\,1,2,1\,\}\:\}.\]
By induction on $n$
and observing the position of the unique maximal term in the sequence, it
is easy to prove that $|D_n|=|D_{n-1}|^2+|D_{n-1}|$.

\begin{all}\label{pro:blank}
For the admissible disc $f$ as above, the
sequence of intermediate labels
$\{\,i_1,i_2,\ldots,i_c\,\}$ is an element of $D_m$.
\end{all}

\begin{proof}
There is an intuitive reason for the claim: if the index $i$ was repeated
without the boundary of the disc climbing higher between the two
occurrences, then the disk would pinch off at $a_i$. We shall give a more
rigorous proof using Blank's theorem
(see \cite{blank}, or the review
in \cite{franc}, from where we'll borrow our terminology).

Let $\tilde f$ be a small generic perturbation of $f$.
This has the corners rounded so that $\tilde f\big|_{\partial\Pi_k}$ is an
immersion; moreover, let us note that
because $f(\partial\Pi_k)$ may cover parts of $\gamma$ more than once,
the complement of $\tilde f(\partial\Pi_k)$ typically contains more
regions than that of $f(\partial\Pi_k)$.
We are going to apply Blank's theorem to
$\tilde f\big|_{\partial\Pi_k}$.
We extend rays from each bounded component of
$\R^2\setminus \tilde f(\partial\Pi_k)$ to infinity.
If the region is not one of those obtained from the $T_i$
(let us denote these by $T'_i$), or one of the
small ones resulting from perturbation
close to the arc bounding $T_i$, this can be done so that $\gamma$
only intersects the ray in a positive manner (ie, if the ray is
oriented toward $\infty$, then $\gamma$ crosses from the right side to
the left side). By our observation on the compatibility of $f$, this
implies that $\tilde f(\partial\Pi_k)$ also intersects those rays positively.
Let us draw rays starting from the regions $T'_1,\ldots,T'_q$ as in 
\figref{fig:pic}. Label all rays, in particular label these last $q$ with
the symbols $\alpha_1,\ldots,\alpha_q$.

The \emph{Blank word} of the disc $\tilde f$ is obtained by tracing
$\tilde f(\partial\Pi_k)$, starting from, say, $a_m$, and writing down the
labels of
the rays we meet with exponents $\pm1$
according to positive and negative
intersections (we'll call these
\emph{positive} and \emph{negative symbols} or \emph{letters}). A
\emph{grouping} of the Blank word is a set of
properly nested disjoint unordered pairs of
the form $\{\,\alpha,\alpha^{-1}\,\}$ so that each negative symbol
is part of exactly one pair.
Blank's theorem states that the set of groupings of the Blank word is in
a one-to-one correspondence with non-equivalent extensions of the
immersion $\tilde f\big|_{\partial\Pi_k}$ to immersions of $\Pi_k$.

As $\tilde f$ is such an extension, a grouping exists and we claim that this
implies the Proposition. We delete labels different from the $\alpha_s$
from the Blank word of $\tilde f$ and concentrate on
the remaining word, which also inherits a grouping (note that most of the
deleted labels only appeared positively anyway, except for certain ones
that belong to some regions that were the result of perturbation).
We'll call this the \emph{Blank word of $f$} and
denote it by $W_f$. Note that if we only keep the negative letters
$\alpha_s^{-1}$ from $W_f$, then the sequence of their indices is exactly
$\{\,i_1,i_2,\ldots,i_c\,\}$. In fact, $W_f$ is decomposed into segments
$S_1,\ldots,S_c$ ended by these negative symbols, and a
final segment $S$:
\begin{multline*}
W_f=
\overbrace{\alpha_1\alpha_2\ldots\alpha_{i_1-1}\alpha_{i_1}^{-1}}^{S_1}
\cdot
\overbrace{\alpha_1\alpha_2\ldots\alpha_{i_2-1}\alpha_{i_2}^{-1}}^{S_2}
\cdot\ldots\\
\cdot\underbrace{\alpha_1\alpha_2\ldots\alpha_{i_c-1}\alpha_{i_c}^{-1}}
_{S_c}\cdot
\underbrace{\alpha_1\alpha_2\ldots\alpha_{m-1}}_S.
\end{multline*}

We'll prove the following statements by induction on $j$:
\begin{enumerate}[\rm(1)]
\item\label{rend}
Each copy of $\alpha^{-1}_{m-j}$ is paired in the grouping with a copy
of $\alpha_{m-j}$ which is located to the right of it.
\item\label{benn}
None of the pairs $\{\,\alpha^{-1}_{m-j},\alpha_{m-j}\,\}$ is nested in
a pair $\{\,\alpha^{-1}_{m-l},\alpha_{m-l}\,\}$ for any $l>j$.
\item\label{jo}
Any two copies of $\alpha^{-1}_{m-j}$ are separated in $W_f$ by a copy
of $\alpha^{-1}_{m-n}$ for some $n<j$.
\end{enumerate}
The last statement of course directly implies the Proposition.

The last letter of $W_f$ is $\alpha_{m-1}$ and this is the only positive
occurrence of this letter.
Hence \eqref{rend} and \eqref{benn} are obvious for $j=1$.
Also, since a second one wouldn't find a pair,
there can be at most one copy of $\alpha^{-1}_{m-1}$ in $W_f$, thus
\eqref{jo} is vacuously true for $j=1$.

Assume that the statements hold for all $j'=1,\ldots,j-1$. To prove
\eqref{rend}, assume that a certain copy $\alpha^{-1}$ of
$\alpha^{-1}_{m-j}$ forms the
pair $P$ with a copy $\alpha$ of $\alpha_{m-j}$, which is to the
left of it. Then, as
$\alpha$ can't be part of the final segment $S$, it belongs to some
$S_b$ with last letter $\alpha^{-1}_{i_b}$, where $i_b>m-j$, ie
$m-i_b<j$. This copy of $\alpha^{-1}_{i_b}=\alpha^{-1}_{m-(m-i_b)}$ is
to the left of $\alpha^{-1}$, hence it is part of a pair which is nested
inside $P$, which contradicts the hypothesis \eqref{benn} for $j'=m-i_b$.

To prove \eqref{benn}, suppose that a certain pair
$\{\,\alpha^{-1}_{m-j},\alpha_{m-j}\,\}$ is nested in $P=\{\,\alpha,
\alpha^{-1}\,\}$ where $\alpha$ has an index less than $m-j$. This copy
of $\alpha_{m-j}$ can't be part of $S$ because it has a symbol with lower
index (namely, $\alpha$ or $\alpha^{-1}$) to the right of it. It is clear
then that the first negative letter after $\alpha_{m-j}$ has index
$m-j'$, which is higher
than $m-j$ (ie\ $j'<j$), and it is still nested in $P$. Then so is the
pair containing it, which contradicts \eqref{benn} for $j'$.

Finally for \eqref{jo}, assume the contrary again, namely that there are
two copies of $\alpha^{-1}_{m-j}$ in $W_f$ that are not separated by any
higher index negative symbol. Then the pair of the first
$\alpha^{-1}_{m-j}$ can not lie between them either, because it would be
part of some $S_b$ and then $\alpha^{-1}_{i_b}$ would separate. Thus, we
have two pairs of the form $\{\,\alpha^{-1}_{m-j},\alpha_{m-j}\,\}$ nested
in one another. But then the positive symbol of the inner pair would be
followed by some $\alpha^{-1}_{i_b}$, which is part of a pair that is
nested in the outer pair, and that contradicts \eqref{benn} for
$j'=m-i_b$.
\end{proof}

\begin{Def}\label{def:B}
Let $1\le i,j\le q$. The element $B_{i,j}$ of the DGA of $\gamma_\beta$
is the sum of the following
products. For each path composed of parts of the strands of the braid $\beta$
that connects the left endpoint labeled $i$ to the
right endpoint labeled $j$ so that it only turns at quadrants facing up,
take the product of the crossings from left to right that it turns at.
\end{Def}

For example, $B_{i,j}$ contains the constant term $1$ if and only if
$j=\sigma(i)$. We will need the following polynomials of the $B_{i,j}$:

\begin{Def}\label{def:C}
Let $q\ge i>j\ge1$ and let
\[C_{i,j}=\sum_{\{\,i_1,\ldots,i_c,j\,\}\in D_i}
B_{i,i_1}B_{i_1,i_2}B_{i_2,i_3}\ldots B_{i_{c-1},i_c}B_{i_c,j}.\]
Similarly, for $1\le m\le q$, let
\begin{equation}\label{eq:Cii}
C_{m,m}=\sum_{\{\,i_1,\ldots,i_c\,\}\in D_m}
B_{m,i_1}B_{i_1,i_2}B_{i_2,i_3}\ldots B_{i_{c-1},i_c}B_{i_c,m}.
\end{equation}
Finally, for any $i$ and $j$, let
\begin{equation*}
M_{i,j}=\sum_{\{\,i_1,\ldots,i_c\,\}\in D_{\min\{i,j\}}}
B_{i,i_1}B_{i_1,i_2}B_{i_2,i_3}\ldots B_{i_{c-1},i_c}B_{i_c,j}.
\end{equation*}
\end{Def}

In particular, for each summand in $C_{i,j}$, $j$ is the last element of
the admissible sequence but it may occur elsewhere, too. For example,
$C_{1,1}=B_{1,1}$, $C_{2,1}=B_{2,1}$, and
$C_{3,1}=B_{3,1}+B_{3,2}B_{2,1}+B_{3,1}B_{1,2}B_{2,1}$.
Note also that $M_{1,j}=B_{1,j}$, $M_{i,1}=B_{i,1}$, $M_{m,m}=C_{m,m}$ and
$M_{i,i-1}=C_{i,i-1}$, whenever these expressions are defined.

\begin{tetel}\label{thm:relacio}
$\partial(a_m)=1+C_{m,m}$. Consequently, the index $0$ part
$H_0(L_\beta)$ of the contact
homology $H(L_\beta)$ has a presentation where the generators are the
crossings of $\beta$ and the relations are $C_{m,m}=1$ for
$m=1,\ldots,q$.
\end{tetel}

\begin{proof}
The $1$ in the formula comes from the $m^{\rm th}$ trivial disc. We claim that
the rest of the
contributions add up to $C_{m,m}$. From Proposition \ref{pro:blank} and the
paragraph preceding Definition \ref{def:sorozat}, we know that there can't
be any such terms other than the ones included in $C_{m,m}$. To see that all
such monomials actually arise from admissible discs, we just need to find
those discs. This can either be done by an inductive construction on $c$ (for
the inductive step, remove the smallest number from the sequence), or by
applying Blank's theorem.
\end{proof}

It is not clear whether $H(L_\beta)$ contains any non-zero higher index
part at all (except for the case of the unknot, when the single index $1$
crossing is a non-nullhomologous cycle). This is mainly why we only
work with $H_0(L_\beta)$ in this paper. This difficulty in handling
contact homology also underlines the importance of the augmentation that
we construct in the next section.

\section{Augmentations of braid closures}\label{sec:aug}

\begin{Def}
An \emph{augmentation} of the Lagrangian diagram $\gamma$ of a Legendrian
link $L$ is a subset $X$ of its crossings with the following properties.
\begin{itemize}
\item All elements of $X$ are proper crossings of $L$
(ie, intersections of different components are not allowed in $X$).
\item The index of each element of $X$ in any admissible grading is $0$
(in fact, this requirement implies the previous one).
\item For each generator $a$, the number of admissible discs with positive
corner $a$ and all negative corners in $X$ is even.
\end{itemize}
\end{Def}

The last requirement implies that the evaluation homomorphism (which is
defined on the link DGA, and which is also called an augmentation)
$\varepsilon_X\colon\mathscr A\to\Z_2$ that sends elements of $X$
to $1$ and other generators to $0$, gives rise to an algebra homomorphism
$(\varepsilon_X)_*\colon H(L)\to\Z_2$.

\begin{pelda}\label{ex:aug}
Consider the right-handed Legendrian trefoil knot diagram of
\figref{fig:trefli}. We claim that the set $\{\,b_3\,\}$ is an
augmentation. Indeed, the only two non-zero differentials (see
Example \ref{ex:trefli2} or the previous section) are
$\partial(a_1)=1+b_1+b_3+b_1b_2b_3$ and
$\partial(a_2)=b_2+b_2b_3+b_1b_2+b_2b_3b_1b_2$, and even these
vanish after mapping $b_1$ and $b_2$ to $0$ and $b_3$ to $1$.
\end{pelda}

\begin{megj}\label{rem:fuchs}
Let us return to the front diagram in the upper half of \figref{fig:pic}.
Such a diagram always has an \emph{admissible decomposition} (or \emph{ruling})
in the sense of \cite{chek2}: the only two values
of the Maslov potential (even though it's $\Z$--valued)
are $1$ on the upper strands and $0$ on the strands of the original braid
$\beta$, thus all crossings are Maslov, and
we may declare all of them switching (in the multi-component case,
consider only proper crossings). This gives
rise to a decomposition where the discs are nested in one another, so it's
admissible.

The existence of a ruling implies that $L_\beta$ is not a stabilized
link type for any positive braid $\beta$ \cite{chek2}. 
Also, by 
a theorem of Fuchs \cite{fuchs}, it implies that the diagram has an
augmentation\footnote{Fuchs and Ishkhanov \cite{masikirany} and
independently Sabloff \cite{josh} have also proven that the
existence of an augmentation implies the existence of a ruling.} 
(note this is
by no means unique). In his proof, Fuchs constructs an augmentation
of a diagram which is equivalent to the original but has a lot more
crossings. This means that the original diagram can also be
augmented: the pull-back of an augmentation by a DGA morphism, like
the ones listed in section \ref{sec:z2}, is again an augmentation.
\end{megj}

\begin{figure}[ht!]\tiny\anchor{fig:fonat}
\psfrag{1}{$1$}
\psfrag{2}{$2$}
\psfrag{3}{$3$}
\psfrag{4}{$4$}
\psfrag{5}{$5$}
\psfrag{6}{$6$}
\psfrag{7}{$7$}
\psfrag{8}{$8$}
\psfrag{a}{$1$}
\psfrag{b}{$2$}
\psfrag{c}{$3$}
\psfrag{d}{$4$}
\psfrag{e}{$5$}
\psfrag{f}{$7$}
\cl{\includegraphics[width=.9\hsize]{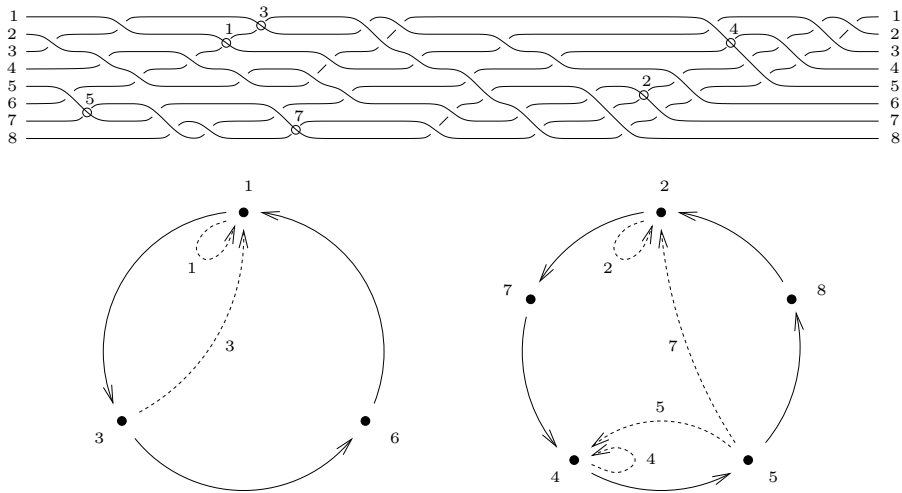}}
\caption{Constructing an augmentation for the closure of a
braid}\label{fig:fonat}
\end{figure}

In the case of a braid closure, it would be
impractical to pull back Fuchs' augmentation to the original diagram. Instead,
we'll start from scratch and construct an augmentation of the Legendrian
closure $\gamma_\beta$ of an arbitrary positive braid $\beta$.
When selecting crossings into $X$, it will suffice to work with the
braid itself, as illustrated in \figref{fig:fonat}, for the index 
$0$
proper crossings of $\gamma_\beta$ are all crossings of $\beta$.
(We will call these the \emph{proper crossings of the braid $\beta$}.)

First, we associate an oriented graph to an arbitrary
permutation $\sigma\in S_q$. Let $o$ be one of the cycles of
$\sigma$. Let us write the elements of $o$ around the perimeter of a
circle in the
cyclic order suggested by $\sigma$, directing an edge from $s$ to
$\sigma(s)$ for all $s$. (See \figref{fig:fonat} for an example. The
permutation in the diagram is the one underlying the braid).
If $p$ is a non-maximal element of
$o$, then follow the cycle in the forward direction starting from $p$
until it hits the first number in $o$ which is bigger then $p$. Let
$p_+$ be the number immediately before that. In particular, $p_+\le
p<\sigma(p_+)$. Do the same in the backward direction, resulting in the
number $p_-$ such that $p_-\le p<\sigma^{-1}(p_-)$.
(For example, in \figref{fig:fonat}, $5$ is an element of a 
$5$--cycle
$o$, with $5_+=5$ and $5_-=4$.) Then, connect $p_+$ to
$p_-$ by a directed chord of the circle and label the chord by $p$.
If $p_+=p_-=p$, then instead, we attach a loop edge (labeled $p$)
at $p$ to the graph. This will also be called a chord. Draw this loop
edge inside the circle, right next to the perimeter,
on the side of $p$ where the smaller of its two neighbors lies.

\begin{Def}
If $p$ is non-maximal in its cycle $o$, then the oriented loop
$\Gamma_p$ that starts from $p$, goes along the circle of $o$ to
$p_+$, then goes to $p_-$ on the chord labeled $p$, then follows
the circle again back to $p$ will be called the \emph{loop of
$p$}. If $\hat p$ is the largest number in $o$, then let the
\emph{loop of $\hat{p}$} be the loop $\Gamma_{\hat p}$ that
travels around the original circle once.
\end{Def}

\begin{lemma}\label{lem:diszj}
If $p<r$, then
\begin{enumerate}[\rm(a)]
\item\label{ize} $\Gamma_p$ doesn't contain $r$
\item\label{bize} the discs bounded by $\Gamma_p$ and $\Gamma_r$ are either
disjoint or the latter contains the former.
\end{enumerate}
In particular, the $|o|-1$ (oriented) chords obtained in the construction are
pairwise disjoint. They can't be parallel to the original edges and
they differ from each other as well.
\end{lemma}

\begin{proof}
Chords differ from edges because it's impossible that each end of an edge
be smaller than the other end (a chord and an edge can be parallel in the
non-oriented sense, as in \figref{fig:fonat}). The other statements
follow from \eqref{ize} and \eqref{bize}. Statement \eqref{ize} is obvious
from the construction, and so is \eqref{bize} if $r$ is maximal in $o$.

The directed
arc of the circle stretching from $p_-$ to $p_+$ only contains numbers less
than $r$. If this arc is disjoint from $\Gamma_r$, then of course so is
$\Gamma_p$. Otherwise, the whole arc must be contained in $\Gamma_r$. From
this, statement \eqref{bize} is clear, except if the chord labeled $p$ is a
loop edge and either $p=r_+$ or $p=r_-$. In the first case, $\sigma(p)>r$, but
$\sigma^{-1}(p)\le r$, so \eqref{bize} holds by the construction of the
loop edge. The other case is handled analogously.
\end{proof}

\begin{Def}
Let $\Gamma_\sigma$ be the disjoint union of
the graphs constructed above over all cycles $o$ of $\sigma\in S_q$.
This oriented planar graph, with vertices labeled
by the numbers $1,\ldots,q$, is called the \emph{augmented graph of the
permutation $\sigma$}.
\end{Def}

\begin{lemma}\label{lem:graf-szint}
For all $p\in\{\,1,\ldots,q\,\}$, the loop of $p$ is the
unique directed loop in $\Gamma_\sigma$
starting and ending at $p$ so that the sequence of the vertices (apart from
$p$) visited by it is in $D_p$.
\end{lemma}

\begin{proof}
It is enough to prove the statement
for a connected component associated to a cycle $o$.
The loop $\Gamma_p$
has the required property for all $p$, because all the vertices visited
by it are less than $p$ and no repetition occurs.

We have to rule out the existence of other loops.
Since $\sigma(p_+)$ and $\sigma^{-1}(p_-)$ are larger
than $p$, by the disjointness statement in Lemma \ref{lem:diszj},
all loops in question are trapped
in the disc bounded by $\Gamma_p$ (and this is obviously true
when $p$ is maximal in $o$). Then by statement \eqref{bize} of
Lemma \ref{lem:diszj}, apart from edges of $\Gamma_p$, they may
only contain chords labeled by numbers less than $p$.
Now, if the chord labeled by $s$ occurred in the loop and $s<p$ was the
smallest such number, then because $p$ is not on $\Gamma_s$,
the sequence of vertices on $\Gamma_s$ would appear as
a subsequence of the original vertex sequence of the loop. In that case,
$s$ would be
repeated in the sequence without the two occurrences separated by a
larger number, which is a contradiction.
\end{proof}

The following ``edge reversal lemma,'' which we will need in section
\ref{sec:moreaug}, can be proven very similarly.

\begin{lemma}\label{lem:masodik}
For any $1\le p<r\le q$ so that $\sigma(r)=p$, if $p$ is the second
largest vertex (after $r$) along $\Gamma_r$, then there is a unique
oriented path in $\Gamma_\sigma$ from $p$ to $r$ so that the sequence of
the intermediate vertices is in $D_p$. Otherwise, there is no such path.
\end{lemma}

\begin{Def}
Let $Y$ be a set of crossings of the positive braid $\beta$. By the
\emph{graph realized by
$Y$} we mean the oriented graph with vertices $1,\ldots,q$ so that a
directed edge
connects $i$ to $j$ if and only if $B_{i,j}\big|_Y=1$. Here, $B_{i,j}$ is
as in
Definition \ref{def:B} and by $B_{i,j}\big|_Y$ we mean the element of
$\Z_2$ obtained by substituting $1$ for elements of $Y$ and $0$ for other
generators in $B_{i,j}$.
\end{Def}

\begin{lemma}\label{lem:grafok}
Let $\sigma$ be the underlying permutation of $\beta$ and $Y$ a set of
proper crossings of $\beta$. If the graph realized by $Y$ agrees with the
augmented graph $\Gamma_\sigma$
of $\sigma$, then $Y$ is an augmentation of the Legendrian
closure of $\beta$.
\end{lemma}

\begin{proof}
Assume the two graphs do agree. Then by Lemma \ref{lem:graf-szint}, exactly
one of the summands of
$C_{m,m}$ (see equation \eqref{eq:Cii}) contributes $1$ to
the sum $C_{m,m}\big|_Y$,
namely the one that belongs to the sequence of vertices on
$\Gamma_m$. Therefore by Theorem \ref{thm:relacio}, $Y$ is an
augmentation.
\end{proof}

Next, based on $\Gamma_\sigma$, we construct a candidate $X$, and then we
will use Lemma \ref{lem:grafok} to prove that it's an augmentation.
Loosely speaking, the edges of $\Gamma_\sigma$
connecting $s$ to $\sigma(s)$ are always realized, even by the
empty set. To realize the chord labeled $p$, we'll select the crossing
$b_{p_+,p_-,1}$ into $X$. We can do this because it always exists:
$p_+\le
p<\sigma^{-1}(p_-)$ (for left labels) and $p_-\le p<\sigma(p_+)$ (for
right labels), therefore the strand connecting $p_+$ to $\sigma(p_+)$
always meets the strand connecting $\sigma^{-1}(p_-)$ to $p_-$. If there
were more than one points with the first two labels $p_+,p_-$, we could have
selected any of them\footnote{For instance, the crossing $b_3$ of Example
\ref{ex:aug} is denoted by $b_{1,1,2}$ in the general labeling system.};
we used the third label $1$ for concreteness
and for ease in the proof of Theorem \ref{thm:osztja}. In
\figref{fig:fonat}, we marked the selected crossings and labeled them
with the label of the chord that they realize.

\begin{all}\label{pro:aug}
The set
\[X=\{\:b_{p_+,p_-,1}\mid p\in\{\,1,\ldots,q\,\}\text{ is not a
maximal element of a cycle of }\sigma\:\}\]
is an augmentation of the Legendrian closure of the positive braid $\beta$
with underlying permutation $\sigma$.
\end{all}

In particular, for a pure braid $\beta$, the empty set is an augmentation.
In other words, the DGA of the Legendrian closure of a pure braid is
augmented, ie,\ the boundary of each generator is a polynomial without
a constant term. If $\beta$ is not pure, then $X\ne\varnothing$, hence
$\varepsilon_X\ne 0$, and it follows that $H(L_\beta)\ne 0$.

\begin{proof}
It is clear from the construction that all the crossings in $X$ are proper.
By Lemma \ref{lem:grafok}, it suffices to prove that the graph $G$
realized by $X$ is the graph $\Gamma_\sigma$. For
this, the chief claim is that no two points of $X$ are connected with a
part of a strand so that it arrives at both points from above. In other
words, the situation of \figref{fig:baj} can not arise: there is no
pair of numbers $p,r$ so that $\sigma(r_+)=p_-$. Indeed, then we'd have
$r<\sigma(r_+)=p_-\le p$ and similarly, $p<\sigma^{-1}(p_-)=r_+\le r$,
which would be a contradiction.

\begin{figure}[ht!]\small\anchor{fig:baj}
\psfraga <-2pt, 0pt> {p+}{$p_+$}
\psfraga <-2pt, 0pt> {p-}{$p_-$}
\psfraga <-2pt, 0pt> {r+}{$r_+$}
\psfraga <-2pt, 0pt> {r-}{$r_-$}
\psfraga <-3pt, 0pt> {p}{$p_+,p_-$}
\psfraga <-2pt, 0pt> {r}{$r_+,r_-$}
\cl{\includegraphics[width=.8\hsize]{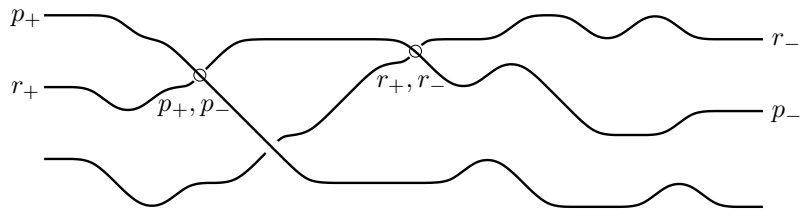}}
\caption{A situation that we have to rule out}\label{fig:baj}
\end{figure}

We know then that any path that is to contribute a non-zero
summand to a certain $B_{i,j}\big|_X$ can have at most one corner,
which of course has to be in $X$. Thus in the discussion before
the Proposition we exhausted all such contributions: paths with no
corners are responsible for the edges and paths with one corner
are responsible for the chords of $\Gamma_\sigma$. So indeed,
$G=\Gamma_\sigma$.
\end{proof}

\begin{figure}[ht!]\Large\anchor{fig:eukl}
\psfraga <-12pt, 3pt> {1}{$1$}
\psfraga <-12pt, 3pt> {2}{$2$}
\psfraga <-12pt, 3pt> {3}{$3$}
\psfraga <-12pt, 3pt> {4}{$4$}
\psfraga <-12pt, 3pt> {5}{$5$}
\psfraga <-12pt, 3pt> {6}{$6$}
\psfraga <-12pt, 3pt> {7}{$7$}
\psfraga <-12pt, 3pt> {8}{$8$}
\psfraga <-12pt, 3pt> {9}{$9$}
\psfraga <-12pt, 3pt> {10}{$10$}
\psfraga <-12pt, 3pt> {11}{$11$}
\psfraga <-12pt, 3pt> {12}{$12$}
\psfraga <-12pt, 3pt> {13}{$13$}
\psfraga <-12pt, 3pt> {14}{$14$}
\psfraga <-12pt, 3pt> {15}{$15$}
\psfraga <-12pt, 3pt> {16}{$16$}
\psfraga <-12pt, 3pt> {17}{$17$}
\psfraga <-12pt, 3pt> {18}{$18$}
\psfraga <-12pt, 3pt> {19}{$19$}
\psfraga <-12pt, 3pt> {20}{$20$}
\psfraga <-12pt, 3pt> {21}{$21$}
\psfraga <-12pt, 3pt> {22}{$22$}
\psfraga <-12pt, 3pt> {23}{$23$}
\psfraga <-12pt, 3pt> {24}{$24$}
\psfraga <-12pt, 3pt> {25}{$25$}
\psfraga <-12pt, 3pt> {26}{$26$}
\psfraga <0pt, 3pt> {1'}{$1$}
\psfraga <0pt, 3pt> {2'}{$2$}
\psfraga <0pt, 3pt> {3'}{$3$}
\psfraga <0pt, 3pt> {4'}{$4$}
\psfraga <0pt, 3pt> {5'}{$5$}
\psfraga <0pt, 3pt> {6'}{$6$}
\psfraga <0pt, 3pt> {7'}{$7$}
\psfraga <0pt, 3pt> {8'}{$8$}
\psfraga <0pt, 3pt> {9'}{$9$}
\psfraga <0pt, 3pt> {10'}{$10$}
\psfraga <0pt, 3pt> {11'}{$11$}
\psfraga <0pt, 3pt> {12'}{$12$}
\psfraga <0pt, 3pt> {13'}{$13$}
\psfraga <0pt, 3pt> {14'}{$14$}
\psfraga <0pt, 3pt> {15'}{$15$}
\psfraga <0pt, 3pt> {16'}{$16$}
\psfraga <0pt, 3pt> {17'}{$17$}
\psfraga <0pt, 3pt> {18'}{$18$}
\psfraga <0pt, 3pt> {19'}{$19$}
\psfraga <0pt, 3pt> {20'}{$20$}
\psfraga <0pt, 3pt> {21'}{$21$}
\psfraga <0pt, 3pt> {22'}{$22$}
\psfraga <0pt, 3pt> {23'}{$23$}
\psfraga <0pt, 3pt> {24'}{$24$}
\psfraga <0pt, 3pt> {25'}{$25$}
\psfraga <0pt, 3pt> {26'}{$26$}
\vspace{-1cm}\newline
\cl{\includegraphics[angle=45,width=.9\linewidth]{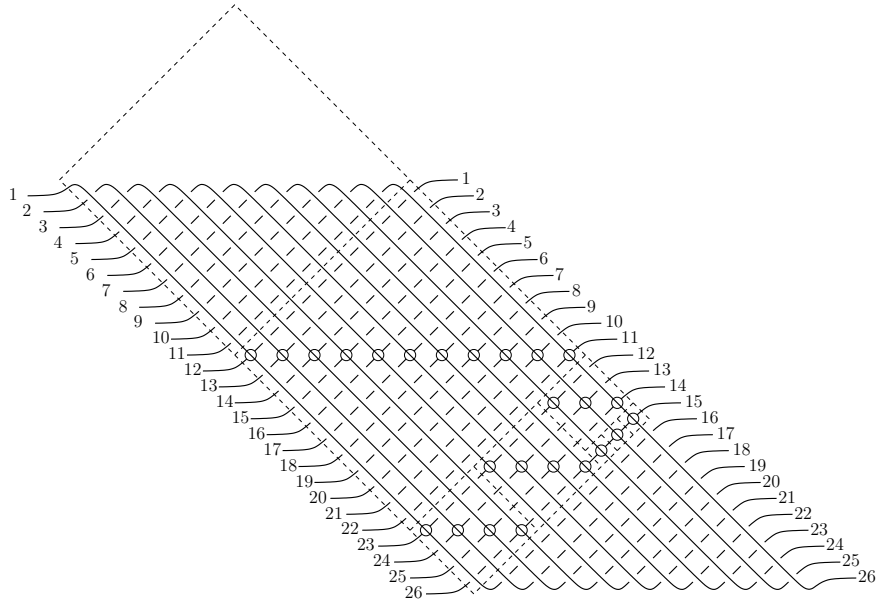}}
\vspace{-3.5cm}\caption{The augmentation of a $(p,q)$ torus link
implements the Euclidean algorithm; on the diagram,
$p=11$ and $q=26$.}\label{fig:eukl}
\end{figure}

\begin{megj}
Let $\beta$ be the standard positive braid whose Legendrian closure is a
positive $(p,q)$ torus link (the braid used to produce
\figref{fig:pic} is an example with $p=5$ and $q=4$).
If we apply our construction to it, we find
an interesting connection of the resulting augmentation to the Euclidean
algorithm. We mention this here
without proof; we will only need a small part of the statement which is 
hidden in the proof of Proposition
\ref{pro:seq}.

Let us denote the quotients and residues in the Euclidean algorithm
(with input $p$ and $q$) as follows:
\begin{align*}
p&=k_{-1}q+r_0 & (0&\le r_0<q)\\
q&=k_0r_0+r_1 & (0&\le r_1<r_0)\\
r_0&=k_1r_1+r_2 & (0&\le r_2<r_1)\\
r_1&=k_2r_2+r_3 & (0&\le r_3<r_2)\\
&\vdots&&\\
r_{l-3}&=k_{l-2}r_{l-2}+r_{l-1} & (0&\le r_{l-1}<r_{l-2})\\
r_{l-2}&=k_{l-1}r_{l-1}+r_{l} & (0&\le r_{l}<r_{l-1})\\
r_{l-1}&=k_lr_l+0.&&
\end{align*}
(Of course, $r_l=\gcd\{\,p,q\,\}$.) The points of the augmentation $X$ are
arranged
in blocks of the following sizes: $(k_0-1)$ blocks of size $r_0$; $k_1$
blocks of size $r_1$; $k_2$ blocks of size $r_2$ and so on until the last
$k_l$ blocks of size $r_l$. If we draw the diagram of the braid as in 
\figref{fig:eukl}, every block can be viewed as the diagonal of a square, and
the squares can in turn be seen to be placed inside a $p\times q$ rectangular
box so that they realize a `graphic implementation' of the
Euclidean algorithm.
\end{megj}

\section{A loop of positive links}\label{sec:loop}

\begin{figure}[ht!]\small\anchor{fig:97front}
\cl{\includegraphics[width=.9\hsize]{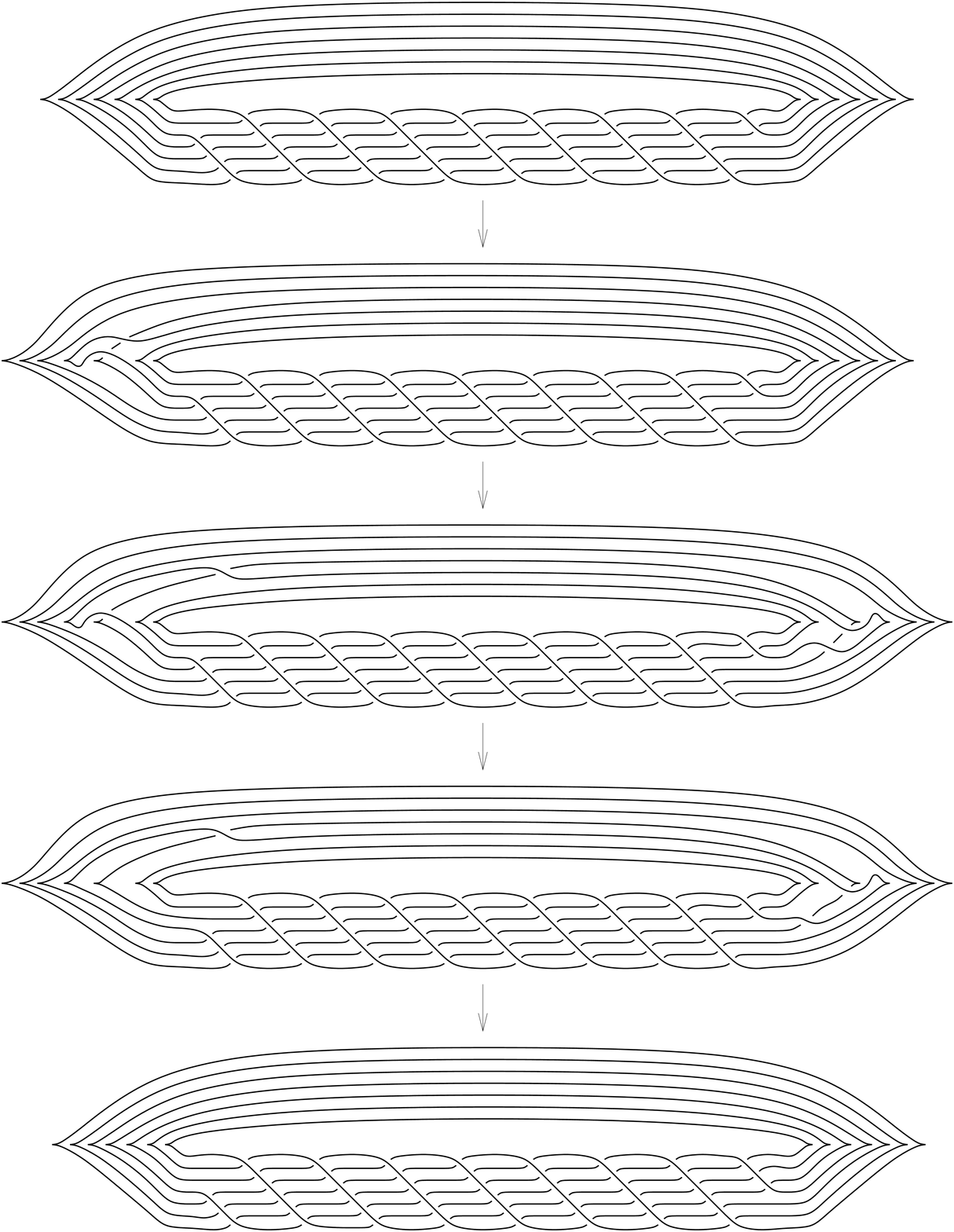}}
\caption{A path that corresponds to conjugating a braid, viewed in
the front projection. The moves are similar to those in 
Figure \ref{fig:hurok}.}\label{fig:97front}
\end{figure}

\begin{figure}[ht!]\small\anchor{fig:97lag}
\psfrag{a}{$a_{\text{temp}}$}
\psfrag{al}{$a_{\kern-0.7pt m}$}
\psfrag{all}{$a_{m\kern-0.7pt{+}\kern-0.7pt1}$}
\psfrag{b}{$b_{m,1}$}
\psfrag{c}{$c_{m,1}$}
\psfraga <-3pt, 1pt> {e}{$b_{m+1,1}$}
\psfraga <-3pt, 0pt> {d}{$\ddots$}
\psfraga <-1pt, 1pt> {u}{$b_{q-1,1}$}
\psfrag{U}{$U_{m+2}$}
\cl{\includegraphics[width=.9\hsize]{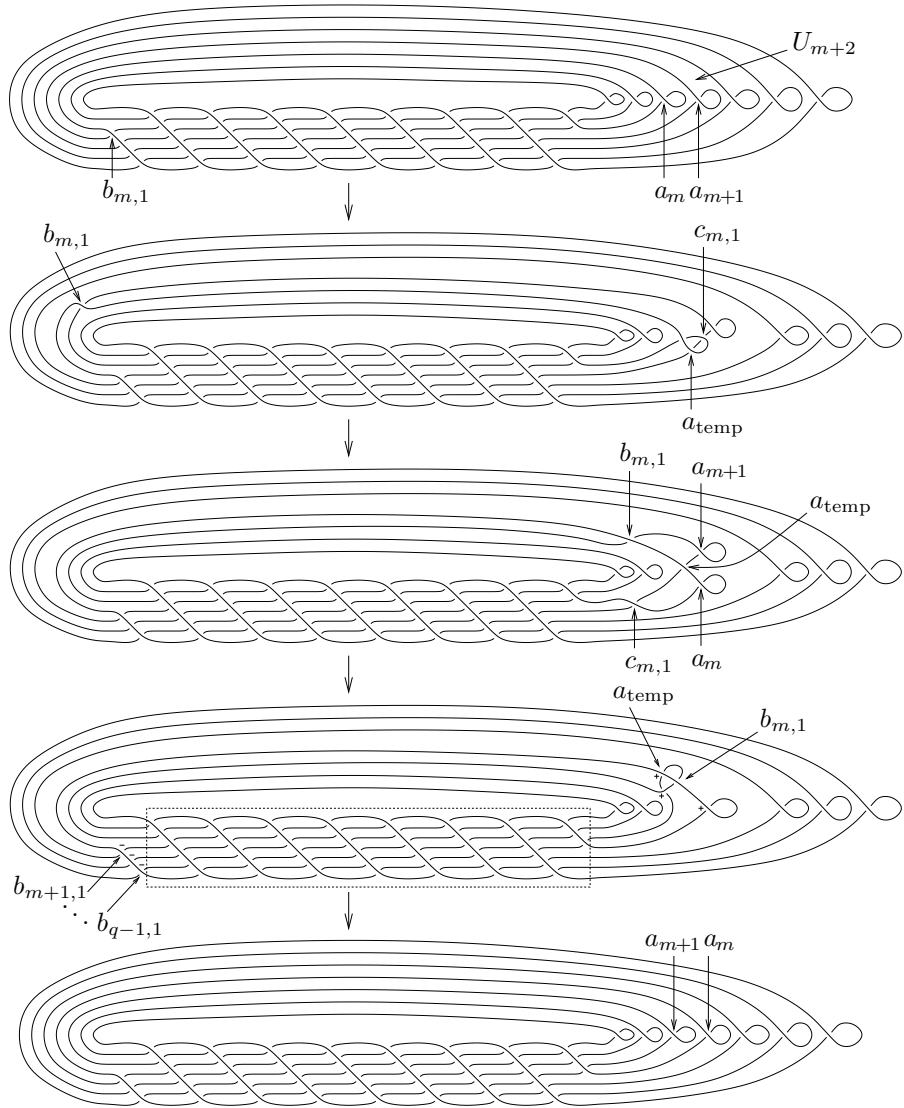}}
\caption{A path
that corresponds to conjugating a braid, viewed in the Lagrangian
projection. The individual Reidemeister moves \emph{do not} correspond
to those in Figure \ref{fig:97front}.}\label{fig:97lag}
\end{figure}

Let $L_\beta$ be the Legendrian closure of the positive braid
$\beta$. Then there exists a natural closed loop in the connected
component $\mathscr L_\beta$ of the space of Legendrian links that
contains $L_\beta$, as follows. Let us write
$\beta=\lambda_1\ldots\lambda_w$ as a product of the braid group
generators. In \figref{fig:97front}, we show through an
example how $L_\beta$ can be changed into the Legendrian closure
of the conjugate braid that results from moving the first factor
$\lambda_1$ to the end of the word: if $\lambda_1$ is a half-twist
of the $m^{\rm th}$ and $(m+1)^{\rm st}$ strands of the braid, then one
interchanges the $m^{\rm th}$ and $(m+1)^{\rm st}$ strands \emph{above} the
braid. In \figref{fig:97lag}, the same path $\Phi_{\lambda_1}$
is shown, but in the Lagrangian projection. Note that the index
$1$ crossings $a_m$ and $a_{m+1}$ trade places. (The notation used
on the diagram for the index $0$ crossings is the one that we will
introduce below for the special case of torus links.) The
Lagrangian diagrams of the endpoints are clearly obtained by
resolution of the corresponding fronts. However, we will not prove
that the paths themselves agree, too (up to homotopy). Instead, we
will content ourselves with checking (using Theorem
\ref{thm:linprog}) that the four Reidemeister moves in
\figref{fig:97lag} are consistent, and thereafter use the Lagrangian
construction as our definition of $\Phi_{\lambda_1}$. (There are
no such consistency issues with Reidemeister moves of fronts, but
we need the Lagrangian diagrams to compute holonomies.)

\begin{tetel}\label{thm:consistent}
The sequence of Reidemeister moves in \figref{fig:97lag}, defining
$\Phi_{\lambda_1}$, is consistent.
\end{tetel}

\begin{proof}
This is a generalization of Example \ref{ex:trefli4}. Isotope the
diagram just like in the proof of Lemma \ref{lem:kicsi}. Then by
choosing a small enough $\varepsilon$, the height $h(a_{m+1})$
will dominate the expression whose positivity is needed (by
Theorem \ref{thm:linprog}) in order for the first move (which is a
Reidemeister II move) to be consistent. To carry out the
III$_{\text{b}}$ move that follows, we need to isotope the second
diagram from the top so that what remains from the region
$U_{m+2}$ after the first move has larger area than the vanishing
triangle. This can be achieved by the same trick (this time,
moving away the `outer' $q-m-1$ strands). Next, the newborn
triangle needs to be blown up so that it has larger area than the
triangle which is due to vanish in the second III$_{\text{b}}$
move (which is the third move altogether). For this, the same
trick with the bulges still works: apply it to the $q-m-1$ outer
strands and the one that crosses itself at $a_m$. Finally, the
exact same argument guarantees that the fourth move, of type
II$^{-1}$, is consistent too.
\end{proof}

Let us use the formulas of section \ref{sec:z2} to compute the
holonomy $\mu_1$ of $\Phi_{\lambda_1}$ (cf Example
\ref{ex:trefli3}). More precisely, we will compute the action of
$\mu_1$ on the index $0$ crossings which generate the index $0$
contact homology $H_0(L_\beta)$. In the first move, two new
crossings appear; using Lemma \ref{lem:fokok} it is easy to show that 
their indices are $0$ and $1$.
Let us denote the one with index $0$
by\footnote{In other words, let it inherit the labels of the
crossing that is being moved to the other end of the braid. This
is what we've done in \figref{fig:97lag} too, except that
there, a different notation is used for index $0$ crossings.}
$c_{m,\sigma(m+1),1}$ and the one with index $1$ by
$a_{\text{temp}}$. The old index $0$ crossings are not affected by
this move (ie, the holonomy maps them trivially). This is true
by Proposition \ref{pro:easy}: for index reasons, the boundary of
any index $0$ crossing is $0$. It is easy to see that the
following two triangle moves don't affect the old index $0$
crossings, either. 

In the fourth, Reidemeister II$^{-1}$ move
however the crossing $b_{m,\sigma(m+1),1}$ (together with
$a_{\text{temp}}$) vanishes, and its image in the holonomy becomes
the polynomial $M'_{m+1,m}=C'_{m+1,m}$. We used primed symbols to
remind us that these are to be computed with respect to the
conjugated braid.
The fact that for each
admissible disc that turns at the positive
quadrant at $a_{\text{temp}}$ that faces away from $b_{m,\sigma(m+1),1}$,
the intermediate sequence of labels
(of other index $1$ generators through which the boundary of the admissible
disk passes) doesn't contain
any number larger than $m$ can be shown just like in the proof of Lemma
\ref{lem:kicsi}. The fact that the sequence is admissible follows by an
argument very similar to the proof of Proposition \ref{pro:blank}.
Finally, the fact that each such admissible sequence does contribute the said
terms can be established as in the proof of Theorem \ref{thm:relacio}.
We have proven:

\begin{all}\label{pro:mualt}
The holonomy $\mu_1$ of 
$\Phi_{\lambda_1}$ 
maps each crossing of the braid
trivially except the first one from the left, which is mapped to the
polynomial $M'_{m+1,m}=C'_{m+1,m}$ (if the first crossing is between the
$m^{\rm th}$ and $(m+1)^{\rm st}$ strands). This expression is to be computed as in
Definition \ref{def:C}, with respect to the conjugate braid
$\lambda_2\lambda_3\ldots\lambda_w\lambda_1$.
\end{all}

Note also that at the end stage of $\Phi_{\lambda_1}$ each crossing of the 
(conjugate) braid has index $0$, thus in the multi-component case our 
choice of preferred grading is justified.

Now, it is clear that the concatenation
$\Omega_\beta=\Phi_{\lambda_1}\ldots\Phi_{\lambda_w}$ is a closed loop in
$\mathscr L_\beta$ and by definition, its monodromy is the composition of
holonomies $\mu_w\circ\ldots\circ\mu_1$, followed by a re-labeling to restore
the original labels. Namely, each symbol $c$ labeling an index $0$ crossing
needs to be changed back to $b$;
in fact, $a_i$ would have to be changed to $a_{\sigma(i)}$, but because
we only concern ourselves with $H_0(L_\beta)$, this can be ignored.

\section[mu_0 divides p+q]{\boldmath $|\mu_0|$ divides $p+q$}\label{sec:moreBC}

The last two sections of the paper contain the proof of Theorem
\ref{thm:monorendje}. The argument works for any $p,q$ except when
$q$ divides $p$ (the case of a pure braid), or $p$ divides $q$.
The reason why we don't claim Theorems \ref{thm:monorendje} and
\ref{thm:trefli} for multi-component torus links is that we only
proved Theorems \ref{thm:moves} and \ref{thm:inv} (and thus
Theorem \ref{thm:hurkok}) for knots. (However, the extension of
those proofs should only be a matter of changing the formalism to
that of link contact homology.)

Let us revisit the loop $\Omega_{p,q}$ of Legendrian $(p,q)$ torus
knots defined in the introduction. As the braid $\beta$ is now
composed of $p$ periods, the general loop $\Omega_\beta$ described in 
section
\ref{sec:loop} is the $p$-fold concatenation of another, and the
latter is easy to identify as homotopic\footnote{Because we'll
omit the rigorous justification of this fact, the reader may treat
the new description of $\Omega_{p,q}$ as the definition.} to
$\Omega_{p,q}$. In particular, by Theorem \ref{thm:inv}, the
`full' monodromy takes the form $\mu^p$, where $\mu$ is the
monodromy of $\Omega_{p,q}$. From now on, we will concentrate on
this map $\mu$, and especially on its restriction $\mu_0$ to the
index $0$ part of the contact homology of the standard torus link
diagram $\gamma$ (shown in \figref{fig:pic}), representing the
base point $L$.

We will adjust our notation to this special situation. The crossings of
$\beta$ will be indexed with two integers (as opposed to three), namely
$b_{m,n}$ ($m=1,\ldots,q-1$, $n=1,\ldots,p$) will denote the $m^{\rm th}$ crossing
counted from the top in the $n^{\rm th}$ period of the braid.
Note that in the definition of $\mu$, after a full period of the
braid has been moved from the left end to the right end, a re-labeling
takes place, too: the second label of each crossing in the other $(p-1)$
periods is reduced by $1$, and the labels $c_{m,1}$ in the now last period
are changed to $b_{m,p}$.

\begin{all}\label{pro:muspec}
The monodromy $\mu$ of the loop $\Omega_{p,q}$ of Legendrian torus links
acts on the index $0$ generators as follows:
\[\mu(b_{m,n})=\left\{\begin{array}{ll}
b_{m,n-1} & \text{if }2\le n\le p\\
C_{q,m} & \text{if }n=1
\end{array}\right..\]
\end{all}

\begin{proof}
The claim is clear for those crossings not in the first period: by
Proposition \ref{pro:mualt}, they are only affected, and in the described way,
by the re-labeling. The rest of the statement will be proven by induction on
$q-m$. When this value is $1$, ie\ $m=q-1$, this is just the statement of
Proposition \ref{pro:mualt} (the conjugate braid in this case is the
original $\beta$ again, and the re-labeling changes $C'_{q-1+1,q-1}$ into
$C_{q,q-1}$). Assume the statement holds for $b_{q-1,1},\ldots,
b_{m+1,1}$.
Right after the conjugation that removes it from the left end of the braid,
the image of $b_{m,1}$ is $M''_{m+1,m}$, computed with respect to the braid
after this conjugation. This can be re-written (by grouping terms with
respect to the first factor in the product) as
\[M''_{m+1,m}=b_{m+1,1}M'_{m+1,m}+b_{m+2,1}M'_{m+2,m}+\ldots+
b_{q-1,1}M'_{q-1,m}+M'_{q,m},\]
where the terms labeled $M'$ on the right are to be computed in the braid
indicated
by the box in \figref{fig:97lag}. Note however that by Definition
\ref{def:C}, the same expressions are obtained
if we use the whole braid $\beta$
(before re-labeling). If we apply the holonomies of the remaining $q-1-m$
conjugations and the re-labeling to this expression, we get (by the inductive
hypothesis)
\[\mu(b_{m,1})=C_{q,m+1}M_{m+1,m}+C_{q,m+2}M_{m+2,m}+\ldots+
C_{q,q-1}M_{q-1,m}+M_{q,m}=C_{q,m}.\]
The last equality is true because the middle expression is exactly what
results if we group terms in $C_{q,m}$ with respect to the last label
in the admissible sequence which is more than $m$.
\end{proof}

\begin{all}\label{pro:muBC}
In the contact homology ring $H(L)$,
we have:
\begin{equation}\label{eq:muB}\mu(B_{i,j})=\left\{\begin{array}{ll}
B_{i-1,j-1}+B_{i-1,q}b_{j-1,p}&\text{if }i,j\ge2\\
B_{i-1,q}&\text{if }i\ge2\text{ and }j=1\\
b_{j-1,p}&\text{if }i=1\text{ and }j\ge2
\end{array}\right.,\end{equation}
\begin{equation}\label{eq:muC}\mu(C_{i,j})=\left\{\begin{array}{ll}
C_{i-1,j-1}&\text{if }j\ge2\\
M_{i-1,q}&\text{if }j=1
\end{array}\right.,\end{equation}
and
\begin{equation}\label{eq:muM}
\mu(M_{i,j})=M_{i-1,j-1},\text{ whenever }i,j\ge2.\end{equation}
\end{all}

We omitted $B_{1,1}$ and $C_{i,i}$ because they (and hence their images)
are equal to $1$ in the contact homology. Recall also that $M_{i,1}=B_{i,1}$
and $M_{1,j}=B_{1,j}$.

\begin{proof}
When $i\ge2$, none of the terms in $B_{i,j}$ contains any of
$b_{1,1},\ldots,b_{q-1,1}$, so they only need to be re-labeled. This means
that all crossings that the path through the braid which generated the term
turned at, are shifted to the
left by a unit. This operation changes the entry point from the one labeled
$i$ to the one labeled $i-1$.
If $j=1$, then the shifted path can be completed
by the overcrossing strand of the last period of the braid, which shows that
the re-labeled expression is a summand in $B_{i-1,q}$. Moreover, all such
summands
are obtained in this way exactly once. When $j\ge2$, the re-labeling results
in a summand of $B_{i-1,j-1}$, but not all such are obtained: we miss
contributions from paths that turn at the crossing (the last one on the
strand with right endpoint $j-1$) $b_{j-1,p}$. Hence the correction term in
the top row of \eqref{eq:muB} (note that the paths turning at $b_{j-1,p}$
are exactly those that would otherwise have arrived at $q$).

When $i=1$ and $j\ge2$, we have
\begin{equation*}\begin{split}
\mu(B_{1,j})&=
\mu(b_{1,1}B_{2,j}+b_{2,1}B_{3,j}+\ldots+b_{q-1,1}B_{q,j}+R)\\
&=C_{q,1}(B_{1,j-1}+B_{1,q}b_{j-1,p})+C_{q,2}(B_{2,j-1}+B_{2,q}b_{j-1,p})+
\ldots\\
&+C_{q,q-1}(B_{q-1,j-1}+B_{q-1,q}b_{j-1,p})+\mu(R)\\
&=C_{q,1}B_{1,j-1}+C_{q,2}B_{2,j-1}+\ldots+C_{q,q-1}B_{q-1,j-1}\\
&+(C_{q,1}B_{1,q}+C_{q,2}B_{2,q}+\ldots+C_{q,q-1}B_{q-1,q})b_{j-1,p}\\
&+B_{q,j-1}+B_{q,q}b_{j-1,p}\\
&=C_{q,j-1}+C_{q,j-1}C_{j-1,j-1}+C_{q,q}b_{j-1,p}
=C_{q,j-1}+C_{q,j-1}+b_{j-1,p}\\
&=b_{j-1,p}.
\end{split}\end{equation*}
Here, $R$ is the sum of the contributions to $B_{1,j}$ that don't contain
crossings of the first period. These terms only have to be re-labeled and that
can be done just like in the argument above for $B_{i,j}$ when $i\ge2$.
In the sum of sums
\[C_{q,1}B_{1,j-1}+C_{q,2}B_{2,j-1}+\ldots+C_{q,q-1}B_{q-1,j-1}+B_{q,j-1},\]
we re-grouped the terms; those with an admissible sequence of labels
formed $C_{q,j-1}$, and the rest, where $j-1$ was repeated `illegally,'
formed $C_{q,j-1}C_{j-1,j-1}$.

Note that by the now proven \eqref{eq:muB}, for all $i,j\ge2$,
$\mu(B_{i,j}+B_{i,1}B_{1,j})=B_{i-1,j-1}$. Therefore, when $i>j\ge2$,
\begin{equation*}\begin{split}
\mu(C_{i,j})&=
\mu\left(\sum_{\{\,i_1,\ldots,i_c,j\,\}\in D_i}
B_{i,i_1}B_{i_1,i_2}B_{i_2,i_3}\ldots B_{i_{c-1},i_c}B_{i_c,j}\right)\\
&=\mu\left(\sum_{\text{\scriptsize$\begin{gathered}\{\,j_1,\ldots,j_d,j\,\}\in D_i\\
j_1,\ldots,j_d\ge2\end{gathered}$}}
\begin{aligned}[t]
(B_{i,j_1}+B_{i,1}B_{1,j_1})(B_{j_1,j_2}+&B_{j_1,1}B_{1,j_2})\ldots\\
(B_{j_{d-1},j_d}+B_{j_{d-1},1}B_{1,j_d})&(B_{j_d,j}+B_{j_d,1}B_{1,j})
\end{aligned}\right)\\
&=\sum_{\text{\scriptsize$\begin{gathered}\{\,j_1,\ldots,j_d,j\,\}\in D_i\\
j_1,\ldots,j_d\ge2\end{gathered}$}}B_{i-1,j_1-1}B_{j_1-1,j_2-1}\ldots
B_{j_{d-1}-1,j_d-1}B_{j_d-1,j-1}\\
&=C_{i-1,j-1}.
\end{split}\end{equation*}
As a consequence of this and \eqref{eq:muB}, for all $i\ge2$,
\[\mu(C_{i,1})=\mu\left(B_{i,1}+\sum_{j=2}^{i-1}C_{i,j}B_{j,1}\right)
=B_{i-1,q}+\sum_{j=2}^{i-1}C_{i-1,j-1}B_{j-1,q}=M_{i-1,q}.\]
Finally, when $2\le i<j$,
\begin{equation*}\begin{split}
\mu(M_{i,j})&=\mu\left(B_{i,j}+\sum_{k=1}^{i-1} C_{i,k}B_{k,j}\right)\\
&=B_{i-1,j-1}+B_{i-1,q}b_{j-1,p}+M_{i-1,q}b_{j-1,p}\\
&+\sum_{k=2}^{i-1} C_{i-1,k-1}(B_{k-1,j-1}+B_{k-1,q}b_{j-1,p})\\
&=M_{i-1,q}b_{j-1,p}+\left(B_{i-1,j-1}
+\sum_{k=2}^{i-1}C_{i-1,k-1}B_{k-1,j-1}\right)\\
&+\left(B_{i-1,q}+\sum_{k=2}^{i-1}C_{i-1,k-1}B_{k-1,q}\right)b_{j-1,p}\\
&=M_{i-1,q}b_{j-1,p}+M_{i-1,j-1}+M_{i-1,q}b_{j-1,p}\\
&=M_{i-1,j-1},
\end{split}\end{equation*}
and when $i>j$, the argument is very similar to the one we gave for $C_{i,j}$.
\end{proof}

\begin{tetel}\label{thm:szamol}
The order of the (restricted)
monodromy $\mu_0=\mu\big|_{H_0(L)}$ of the loop $\Omega_{p,q}$ of
Legendrian $(p,q)$ torus knots divides $p+q$.
\end{tetel}

\begin{proof}
This is now a straightforward computation, generalizing the first
paragraph of the proof of Proposition \ref{pro:5}. Consider the
generator $b_{m,p}$ ($m=1,\ldots,q-1$). By Proposition
\ref{pro:muspec}, the first $p$ iterations of $\mu$ act on it as
follows:
\[\mu(b_{m,p})=b_{m,p-1},\text{ }\mu^2(b_{m,p})=b_{m,p-2},\dots,
\mu^{p-1}(b_{m,p})=b_{m,1},\text{ }\mu^p(b_{m,p})=C_{q,m}.\]
Then by \eqref{eq:muC} of Proposition \ref{pro:muBC},
the next $m$ iterations are as follows:
\[\mu^{p+1}(b_{m,p})=C_{q-1,m-1},\ldots,\mu^{p+m-1}(b_{m,p})=C_{q-m+1,1},
\text{ }\mu^{p+m}(b_{m,p})=M_{q-m,q}.\]
Now by \eqref{eq:muM}, the next $q-m-1$ iterations are
\[\mu^{p+m+1}(b_{m,p})=M_{q-m-1,q-1},\ldots,\mu^{p+q-1}(b_{m,p})=M_{1,m+1}.\]
Finally, because $M_{1,m+1}=B_{1,m+1}$, \eqref{eq:muB} yields
\[\mu^{p+q}(b_{m,p})=b_{m,p}.\]
Because $b_{m,n}$ is on the orbit of $b_{m,p}$ for all $n=1,\ldots,p$, we see
that $\mu^{p+q}$ is the identity on all of the degree $0$ generators.
\end{proof}

\section[(p+q) divides |mu_0|]{\boldmath $(p+q)$ divides $|\mu_0|$}\label{sec:moreaug}

\begin{tetel}\label{thm:osztja}
The number $p+q$ divides the order of the monodromy $\mu$.
\end{tetel}

In the proof of Theorem \ref{thm:szamol}, we described 
explicitly 
the $(p+q)$-element orbit
of each index $0$ generator $b_{m,n}$
(altogether $q-1$ orbits). Recall that all of those orbits
contain a \emph{$b$--sequence} $b_{m,p}, b_{m,p-1},\ldots,b_{m,1}$ of 
length $p$, a
\emph{$C$--sequence} $C_{q,m}, C_{q-1,m-1},\ldots,C_{q-m+1,1}$ of length 
$m$ and an
\emph{$M$--sequence} $M_{q-m,q}, M_{q-m-1,q-1},\ldots,M_{1,m+1}$ of 
length $q-m$.
(These expressions are cycles in the chain complex
$\mathscr A$, but we really mean the homology classes represented by
them.) We will choose
one of these orbits
and show that it doesn't have a period shorter than $p+q$. For 
this, we will evaluate the augmentation
$\varepsilon=\varepsilon_X$ (see Proposition \ref{pro:aug}) on elements of
the orbit, and prove that the resulting sequence of $0$'s and $1$'s,
which we will call the \emph{$0$-$1$--sequence of the orbit}, has
no such 
shorter period. In fact, we claim the following:

\begin{all}\label{pro:01}
If $q>p$ but $p\nmid q$, the $0$-$1$--sequence $S$ of the orbit of $b_{p,p}$
consists of $p$ consecutive $0$'s and $q$ consecutive $1$'s.
If $q<p$ but $q\nmid p$, then the same holds for the orbit of
$b_{[p\pmod q],p}$.
\end{all}

In the latter case, we will denote the value $1\le[p\pmod q]\le q-1$ by
$r_0$. For the rest of the section, either this value $r_0$ or $p$, as the
case may be, should be substituted for $m$ in the formulas for the $b$--,
$C$--, and $M$--sequences.
Note that if $q\mid p$, the said orbit doesn't even
exist (if $q>p$ and $p\mid q$, then its $0$-$1$--sequence consists only
of $1$'s).

Recall that $X$ was constructed so that the graph realized by $X$ was
the augmented graph of the underlying permutation $\sigma$ of the braid.
Hence, this oriented graph
$\Gamma_\sigma$ has adjacency matrix $\left[\varepsilon_X(B_{i,j})\right]$
and
therefore it contains all the information we need to evaluate
the algebra homomorphism $\varepsilon=\varepsilon_X$ on the polynomial
expressions of the $C$-- and $M$--sequences. We will only need to refer to the
actual braid in the case of the $b$--sequence. In our situation,
\[\sigma(i)=[(i-p)\pmod q],\text{ }i=1,\ldots,q.\]
When $q<p$, we could equivalently write $\sigma(i)=[(i-r_0)\pmod q]$.
This explains why our choice of orbit in Proposition \ref{pro:01} is
reasonable: both in the $C$--sequence and in the $M$--sequence the two
lower indices are always the endpoints of an edge of $\Gamma_\sigma$, but
\emph{they are listed in the reverse order}. So when we evaluate
$\varepsilon$ on these polynomials, what we need to examine
is whether
the given edge of the graph can be ``reversed,'' ie,\ if it is part
of an oriented
loop (and how many loops) with an admissible sequence of vertices.

We will re-state and prove Proposition \ref{pro:01} in a more detailed
version.

\begin{all}\label{pro:seq}
The orbit specified in Proposition \ref{pro:01} contributes
$0$'s and $1$'s to the sequence $S$ as follows.
\begin{enumerate}[\rm(1)]

\item\label{egy} If $2p\le q$ but $p\nmid q$ (hence in fact $2p<q$),
then we have
\begin{multline*}
\overbrace{b_{p,p},\ldots,b_{p,1}}^{p\text{ copies of }1};
\overbrace{C_{q,p},\ldots,C_{q-p+1,1}}^{p\text{ copies of }1};\\
\overbrace{M_{q-p,q},\ldots,M_{q-2p+1,q-p+1}}^{p\text{ copies of }0},
\overbrace{M_{q-2p,q-p},\ldots,M_{1,p+1}}^{q-2p\text{ copies of }1}.
\end{multline*}

\item\label{ket} If $2p>q$ but $q>p$, then the sequence is
\[
\overbrace{b_{p,p},\ldots,b_{p,q-p+1}}^{2p-q\text{ copies of }0},
\overbrace{b_{p,q-p},\ldots,b_{p,1}}^{q-p\text{ copies of }1};
\overbrace{C_{q,p},\ldots,C_{q-p+1,1}}^{p\text{ copies of }1};
\overbrace{M_{q-p,q},\ldots,M_{1,p+1}}^{q-p\text{ copies of }0}.
\]

\item\label{ha} If $q<p$, $q\nmid p$, and $2r_0<q$, then we get
\begin{multline*}
\overbrace{b_{r_0,p},\ldots,b_{r_0,r_0+1}}^{p-r_0\text{ copies of }0},
\overbrace{b_{r_0,r_0},\ldots,b_{r_0,1}}^{r_0\text{ copies of }1};
\overbrace{C_{q,r_0},\ldots,C_{q-r_0+1,1}}^{r_0\text{ copies of }1};\\
\overbrace{M_{q-r_0,q},\ldots,M_{r_0+1,2r_0+1}}
^{q-2r_0\text{ copies of }1},
\overbrace{M_{r_0,2r_0},\ldots,M_{1,r_0+1}}^{r_0\text{ copies of }0}.
\end{multline*}

\item\label{negy} Finally, if $q<p$, $q\nmid p$, and $2r_0\ge q$, then
we have
\begin{multline*}
\overbrace{b_{r_0,p},\ldots,b_{r_0,q-r_0+1}}
^{p-q+r_0\text{ copies of }0},
\overbrace{b_{r_0,q-r_0},\ldots,b_{r_0,1}}^{q-r_0
\text{ copies of }1};\\
\overbrace{C_{q,r_0},\ldots,C_{q-r_0+1,1}}^{r_0\text{ copies of }1};
\overbrace{M_{q-r_0,q},\ldots,M_{1,r_0+1}}^{q-r_0\text{ copies of }0}.
\end{multline*}
\end{enumerate}
\end{all}

\begin{proof}
\textbf{\boldmath $b$--sequence}\qua In case \eqref{egy}, all numbers
$1\le j\le p$ are so that
$j<[(j+p)\pmod q]=j+p$ and $j<[(j-p)\pmod q]$. This means that $j_+=j_-=j$,
ie, that there is a loop edge attached to $j$ in $\Gamma_\sigma$. It is
easy to check that in the construction of $X$, the crossing that realizes
this loop edge is exactly $b_{p,j}$.

In case \eqref{ket}, we similarly find loop edges but only attached to
the
numbers $1,\ldots,q-p$. These are realized by the crossings
$b_{p,1},\ldots,b_{p,q-p}$. For the crossings $b_{p,q-p+1},\ldots,b_{p,p}$,
we find that their second labels in the system that we used to label
crossings
of general braids in section \ref{sec:BC} are also $q-p+1,\ldots,p$. These
numbers can not be the endpoints of a chord because the numbers preceding
them
in the permutation (namely, $1,\ldots,2p-q$) are smaller than them.
Therefore these crossings are indeed not selected into $X$.

In cases \eqref{ha} and \eqref{negy}, first note that the crossings
$b_{r_0,q+1},\ldots,b_{r_0,p}$ have third labels greater than $1$ in
the labeling system of section \ref{sec:BC}, so they never get
selected into $X$. Neither do $b_{r_0,r_0+1},\ldots,b_{r_0,q}$,
because their `old' first labels are $r_0+1,\ldots,q$ and these are
taken to the smaller values $1,\ldots,q-r_0$ by $\sigma$ (ie,
they'll never be the startpoint of a chord). After this, the rest of
the
$b$--sequence can be sorted out just like in the first two cases.

\textbf{\boldmath $C$--sequence}\qua Here the claim is that it always contributes
only $1$'s to $S$. This is because $i>\sigma^{-1}(i)=[(i+p)\pmod q]$
implies
$\varepsilon(C_{i,[(i+p)\pmod q]})=1$. Indeed, since it is preceded in
the permutation by a smaller number, $i$ can't
be the endpoint of a chord, only of the single edge coming from
$\sigma^{-1}(i)$.
So there is a unique term in $C_{i,[(i+p)\pmod q]}$ that
contributes $1$ to $\varepsilon(C_{i,[(i+p)\pmod q]})$,
namely the one which, when multiplied by $B_{[(i+p)\pmod q],i}$ on the
right,
produces the term corresponding to the unique loop described in
Lemma \ref{lem:graf-szint}.

\textbf{\boldmath $M$--sequence}\qua All the four claims in this case follow from
Lemma \ref{lem:masodik}. Note in particular that if we allowed
$p\mid q$ in case \eqref{egy}, then $j-p$ would be the second largest
vertex on $\Gamma_j$ for all $j=p+1,\ldots,q$, and therefore the first
$p$ elements of the $M$--sequence wouldn't be mapped to $0$ by
$\varepsilon$. However if $p\nmid q$, then the first number $x$
after the
sequence $j,j-p,j-2p,\ldots$ `wraps around' the circle $\Z_q$ is
different from $j$. If $x$ is also smaller than $j$ (and this will be
the case exactly when $j=q,q-1,\ldots,q-[q\pmod p]+1$), then it falls
between $j$ and $j-p$, so $j-p$
is not second largest in $\Gamma_j$. If $x$ is larger than $j$, then
it is recognized as what we called $\sigma(j_+)$ in section
\ref{sec:aug}. But then there exists a chord in $\Gamma_\sigma$,
starting from $j_+$, and ending at the element $j_-$ of
$\Gamma_j$. This $j_-$ is by construction such that $[(j_-+p)\pmod
q]>j$. But because $p<q/2$, if we assume that $j>q-p$,
then this is
only possible if $j_-+p>j$, ie if $j_->j-p$. This again means that
in these cases, $j-p$ is not second largest on $\Gamma_j$. This proves
the claim about the first part of the $M$--sequence. Finally, if
$j\le q-p$ (which implies $j_-=j$), then $\Gamma_j$ only visits the
positive elements of the arithmetic
progression $j,j-p,j-2p,\ldots$  (ie, there is no wrapping around)
and $j-p$ is obviously second largest among these. Cases \eqref{ket},
\eqref{ha} and \eqref{negy} can be handled similarly.
\end{proof}

Note that this was a generalization of the second paragraph in the
proof of Proposition \ref{pro:5}; there, we used the augmentation
$X=\{\,b_3\,\}$.

\end{document}

%% file: gtspec.tex
\def\ifplaintex{\expandafter\ifx\csname documentclass\endcsname\relax}

\def\ifplaintex{\expandafter\ifx\csname documentclass\endcsname\relax}

\ifplaintex 
\else
\fi

\expandafter\ifx\csname epsfbox\endcsname\relax\input epsf\fi

\def\gt{{\mathsurround=0pt\it $\cal G\mskip-2mu$eometry \&\ 
$\cal T\!\!$opology}}        

\def\gtp{{\mathsurround=0pt\it $\cal G\mskip-2mu$eometry \&\ 
$\cal T\!\!$opology $\cal P\!$ublications}}  

\def\lognumber#1{\def\thelognumber{#1}}
\def\volumenumber#1{\def\thevolumenumber{#1}}
\def\papernumber#1{\def\thepapernumber{#1}}
\def\volumeyear#1{\def\thevolumeyear{#1}}

\def\pagenumbers#1#2{\def\startpage{#1}\def\finishpage{#2}}
\def\published#1{\def\publishdate{#1}}
\def\proposed#1{\def\theproposer{#1}}
\def\seconded#1{\def\theseconders{#1}}
\def\received#1{\def\receiveddate{#1}}
\def\revised#1{\def\reviseddate{#1}}
\def\accepted#1{\def\accepteddate{#1}}

\def\coverauthors#1{\def\thecoverauthors{#1}}
\def\asciiauthors#1{\def\theasciiauthors{#1}}

\long\def\asciiabstract#1{\long\def\theasciiabstract{#1}}

\let\\\par\let\thelognumber\relax
\let\thevolumenumber\relax\let\thepapernumber\relax
\let\thevolumeyear\relax\let\thesamplenumber\relax\let\startpage\relax
\let\finishpage\relax\let\publishdate\relax\let\receiveddate\relax
\let\reviseddate\relax\let\accepteddate\relax\let\theasciititle\relax
\let\theasciiauthors\relax
\let\theasciiabstract\relax
\let\theasciiemail\relax\let\theshortauthors\relax\let\theshorttitle\relax
\let\thecoverauthors\relax

\long\def\maketitlep{   

\count0=\startpage

\gt\hfill      
\hbox to 77pt{\vbox to 0pt{\vglue -15pt\epsfbox{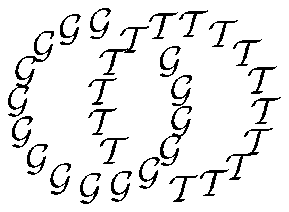}\vss}\hss}
\break
{\small\ifx\thesamplenumber\relax 
Volume \else Sample
\fi\thevolumenumber\ (\thevolumeyear)
\startpage--\finishpage\nl
Published: \publishdate}
\vglue 0.5truein plus 0.4fil minus 0.1truein

{\parskip=0pt\leftskip 0pt plus 1fil\def\\{\par\smallskip}{\ifplaintex\large
\else\Large\fi\bf\thetitle}\par\medskip}   

\vglue 0pt plus 0.1fil 

{\parskip=0pt\leftskip 0pt plus 1fil\def\\{\par}{\sc\theauthors}
\par\medskip}

\vglue 0pt plus 0.1fil 

{\small\parskip=0pt\let\newline\\
{\leftskip 0pt plus 1fil\def\\{\par}{\sl\theaddress}\par}
\expandafter\ifx\theemail\relax    
\relax\else\vglue 5pt plus 0.02fil minus 2pt\def\\{\stdspace{\rm 
and}\stdspace} 
\cl{Email:\stdspace\tt\theemail}\fi
\ifx\theurl\relax                  
\relax\else\vglue 5pt plus 0.02fil minus 2pt\def\\{\stdspace{\rm 
and}\stdspace}
\cl{URL:\stdspace\tt\theurl}\fi\par}

\vglue 7pt plus 0.3fil minus 3pt

{\bf Abstract}
\vglue 5pt plus 0.1fil minus 2pt

\theabstract

\vglue 7pt plus 0.3fil minus 3pt

{\bf AMS Classification numbers}\quad Primary:\quad \theprimaryclass

Secondary:\quad \thesecondaryclass

\vglue 5pt plus 0.3fil minus 2pt

{\bf Keywords:}\quad \thekeywords

\vglue 10pt plus 0.5fil minus 5pt

{\small  Proposed: \theproposer\hfill Received: \receiveddate\nl
Seconded: \theseconders\hfill 
\ifx\reviseddate\relax                         
Accepted: \accepteddate                        
\else
Revised: \reviseddate                          
\fi}
\eject
}       

\font\phead=cmsl9 scaled 950
\font=cmsl9 scaled 1050
\font\pnum=cmbx10 scaled 913
\font\lnum=cmbx10 
\font\pfoot=cmsl9 scaled 950
\font=cmsl9 scaled 1050
\ifplaintex
\headline{\vbox to 0pt{\vskip -4.5mm\line{\small\phead\ifnum
\count0=\startpage ISSN 1364-0380 (on line)
1465-3060 (printed) \hfill {\pnum\folio}\else\ifodd\count0\def\\{ }%
\ifx\theshorttitle\relax\thetitle\else\theshorttitle\fi\hfill{\pnum\folio}
\else\def\\{ and }{\pnum\folio}\hfill\ifx\theshortauthors\relax\theauthors
\else\theshortauthors\fi\fi\fi}\vss}}
\footline{\vbox to 0pt{\vglue 0mm\line{\small\pfoot\ifnum\count0=\startpage
\copyright\ \gtp\hfill\else
\gt, Volume \thevolumenumber\ (\thevolumeyear)\hfill\fi}\vss
}}
\else
\makeatletter
\def\@oddhead{{\small\lhead\ifnum\count0=\startpage ISSN 1364-0380 (on line)
1465-3060 (printed) \hfill {\lnum\number\count0}\else\ifodd\count0
\def\\{ }\ifx\theshorttitle\relax \thetitle \else\theshorttitle\fi\hfill
{\lnum\number\count0}\else\def\\{ and }{\lnum\number\count0}
\hfill\ifx\theshortauthors\relax 
\theauthors\else\theshortauthors\fi\fi\fi}}\def\@evenhead{\@oddhead}
\def\@oddfoot{\small\lfoot\ifnum\count0=\startpage\copyright\ \gtp\hfill\else
\gt, Volume \thevolumenumber\ (\thevolumeyear)\hfill\fi}
\def\@evenfoot{\@oddfoot}
\makeatother
\fi

\newwrite\gtoutfile
\long\gdef\makeheadfile{  
{\def\\{, }\def\s{ }
\immediate\openout\gtoutfile head.xxx
\immediate\write\gtoutfile{Proxy-for: \ifx\theasciiauthors\relax
\theauthors\else\theasciiauthors\fi\s<\ifx\theasciiemail\relax\theemail\else\theasciiemail\fi>}
\immediate\write\gtoutfile{\noexpand\\}
\immediate\write\gtoutfile{Authors: \ifx\theasciiauthors\relax
\theauthors\else\theasciiauthors\fi}
{\def\\{ }\immediate\write\gtoutfile{Title: \ifx\theasciititle\relax
\thetitle\else\theasciititle\fi}}
\immediate\write\gtoutfile{Subj-class: GT or SG or MG etc}
\immediate\write\gtoutfile{MSC-class: \theprimaryclass\ifx\thesecondaryclass\relax\else, \thesecondaryclass\fi}
\immediate\write\gtoutfile{Journal-ref: Geom. Topol. \thevolumenumber
(\thevolumeyear) \startpage-\finishpage}
\immediate\write\gtoutfile{Comments: Published by Geometry and Topology at}
\immediate\write\gtoutfile{\s\s http://www.maths.warwick.ac.uk/gt/GTVol\thevolumenumber/paper\thepapernumber.abs.html}
\immediate\write\gtoutfile{\noexpand\\}
\immediate\write\gtoutfile{}
\ifx\theasciiabstract\relax
\immediate\write\gtoutfile{\theabstract}\else
\immediate\write\gtoutfile{\theasciiabstract}\fi
\immediate\write\gtoutfile{}
\immediate\write\gtoutfile{\noexpand\\}
\immediate\write\gtoutfile{}
\immediate\closeout\gtoutfile}}  

\def\maketitlepage{\maketitlep\makeheadfile}
\let\maketitle\maketitlepage